%% file: main__4_.tex
\documentclass[final,leqno,letterpaper]{etna}

\usepackage{amsmath, amssymb} 
\usepackage{graphicx}                 
\usepackage{bm}                       
\usepackage{rotating}
\usepackage{hyperref}                 
\usepackage{algorithm}
\usepackage{algpseudocode}
\usepackage{booktabs}
\usepackage{tabularx}
\usepackage{xcolor}
\setbibdata{1}{xx}{46}{2017} 
\usepackage{subcaption}
\usepackage{caption}
\usepackage{longtable}
\usepackage{makecell}

\hypersetup{%
    pdftitle={Comparative analysis of numerical methods for X-ray 3D reconstruction from sparse views},
    pdfauthor={Melitine N Bardi, Konstantinos K Delibasis},
    pdfkeywords={inverse reconstruction problem, large sparse linear systems, x-ray 3D reconstruction, sparse views }
    }

\title{Comparative analysis of numerical methods for X-ray 3D reconstruction from sparse views
}

\author{Melitine N. Bardi \footnotemark[1]
        \and Konstantinos K. Delibasis\footnotemark[1]}

\shorttitle{Numerical Methods for Sparse-View X-ray Reconstruction} 
\shortauthor{M.~BARDI AND K.~DELIBASIS}

\definecolor{projectionBG}{RGB}{224,239,255}
\definecolor{arnoldiBG}{RGB}{231,245,233}
\definecolor{leastSquaresBG}{RGB}{255,243,224}
\definecolor{updateBG}{RGB}{243,229,245}

\definecolor{projectionText}{RGB}{25,80,140}
\definecolor{arnoldiText}{RGB}{35,105,55}
\definecolor{leastSquaresText}{RGB}{150,90,20}
\definecolor{updateText}{RGB}{105,50,130}

\algrenewcommand\algorithmiccomment[1]{%
  \hfill\textcolor{gray!75!black}{$\triangleright$ #1}}

\algrenewcommand\algorithmicrequire{%
  \textcolor{blue!60!black}{\textbf{Input:}}}
\algrenewcommand\algorithmicensure{%
  \textcolor{blue!60!black}{\textbf{Output:}}}

\newcommand{\AlgSection}[3]{%
  \Statex
  \Statex \colorbox{#1}{%
    \parbox{\dimexpr\linewidth-2\fboxsep\relax}{%
      \textcolor{#2}{\textbf{#3}}}}
  \Statex
}

\begin{document}

\maketitle

\renewcommand{\thefootnote}{\fnsymbol{footnote}}

\footnotetext[1]{University of Thessaly, Lamia}

\begin{abstract}
Reconstruction of three-dimensional vascular structures from a limited number of X-ray projections constitutes a severely underdetermined and ill-conditioned inverse problem. This work presents a unified matrix-free computational framework for investigating iterative numerical methods for sparse-view cone-beam tomography. The forward and backprojection operators are implemented using ray–voxel intersection lengths obtained through the Amanatides–Woo traversal algorithm, avoiding explicit construction of the large sparse system matrix. We conduct a systematic comparison of algebraic reconstruction techniques (ART, SART, and SIRT), gradient descent (GD), and Krylov subspace methods, including Conjugate Gradient, LSQR, and GMRES. Particular attention is given to gradient scaling through scalar step sizes, Jacobi preconditioning, and diagonal Hessian majorization, as well as the incorporation of L1, Tikhonov, and total variation regularization.

Numerical experiments are performed on five, three-dimensional (3D) coronary vessel models, each reconstructed from ten simulated X-ray projections under noise and scattering conditions. Real X-ray projections from a static phantom with thin vessels are also used to extend the validation. Reconstruction performance is assessed in terms of convergence, computational cost, volumetric overlap, contrast-to-noise ratio, vessel diameter and centerline errors, and projection consistency at an unseen viewing angle. The results demonstrate that total-variation-regularized gradient descent, combined with diagonal majorization, provides the most favorable overall reconstruction performance under the investigated computational budgets. In contrast, the evaluated Krylov subspace implementations exhibit limited geometric recovery despite their greater computational cost. These findings highlight the importance of appropriate gradient scaling and regularization for stable and computationally efficient sparse-view tomographic reconstruction.
\end{abstract}

\begin{keywords}
inverse reconstruction problem, large sparse linear systems, x-ray 3D reconstruction, sparse views 
\end{keywords}

\begin{AMS}
65R32, 65F10, 92C55
\end{AMS}

\section{Introduction} X-ray fluoroscopy by C-arm is a key imaging modality used in many medical specialties, such as interventional cardiology, providing real-time two-dimensional (2D) visualization of three-dimensional (3D) coronary blood vessels during diagnostic and therapeutic procedures.

As the C-arm is placed in several orientations to complete the diagnostic or therapeutic process, a small number of views can be acquired with known geometry. The posibility of reconstructing the vessel tree in three dimensions (3D), as a quick alternative to cone-beam computer tomography -CBCT, has attracted research interest. 
Table~\ref{tab:literature_review} summarizes the main research contributions to this iterative tomographic reconstruction from sparse views, encompassing algebraic methods, total-variation regularization, first-order optimization, and Krylov subspace techniques. The studies are compared in terms of their numerical formulations, sparse-view acquisition settings, experimental data, and evaluation criteria. Although considerable attention has been devoted to numerical convergence, computational efficiency, and reconstruction quality, very few studies have focused on the preservation of the three-dimensional geometry of thin, branching vascular structures under severe angular undersampling. Also, most of the studies focus on a single reconstruction technique against a base-line, thus a systematic framework that evaluates a wide range of reconstruction methods in terms of sufficiency and extensive error metrics is lacking. 
\include{table_bibliography}
The present work is systematic numerical comparison of algebraic, gradient-based, and Krylov methods, with different regularization and preconditioning schemes, evaluated on both volumetric reconstruction accuracy and vessel-specific geometric errors, such as vessel center-line positioning error and diameter error. The evaluation was performed on patient-derived coronary phantoms and experimentally acquired phantom data.


In previous work \cite{bardi2026forward}, we introduced a simulation framework for the forward problem of coronary X‑ray imaging, by explicitly modeling the cone‑beam C‑arm geometry and its full kinematic degrees of freedom. The simulator performs voxel‑based X‑ray attenuation along each ray using the Beer–Lambert law, and additionally, accounts for the influences of measurement noise and Compton scattering within the model. This framework yields synthetic projections that closely resemble clinical coronary fluoroscopy under clinically relevant viewing angles.

Utilizing this established forward model 
in this work, we implement the 3D Reconstruction from sparse views by defining the vascular structure as a voxel grid that is iteratively updated, so that its projections, obtained by the same C-arm geometry and ray-traversal pipeline, match the observed angiographic images. The discrepancy between observed and estimated volume linear absorption coefficient values becomes the objective function, subject to different modes of regularization. 
The contributions of this work are summarized below.
\begin{itemize}

    \item \textbf{Unified matrix-free reconstruction framework} for sparse-view X-ray tomography, using the same forward and matched backprojection operators without explicitly constructing the very large system matrix $\mathbf{A}$, or the Hessian of the objective function.

    \item \textbf{Systematic comparison of fundamentally different reconstruction paradigms} under identical geometry, projection data, initialization, and stopping criteria, including ART, SART, SIRT, regularized gradient descent (GD), GD with Nesterov's acceleration, conjugate gradient (CG), and Krylov subspace methods GMRES, and LSQR.

    \item \textbf{Investigation of regularized gradient-based reconstruction}, including several regularization terms and a voxel-wise majorization strategy based on $\mathbf{M}= \operatorname{diag} \left( (\mathbf{A}^{T}\mathbf{A})\mathbf{1} \right),$
    aiming at improved scaling and convergence.

    \item \textbf{Quantitative assessment of reconstruction performance in the sparse-view regime}, including reconstruction accuracy such as volume Intersection-over-union, reconstructed vessel CNR (Contrast-to-Noise ratio), absolute value of the residual, and the structural similarity index (SSIM) between a view unseen during reconstruction, generated from the ground truth and the reconstructed volumes. Convergence behavior and computational cost are also evaluated.


    \item \textbf{Vessel-specific geometric validation on patient-derived voxelized volumes and projections from real phantom containing vessels:} 
    The proposed evaluation framework quantifies reconstruction accuracy directly in terms of vascular geometry through the local vessel diameter error and the symmetric skeleton positional error, assessing both vessel caliber preservation and 3-D centerline localization beyond conventional voxel-wise image-error measures.

\end{itemize}

\section{Methodology}

The goal of tomographic reconstruction is to estimate the spatial distribution of linear attenuation coefficients within an object from a set of measured projection data (sinograms).

The inverse problem involves the reconstruction of a 3D volume from a set of projections acquired at various angles. We define the following notation:

\begin{itemize}
    \item Pixels and Voxels: Let ${p} = (u, v)$ represent an image pixel and ${v} = (i, j, k)$ represent a volume voxel.
    \item Ray Geometry: Each pixel $p$ corresponds to a ray intersecting a set of voxels: $\mathcal{R}_p = \{ v \in \text{Ray}(p) \}$.
    \item Physical Properties: $\boldsymbol{\mu}=\{\mu_{{v}}\}$ is the linear attenuation coefficient ($\text{cm}^{-1}$) of voxel $v$, and $A_{pv}$ is the length of the intersection between ray $p=1,\ldots,N_p$ and voxel $v=1,\ldots,N_v$. the number of voxels is equal to the spatial dimensions of $\boldsymbol{\mu}, \text{ }N_v=V_xV_yV_z.$
    \item Measurements: $I$ is the detected image intensity, while $I_0$ is is the incident (unattenuated) intensity.
\end{itemize}

\subsection{Forward Projection Model}

\subsubsection{Continuous Physical Model}
The physical process of X-ray attenuation along a line path is governed by the Beer-Lambert Law:
\begin{equation}
    I = I_0 \exp \left( -\int_{L_p} \mu(\mathbf{x}) \, ds \right)
\end{equation}
where 

\begin{itemize}
    \item $I_0$ denotes the incident X-ray intensity,
    \item $I$ is the measured intensity after traversing the object, 
    \item $\mu(\mathbf{x}) \ge 0$ is the spatial linear attenuation coefficient at spatial coordinate $\mathbf{x}$,
    \item $L_p$ represents the line trajectory of ray $p$.
\end{itemize}

Taking the negative logarithm yields the line integral $y_p$, representing the total integrated attenuation along ray $p$ in the continuous and the deiscrete case:
\begin{equation}
    y_p = -\ln \left( \frac{I}{I_0} \right) = \int_{Ray(p)} \mu(\mathbf{x}) \, ds
\end{equation}

\begin{equation}
Y_p = -\ln\left(\frac{I_p}{I_0}\right) 
= \sum_{v \in \mathcal{R}_p} \mu_{v} A_{pv} + \epsilon_p,
\quad \epsilon_p \sim \mathcal{N}(0,\sigma^2).
\end{equation}

\begin{figure}[htbp]
    \centering
    \includegraphics[width=0.8\textwidth]{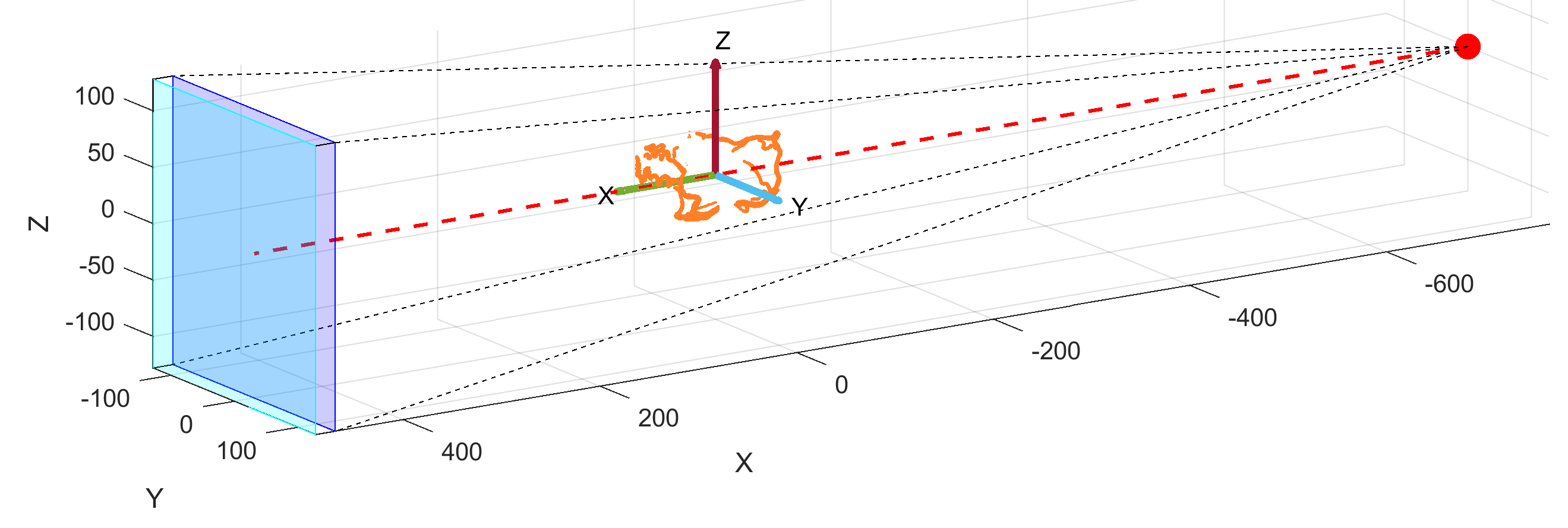}
    \caption{Simulator geometry showing the X-ray source and detector planes relative to the heart vessels and patient coordinate axes ($X, Y, Z$).}
    \label{fig:3d_model}
\end{figure}

\subsubsection{Discrete System Matrix}
To solve this problem computationally, the continuous volume domain is discretized into $N_v$ voxels, each assigned an attenuation value $\mu_v$. Approximating the integral as a finite sum over all intersecting voxels yields:
\begin{equation}
    Y_p = \sum_{v=1}^{N_v} A_{pv} \mu_v, \quad p = 1, \dots, N_p
\end{equation}
where $A_{pv} \ge 0$ represents the intersection segment length of ray $p$ through voxel $v$. In vector-matrix notation, the forward projection model for the full system is:
\begin{equation}
    \mathbf{Y} = \mathbf{A} \boldsymbol{\mu}
    \label{eq:matrix_System}
\end{equation}
where:
\begin{itemize}
    \item $\mathbf{Y} \in \mathbb{R}^{N_p}$ is the vector of measured projection data (sinogram),
    \item $\boldsymbol{\mu} \in \mathbb{R}^{N_v}$ is the discrete volume attenuation vector to be reconstructed,
    \item $\mathbf{A} \in \mathbb{R}^{N_p \times N_v}$ is the system matrix (forward projector), where $A_{pv}$ is the ray-voxel intersection length.
\end{itemize}


\subsubsection {Noise and scattering models}

Let $\boldsymbol{\mu}_0$ denote the ideal binary vessel volume, with vessel
voxels equal to one and background voxels equal to zero. Before forward
projection, independent voxel-wise perturbations are added according to
\begin{equation}
    \boldsymbol{\mu}=\boldsymbol{\mu}_{0} + \boldsymbol{\eta},
    \qquad
    \eta_v\sim\mathcal{U}(-\nu,\nu),
    \label{eq_vol_noise}
\end{equation}
yielding the perturbed volume $\boldsymbol{\mu}$ to increase the realism of the volume. 

Assuming $N_0$ incident photons per detector pixel, the primary photon
intensity is modeled according to the Beer-Lambert law as
\begin{equation}
    \mathbf{I}^{P}=N_0\exp(-\mathbf{Y}).
\end{equation}
Compton scatter is approximated by spatial convolution of the primary
intensity with a normalized Gaussian kernel $K$,
\begin{equation}
    \mathbf{I}^{S}=f_S\left(\mathbf{I}^{P} * K\right),
\end{equation}
where $f_S=0.05$. In the present implementation, $K$ is an $81\times81$
Gaussian kernel with standard deviation $\sigma_K=13$ pixels. The total
expected detector intensity is therefore
\begin{equation}
    \mathbf{I}=\mathbf{I}^{P}+\mathbf{I}^{S}.
\end{equation}

Photon-counting noise is subsequently introduced using a Gaussian approximation of the Poisson distribution,
\begin{equation}
    \mathbf{I}^{N}
    =
    max \left (\mathbf{I}+\sqrt{\mathbf{I}}\,\xi ,0 \right),
    \qquad
    \xi\sim\mathcal{N}(0,1),
\end{equation}
with independent realizations for the detector pixels, with $\mathcal{N}(m,\sigma)$ is the normal distribution with mean $m$ and standard deviation $\sigma.$ 
Finally, the noisy projection data used for reconstruction are obtained through the inverse Beer-Lambert transformation,
\begin{equation}
    \mathbf{Y}^{N}
    =
    -\ln\left(\frac{\mathbf{I}^{N}}{N_0}\right).
\end{equation}
For simplicity, herein-forth we drop the $N$ superscript from $\mathbf{Y}$, assuming that the projections $\mathbf{Y}$ are always noisy.

\subsection{Inverse Problem formulation and efficient optimization methods}

\subsubsection{Mathematical Definition and Adjoint Duality}
While the forward projector maps volume space to projection space ($\mathbf{A}: \mathbb{R}^{N_v} \to \mathbb{R}^{N_p}$), the backprojection operator maps projection data back into volume space ($\mathbf{A}^T: \mathbb{R}^{N_p} \to \mathbb{R}^{N_v}$):
\begin{equation}
    \boldsymbol{\mu} = \mathbf{A}^T \mathbf{Y}
    \label{eq:backprojection}
\end{equation}
For any voxel $v$, backprojection accumulates the weighted attenuation measurements of all rays $p$ that traverse it:
\begin{equation}
    \mu_v = \sum_{p=1}^{N_p} A_{pv} Y_p
\end{equation}
Under the standard Euclidean inner product, the forward and backprojection operators satisfy the Hilbert space adjoint identity:
\begin{equation}
    \langle \mathbf{A}\boldsymbol{\mu}, \mathbf{Y} \rangle_{\mathbb{R}^{N_p}} = \langle \boldsymbol{\mu}, \mathbf{A}^T \mathbf{Y} \rangle_{\mathbb{R}^{N_v}}, \quad \forall \mu \in \mathbb{R}^{N_v}, \; Y \in \mathbb{R}^{N_p}
\end{equation}

\subsection{Memory Optimization of System Matrix}
\label{sec:comp_matrices}
The primary goal is to recover the spatial attenuation vector $\boldsymbol{\mu}$ from $N_\pi$ measured projections $\mathbf{Y}$ (views) by solving system \eqref{eq:matrix_System}, which typically requires matrix-vector multiplications such as $\mathbf{A}x,\mathbf{A}^Tx$ and $\mathbf{A}^TAx$. Because $\mathbf{A}$ is typically ill-conditioned, sparse, and subject to measurement noise, direct matrix construction inversion is impossible. Furthermore, a typical size for the system matrix $A$ could be $\mathbb{R}^{(N_{\pi} N_p) \times N_v}=\mathbb{R}^{(10 \cdot 512^2) \times 256^3}$, therefore, explicit construction is infeasible. Instead, we resorted to an implicit representation of $\mathbf{A}$ that requires substantially less computer memory.

\subsubsection{Auxiliary lookup Tables}
\label{sec:aux_tab}
Instead of explicitly constructing $A$, the implementation of the voxel traversal algorithm has been modified to return the following Tables:
\begin{itemize}
    \item \textbf{PV} ($K_1 \times N_p$): for each ray $p$, it stores a list of the voxel indices that are traversed.
    
    \item \textbf{PL} ($K_1 \times N_p$): for each ray $p$, it stores the segment length inside each traversed voxel (row-wise structure of $A$). Thus, \textbf{PL} collects the non-zero values of $A$, raw-wise.
    
    \item \textbf{VP} ($K_2 \times N_v$): for each voxel $v$, it stores a list of the pixel indices whose rays traverse the voxel.
    
    \item \textbf {VL} ($K_2 \times N_v$):  for each voxel $v$, it stores the segment length of the rays of the previous table (column-wise structure of $A$). Thus, \textbf{VL} collects the non-zero values of $A$, column-wise.
\end{itemize}

$K_1$ dimension of \textbf{PV} and \textbf{PL} is defined by the maximum number of steps of the voxel traversal algorithm is $K_1=V_x+V_y+V_z$ \cite{amanatides1987fast}. 

\begin{figure}{!htbp}
    \centering
    \includegraphics[width=0.9\linewidth]{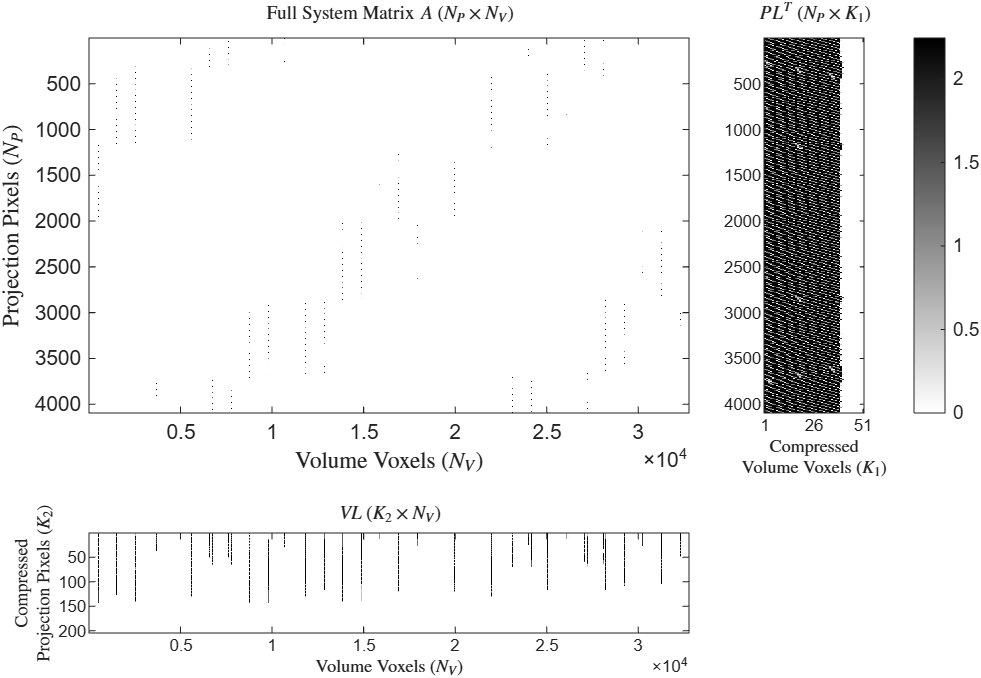}
    \caption{System matrix $A$ and auxiliary $PL, VL$ matrices}
     \label{fig:systemMatrix}
\end{figure}

The number of rows $K_2$ in the last two tables is estimated as follows. A safe and sufficient value for $K_2$ is the maximum number of the image's pixels that a voxel can be projected onto. Let $v_x,v_y,v_z$ be the voxel-size (in mm) and $p_x,p_y$ the pixel size. The worst case is a ray that projects the main diagonal of a voxel on an image. The linear magnification of the projection is equal to $\frac{SID}{SOD}$, where SID is the Source-to-Detector distance and SOD is the Source-to-Object distance. Finally squaring the result, to account for any orientation of the projection on the image, we obtain
\begin{equation}
    K_2= \left \lceil \frac{SID}{SOD} \frac{\sqrt{v_x^2+v_y^2+v_z^2}}{min(p_x,p_y)} \right \rceil ^2,
\end{equation}
where $\lceil x \rceil$ returns the closest larger integer of $x.$ This way the memory requirements are greatly reduced. The auxiliary matrices $\mathbf{PL}$ and $\mathbf{PV}$ in relation to the original system matrix $\mathbf{A} \in \mathbb{R}^{64^2 \times 32^3}$ is visualized in figure \ref{fig:systemMatrix}. The true scale of the compression achieved by this method is described in table \ref{tab:memory_optimization}.

\begin{table}[!htbp]
\caption{Memory Optimization: Original vs. Optimized Representation}
\label{tab:memory_optimization}
\begin{tabular*}{\textwidth}{@{\extracolsep\fill}ccccccc}
\toprule
& & \textit{\textbf{Original}} & \multicolumn{2}{c}{\textit{\textbf{Auxiliary Matrices}}} & \multicolumn{2}{c}{\textit{\textbf{Optim. Metrics}}} \\
\cmidrule(r){3-3} \cmidrule(lr){4-5} \cmidrule(l){6-7}
\shortstack{\textbf{Voxel}\\\textbf{$(N_v$)}} & \shortstack{\textbf{Pixel}\\\textbf{$(N_p)$}} & $\mathbf{A}$ & ${K_1 \times N_p}$ & ${K_2 \times N_v}$ & \textbf{Total} & \textbf{Reduction} \\
\midrule
$32^3$ & $64^2$ & 512 MB & 1.5 MB & 3.2 MB & 9.3 MB & $\sim$55$\times$ \\
$128^3$ & $512^2$ & 2.00 TB & 384 MB & 200 MB & 1.15 GB & $\sim$1,700$\times$ \\
$192^3$ & $512^2$ & 6.75 TB & 576 MB & 243 MB & 1.59 GB & $\sim$4,300$\times$ \\
$256^3$ & $512^2$ & 16.00 TB & 768 MB & 576 MB & 2.62 GB & $\sim$6,200$\times$ \\
\bottomrule
\end{tabular*}
\end{table}

\subsubsection{Forward and Back-Projection Operators using auxiliary tables}
\label{sec:f&bp_c}
For any vector $\mathbf{x}$ with appropriate dimensions, using element-wise product $\odot$ the forward-projection becomes
\begin{equation}
    \mathbf{A} \mathbf{x} = PL_{:,p} \odot x_{PV_{:,p}} \Leftrightarrow (\mathbf{A} \mathbf{x})_p= \sum_{k}PL_{k,p} \cdot x_{PV_{k,p}}
    \label{eq:Ax}
\end{equation}

and the back-projection operator:
\begin{equation}
\label{eq:ATx}
    \mathbf{A^T}\mathbf{x} = VL_{:,v} \odot x_{VP_{:,v}} \Leftrightarrow (\mathbf{A^T}\mathbf{x})_v =\sum_{k}VL_{k,v} \cdot x_{VP_{k,v}}.
\end{equation}

The computation of $A^T Ax$ becomes:
\begin{equation}
\label{Eq:ATAx}
    \mathbf{(A^TA)x}=A^T(Ax)= VL_{:,v} \odot
    (PL_{:,P} \odot x_{PV_{:,p}})_{VP_{:,v}}.
\end{equation}

It is now clear that these arithmetic operators are implemented without explicitly constructing the system matrices.

\input{gd_clean}

\section{ART/SART/SIRT}

\subsection{ART}
The Algebraic Reconstruction Technique (ART) is an iterative solver for the linear system \eqref{eq:matrix_System}, based on the Kaczmarz method. At each iteration, a single ray equation is selected, and the current estimate is orthogonally projected onto the corresponding hyperplane, yielding an updated attenuation volume.

To solve this linear system $\mathbf{A}\boldsymbol{\mu} = \mathbf{Y}$, ART processes the projection equations sequentially. For each ray at pixel $p$ the update rule of the attenuation of the voxels $v \in \mathcal{R}(p)$ is:


\begin{equation}
    \mathbf{\mu}_v^{(t)} = \mathbf{\mu}_v^{(t-1)} + \lambda \frac{Y_p -  \sum_{j} A_{pj} \mathbf{\mu}_{j}^{(t-1)}}{\sum_{j}A_{pj}^2 }  A_{pv} .
\end{equation}

The correction is distributed only among the voxels $v$ intersected by ray $p$, with each voxel update weighted by:
\begin{equation}
    q_p = \sum_j A_{pj}^2
\end{equation}
It is clear that only voxels traversed by ray $p$ are updated, since the corresponding elements of the $\mathbf{A}$ matrix of the remaining voxels are $0$.

In practical implementations, the projection rays are often processed in a randomized order to improve convergence properties and reduce directional artifacts. After each ray update, a non-negativity constraint is applied: $\boldsymbol{\mu}^{(t)}\leftarrow\max(\boldsymbol{\mu}^{(t)},0)$, ensuring physically meaningful attenuation coefficients.

We formulate the operations as:
\begin{align}
    q_p &= \sum_{k} PL_{kp}^2 \label{eq:q_p} \\ 
    r_p &= Y_p - \sum_{k} PL_{kp} \mathbf{\mu}_{PV_{kp}} \label{eq:r_p}
\end{align}
where $PV_{kp}$ denotes the indices of the $k$ voxels intersected by ray $p$.
The update rule of the voxels $\mathbf{v}$  becomes:

\begin{equation}
        \boldsymbol{\mu}_v^{(t)} = \boldsymbol{\mu}_v^{(t-1)} + \lambda \frac{{r}_p^{(t-1)}}{{q}_p } \cdot \mathbf{PL}_{:,p}. 
\end{equation}
This computational implementation is shown in algorithm \ref{alg:art}.


\begin{algorithm}
\caption{Algebraic Reconstruction Technique (ART) Update}
\label{alg:art}
\begin{algorithmic}[1]
\Require Projection data $\mathbf{Y}$, relaxation $\lambda$, current volume $\mu$, PV, PL
\For{$e = 1$ \textbf{to} $N_{epoch}$}
\For{$j = 1$ \textbf{to} $N_{J}$}
\State $N_{p} \gets$ number of pixels in the detector
\State $\text{ray\_order} \gets \text{RandomPermutation}(N_{p})$
\For{$i = 1$ \textbf{to} $N_{p}$}
    \State $p \gets \text{ray\_order}[i]$
    \State $\mathbf{v} \in \mathbb{N}^{K \times 1} \gets PV_{:,p}$ \Comment{Indices of voxels traversed by ray $p$}
    \State $\hat{Y}_p \gets \sum_k \mu_{\mathbf{v}(k)} PL_{k,p}=\mu_{\mathbf{v}} \cdot PL_{:,p}$ \Comment{Forward projection for ray $p$}
    \State $r_p \gets Y_p - \hat{Y}_p $ \Comment{Calculate residual}
    \State $q_p \gets \sum_v PL_{pv}^2$ \Comment{Normalization }
    \If{$q_p > \epsilon$}
    
            \State $\mathbf{\mu}_{\mathbf{v}} \gets \mathbf{\mu}_\mathbf{v} + \lambda \left(\frac{r_p}{q_p}\right) PL_{pv}$ \Comment{Update voxels within Ray}
            \State $\mathbf{\mu}_\mathbf{v} \gets \max(\mathbf{\mu}_\mathbf{v}, 0)$ \Comment{Apply non-negativity}
    \EndIf
\EndFor
\EndFor
\EndFor
\State \Return $\boldsymbol{\mu}$
\end{algorithmic}
\end{algorithm}

\subsection{SART}

The Simultaneous Algebraic Reconstruction Technique (SART) is a variation of the ART algorithm designed to improve convergence properties and reduce the noise and directional artifacts often introduced by updating voxels on a ray-by-ray basis. Instead of immediately applying corrections after evaluating each ray, SART computes the projection errors for all rays, accumulates the corresponding corrections and applies the averaged updates to the whole volume simultaneously, at the end of each projection.

For any given voxel $v$, the corrections from all rays $p$ intersecting $v$ are accumulated and weighted by their path lengths. After all rays have been processed, the updated estimate for voxel $v$, denoted $\mathbf{\mu}_v^{(t)}$, is calculated as:
\begin{equation}
    \mathbf{\mu}_v^{(t)} = \mathbf{\mu}_v^{(t-1)} + \lambda \frac{1}{\sum_{p} A_{pv}} \sum_{p} \frac{A_{pv}}{\sum_{k} A_{pk}} \left( Y_p - \sum_{k} A_{pk} \mathbf{\mu}_{k}^{(t-1)} \right)
    \label{eq}
\end{equation}
where $\lambda$ is the relaxation factor, $k$ iterates over the indices of the voxels traversed by ray $p$, and $A_{pk}$ represents the path length of ray $p$ through voxel $k$. The rightmost term is the current projection residual $r^{(t-1)}_p$ for a given ray $p$. 

This update equation incorporates two normalizations: 
A row-wise normalization which divides each ray's residual by its total path length through the volume and is defined  for ray $p$ as 
\begin{equation}
    q_p = \sum_{k} A_{pk}=\mathbf{A}\textbf{1}
\end{equation}
and a column-wise normalization that divides the accumulated correction for each voxel by the total weight of all rays passing through it and is defined for voxel $v$ as
\begin{equation}
    w_v = \sum_{p} A_{pv}=\textbf{A}^T\textbf{1}.
\end{equation}
Both quantities can be accurately calculated using \eqref{eq:Ax} and \eqref{eq:ATx}.

Using these definitions, the SART update simplifies to:
\begin{equation}
    \mathbf{\mu}_v^{(t)} = \mathbf{\mu}_v^{(t-1)} + \lambda \frac{1}{w_v} \sum_{p} A_{pv} \left( \frac{r^{(t-1)}_p}{q_p} \right),
    \label{eq:smallSART}
\end{equation}
or in matrix notation
\begin{equation}
    \boldsymbol{\mu}^{(t)} = \boldsymbol{\mu}^{(t-1)} + \lambda \mathbf{A} \left( \mathbf{r}^{(t-1)} \oslash{\mathbf{q}} \right) \oslash{\mathbf{w}}
    \label{eq:smallSART}
\end{equation}

In SART, all voxel estimates are updated simultaneously at the end of each projection iteration. Finally, a non-negativity constraint is applied to the entire volume. By updating the valid voxels simultaneously rather than sequentially, SART eliminates the need for randomized ray ordering. 



The aforementioned efficient computational implementation of SART that avoids full matrix construction is provided in Algorithm \ref{alg:sart}. 

\begin{algorithm}
\caption{Simultaneous Algebraic Reconstruction Technique (SART) Update}
\label{alg:sart}
\begin{algorithmic}[1]
\Require Projection data $Y^{(j)}$, relaxation factor $\lambda$, current volume $\mathbf{\mu}$, matrices $PV^{(j)}, PL^{(j)}, VL^{(j)}$ for all views $j$
\For{$e = 1$ \textbf{to} $N_{epoch}$}
    \For{$j = 1$ \textbf{to} $N_{J}$}
        \State $\delta\mathbf{\mu} \gets \mathbf{0} \in \mathbb{R}^{N_v \times 1}$ \Comment{Initialize correction accumulator}
        \State $\mathbf{w} \gets \mathbf{0} \in \mathbb{R}^{N_v \times 1} $ \Comment{Initialize weight accumulator}
        \State $N_{p} \gets$ number of pixels in the detector
        
        \For{$p = 1$ \textbf{to} $N_{p}$}
            \State $\mathbf{v} \in \mathbb{N}^{K \times 1} \gets PV_{:,p}^{(j)}$ \Comment{Indices of voxels traversed by ray $p$}
            \State $\hat{Y}_p \gets \sum_k \mathbf{\mu}_{\mathbf{v}(k)} PL_{k,p}^{(j)}$ \Comment{Forward projection for ray $p$}
            \State $r_p \gets Y_p^{(j)} - \hat{Y}_p$ \Comment{Calculate residual}
            \State $q_p \gets \sum_k PL_{k,p}^{(j)}$ \Comment{Total ray path length}
            
            \If{$q_p > \epsilon$} \Comment{Avoid division by zero}
                \State $\delta\boldsymbol{\mu}_{\mathbf{v}} \gets \delta\boldsymbol{\mu}_{\mathbf{v}} + \left(\frac{r_p}{q_p}\right) \mathbf{PL}_{:,p}^{(j)}$ \Comment{Accumulate corrections}
            \EndIf
        \EndFor
        
        \State $\mathbf{w}_v \gets (\sum_k VL_{k,v}^{(j)}), \quad v=1 \ldots N_v$ \Comment{Column-wise normalization for all $v$}
        \State $\Delta\boldsymbol{\mu} \gets \lambda \left(\delta\boldsymbol{\mu} \oslash \mathbf{w}\right)$ \Comment{Element-wise division}
        \State $\boldsymbol{\mu} \gets \boldsymbol{\mu} + \Delta\boldsymbol{\mu}$ \Comment{Apply simultaneous update}
        \State $\boldsymbol{\mu} \gets \max(\boldsymbol{\mu}, 0)$ \Comment{Enforce non-negativity}
    \EndFor
\EndFor
\State \Return $\boldsymbol{\mu}$
\end{algorithmic}
\end{algorithm}

\subsection{Simultaneous Iterative Reconstruction Technique (SIRT)}
\label{subsec:sirt}

The Simultaneous Iterative Reconstruction Technique (SIRT)
computes a global correction to the reconstructed volume
by processing all projection rays across all views before
updating the voxel attenuation coefficients.
Unlike ART, which updates the volume sequentially after
processing individual rays, and SART, which typically
updates the volume after each projection view, SIRT
accumulates the contributions of all projection residuals
using a fixed volume estimate. The accumulated correction
is then applied simultaneously to all voxels at the end
of each epoch.

For a given ray $p$ in projection view $j$, the forward
projection of the current volume estimate is
\begin{equation}
\hat{Y}_{p,j}^{(t-1)}
=
\sum_{v=1}^{N_v}
A_{pv}^{(j)}\mu_v^{(t-1)},
\label{eq:sirt_forward}
\end{equation}
where $\mathbf{A}^{(j)}$ denotes the system matrix for the $j$ view, whereas the view for projection $\mathbf{Y}$ is denoted vy $j$ in the index.
The projection residual is defined as the difference between 
the measured and estimated projections. Using the auxiliary 
matrices $\mathbf{PV}$ (voxel indices) and $\mathbf{PL}$ (path lengths), 
the residual can be computed matrix-free:
\begin{align}
r_{p,j}^{(t-1)} &= Y_{p,j}-\hat{Y}_{p,j}^{(t-1)} \nonumber \\
&= Y_{p,j} - \sum_{v=1}^{N_v} A_{pv}^{(j)}\mu_v^{(t-1)} \nonumber \\
&= Y_{p,j} - \sum_{k\in\mathcal{I}_{p,j}} PL_{k,p,j} \mu_{PV_{k,p,j}}^{(t-1)},
\label{eq:sirt_residual}
\end{align}
where $\mathcal{I}_{p,j}=\{k:PV_{k,p,j}>0\}$ denotes the valid ray-voxel intersections.

Two normalization factors are introduced.
The row-wise normalization $q_{p,j}$ is the total
intersection length of ray $p$ in view $j$:
\begin{equation}
q_{p,j}
=
\sum_{v=1}^{N_v} A_{pv}^{(j)}
=
\sum_{k\in\mathcal{I}_{p,j}} PL_{k,p,j} \quad \Leftrightarrow \quad \mathbf{q}^{(j)} = \mathbf{A}^{(j)}\mathbf{1}_{N_v}.
\label{eq:sirt_q}
\end{equation}

The column-wise normalization $w_v$ is the total intersection length accumulated by voxel $v$ over all rays and all projection views ($N_J$):
\begin{equation}
w_v = \sum_{j=1}^{N_J} \sum_{p=1}^{N_p} A_{pv}^{(j)} \quad \Leftrightarrow \quad \mathbf{w} = \sum_{j=1}^{N_J} \left(\mathbf{A}^{(j)}\right)^T \mathbf{1}_{N_p}.
\label{eq:sirt_w}
\end{equation}
Here, $\mathbf{1}_{N_v}$ and $\mathbf{1}_{N_p}$
are all-ones vectors of dimensions $N_v$ and
$N_p$, respectively.

For voxels with $w_v>0$, the SIRT update is
\begin{equation}
\mu_v^{(t)} =
\mu_v^{(t-1)} + \frac{\lambda}{w_v}
\sum_{j=1}^{N_J} \sum_{\substack{p=1\\q_{p,j}>0}}^{N_p}
A_{pv}^{(j)} \frac{r_{p,j}^{(t-1)}}{q_{p,j}},
\label{eq:sirt_update}
\end{equation}
where $\lambda>0$ denotes the relaxation factor.
Rays with zero total intersection length do not
contribute to the update.

The numerator accumulates the normalized residual
contributions from all rays intersecting voxel $v$.
The denominator $w_v$ normalizes this accumulated
correction by the total intersection length
associated with the voxel.

For a compact matrix representation, the system
matrices of all projection views are stacked as
\begin{equation}
\mathcal{A}
=
\begin{bmatrix}
\mathbf{A}^{(1)}\\
\mathbf{A}^{(2)}\\
\vdots\\
\mathbf{A}^{(N_J)}
\end{bmatrix}
\in \mathbb{R}^{N_J N_p\times N_v}.
\end{equation}

Let $\mathcal{Y}$ and $\mathcal{R}$ denote the
corresponding stacked measurement and residual vectors:
\begin{equation}
\mathcal{Y} =
\begin{bmatrix}
\mathbf{Y}^{(1)} \\
\mathbf{Y}^{(2)} \\
\vdots \\
\mathbf{Y}^{(N_J)}
\end{bmatrix},
\qquad
\mathcal{R}^{(t-1)} =
\begin{bmatrix}
\mathbf{r}^{(1)} \\
\mathbf{r}^{(2)} \\
\vdots \\
\mathbf{r}^{(N_J)}
\end{bmatrix}
= \mathcal{Y} - \mathcal{A} \boldsymbol{\mu}^{(t-1)}.
\end{equation}

The stacked normalization vectors satisfy
\begin{equation}
\mathbf{q}=\mathcal{A}\mathbf{1}_{N_v},
\qquad
\mathbf{w}=\mathcal{A}^T\mathbf{1}_{N_J N_p}.
\end{equation}

The SIRT update can then be expressed as
\begin{equation}
\boldsymbol{\mu}^{(t)} =
\boldsymbol{\mu}^{(t-1)} + \lambda \left[
\mathcal{A}^T \left(\mathcal{R}^{(t-1)}\oslash\mathbf{q}\right)\right]
\oslash\mathbf{w}.
\label{eq:sirt_matrix}
\end{equation}

The element-wise divisions ($\oslash$) are evaluated only for
positive normalization factors. Unobserved voxels
with $w_v=0$ retain their previous values.

Since all ray residuals are evaluated using the
same volume estimate $\boldsymbol{\mu}^{(t-1)}$,
the update is independent of the order in which
rays and views are processed, apart from
floating-point accumulation effects.

Finally, non-negativity is imposed:
\begin{equation}
\boldsymbol{\mu}^{(t)}
\gets
\max\left(
\boldsymbol{\mu}^{(t)},\mathbf{0}
\right).
\end{equation}

The corresponding matrix-free implementation
is summarized in Algorithm~\ref{alg:sirt}.

\begin{algorithm}
\caption{Simultaneous Iterative Reconstruction Technique (SIRT) Update}
\label{alg:sirt}
\begin{algorithmic}[1]
\Require Cell array of projections $Y$, relaxation $\lambda$, current volume $\boldsymbol{\mu}$, Cell arrays of voxel indices $PV$ and path lengths $PL$
\State $N_J \gets$ number of projection views
\State Initialize correction array $\Delta \boldsymbol{\mu} \gets \mathbf{0}$ \Comment{Matches volume size}
\State Initialize weight array $\mathbf{w} \gets \mathbf{0}$ \Comment{Matches volume size}
\For{$j = 1$ \textbf{to} $N_J$}
    \State $N_{rays} \gets$ number of pixels in view $j$
    \For{$p = 1$ \textbf{to} $N_{rays}$}
        \State $\hat{Y}_{p,j} \gets 0$
        \For{$k \in \mathcal{I}_{p,j}$} \Comment{Forward projection for ray $p$}
            \State $v \gets PV_{k,p,j}$
            \State $\hat{Y}_{p,j} \gets \hat{Y}_{p,j} + \mu_v \cdot PL_{k,p}^{(j)}$ 
        \EndFor
        \State $r_{p,j} \gets Y_{p,j} - \hat{Y}_{p,j}$ \Comment{Calculate residual}
        \State $q_{p,j} \gets \sum_k PL_{k,p}^{(j)}$ \Comment{Row normalization}
        \If{$q_{p,j} > \epsilon$}
            \For{$k \in \mathcal{I}_{p,j}$} \Comment{Accumulate updates}
                \State $v \gets PV_{k,p}^{(j)}$
                \State $\Delta \mu_v \gets \Delta \mu_v + \left(\frac{r_{p,j}}{q_{p,j}}\right) PL_{k,p}^{(j)}$ 
                \State $w_v \gets w_v + PL_{k,p}^{(j)}$ 
            \EndFor
        \EndIf
    \EndFor
\EndFor
\For{\textbf{each} voxel $v$} \Comment{Global Simultaneous Update}
    \If{$w_v > 0$}
        \State $\mu_v \gets \mu_v + \lambda \left(\frac{\Delta \mu_v}{w_v}\right)$ \Comment{Update valid voxels}
    \EndIf
    \State $\mu_v \gets \max(\mu_v, 0)$ \Comment{Apply non-negativity}
\EndFor
\State \Return $\boldsymbol{\mu}$
\end{algorithmic}
\end{algorithm}

\subsection{GMRES}
The generalized minimal residual method (GMRES) is a Krylov subspace algorithm designed to solve square, invertible, and generally non-symmetric linear systems. In the context of tomographic reconstruction, the problem is formulated iteratively, with each iteration solving for the current volume correction $\delta \boldsymbol{\mu}$ given the current residual: $A \delta \boldsymbol{\mu} = \mathbf{r}$, where $\mathbf{r} = \mathbf{Y} - \mathbf{A}\boldsymbol{\mu}$. Here, $\mathbf{A}$ and $\mathbf{r}$ act as placeholders that represent the system matrix and residual for the chosen update scheme. Because the tomographic system matrix is typically overdetermined and rectangular, we transform the problem into the normal equations to render it square and suitable for GMRES:
\begin{equation}
    \mathbf{A}^T\mathbf{A}\delta\boldsymbol{\mu} = \mathbf{A}^T\mathbf{r}.
    \label{eq:gmres_1}
\end{equation}

This GMRES formulation accommodates both single-view and multi-view processing, supporting either a per-view (sequential) or a per-epoch update mode. 
In the single-view approach, the system matrix and residual correspond to the geometry of the current view $j$, such that $\mathbf{A} = \mathbf{A}^{(j)}$ and $\mathbf{r} = \mathbf{r}^{(j)}$. The attenuation volume $\boldsymbol{\mu}$ is updated by solving for $\delta \boldsymbol{\mu}$ before the next view is evaluated, as described in Algorithm \ref{alg:gmres_single}.

Alternatively, in the simultaneous multi-view approach, a single, comprehensive solution for $\delta \boldsymbol{\mu}$ is performed at the end of the full projection sequence. The complete system of equations used in this implementation substitutes $\mathbf{A} = \mathcal{A}$ and $\mathbf{r} = \mathcal{R}$, utilizing the stacked multi-view matrices previously defined in Section \ref{subsec:sirt}. 

We define $\mathbf{b} = \mathbf{A}^T \mathbf{r}$ as the back-projected residual, allowing us to express Equation \eqref{eq:gmres_1} as $\mathbf{A}^T \mathbf{A} \delta\boldsymbol{\mu} = \mathbf{b}$. In the simultaneous multi-view formulation, the block global transpose yields $\mathbf{b}$ as the accumulated back-projections over all views:
\begin{equation}
\label{eq:b_gmres}
    \mathbf{b} = \mathcal{A}^T \mathcal{R} = \sum_{j=1}^{N_J} (\mathbf{A}^{(j)})^T \mathbf{r}^{(j)}.
\end{equation}

The system matrices $\mathbf{A}^{(j)}$ are never constructed, nor are the forward and back-projection operations applied directly via $({\mathbf{A}^{(j)}})^T$. Instead, these operations are implemented implicitly by using the auxiliary tables defined in Section \ref{sec:comp_matrices} following the procedure in Subsubsection \ref{sec:f&bp_c}.

Using the initial back-projected residual $\mathbf{b}$, the algorithm initializes the first orthonormal basis vector $\mathbf{q}_1$. It then constructs the orthonormal basis $Q_k$ of the Krylov subspace via Arnoldi iteration, using our implicitly derived implementation for the operator $\mathbf{A}^T \mathbf{A}$ as in \eqref{eq:ATx}. At each step, it checks convergence by solving the low-dimensional least-squares problem $\tilde{\mathbf{y}} = \arg\min_{\mathbf{y}} \Vert H_k \mathbf{y} - \mathbf{g}_k \Vert_2$. Once the fractional GMRES residual $\rho_k= \frac{\Vert H_k \tilde{\mathbf{y}} - \mathbf{g}_k\Vert_2}{\Vert \mathbf{g}_k \Vert_2} \le \varepsilon$, the method terminates. The solution yields the volume correction $\delta\boldsymbol{\mu} = Q_k \tilde{\mathbf{y}}$, updating $\boldsymbol{\mu} \gets \boldsymbol{\mu} + \delta\boldsymbol{\mu}$. 

The update is calculated and applied at the end of every view, or a global correction is calculated using all views and the volume is updated only once at the end of each epoch, as described in Algorithm \ref{alg:gmres_multi}. 

The computational implementation of the GMRES algorithm using our optimized matrices of Section \ref{sec:comp_matrices}, suitable for both sequential and simultaneous approaches by adjusting the inputs, is detailed in Algorithm \ref{alg:gmres_reconstruction}.

\begin{algorithm}
\caption{GMRES single-view reconstruction}
\label{alg:gmres_single}
\begin{algorithmic}[1]
\State $\boldsymbol{\mu}\gets \mathbf{0}$
    \For{epoch $e$} \Comment{Initialization}
        \For{view $j = 1 \text{ to } N_{\text{J}}$}
            \State $\boldsymbol{\mu} = \text{GMRES}\left (\boldsymbol{\mu}, \mathbf{Y}^{(j)}, \mathbf{PL}^{(j)}, \mathbf{PV}^{(j)}, \mathbf{VL}^{(j)}, \mathbf{VP}^{(j)}, N=1, K, \varepsilon \right )$
        \EndFor
    \EndFor
\end{algorithmic}
\end{algorithm}

\begin{algorithm}
\caption{GMRES multi-view reconstruction}
\label{alg:gmres_multi}
\begin{algorithmic}[1]
\State $\boldsymbol{\mu}\gets \mathbf{0}$
\State $\mathbf{PL} = \{\mathbf{PL}^{(1)},\ldots, \mathbf{PL}^{(N_{J})}\}, \quad \mathbf{PV} = \{\mathbf{PV}^{(1)},\ldots, \mathbf{PV}^{(N_{J})}\}$
\State $\mathbf{VL} = \{\mathbf{VL}^{(1)}, \ldots, \mathbf{VL}^{(N_{J})}\}, \quad \mathbf{VP} = \{\mathbf{VP}^{(1)},\ldots, \mathbf{VP}^{(N_{J})}\}$
    \For{epoch $e$} \Comment{Initialization}
            \State $\boldsymbol{\mu} = \text{GMRES}\left (\boldsymbol{\mu}, \mathbf{Y}, \mathbf{PL}, \mathbf{PV}, \mathbf{VL}, \mathbf{VP}, N=N_{J}, K, \varepsilon \right )$
    \EndFor
\end{algorithmic}
\end{algorithm}

\begin{algorithm}
\caption{GMRES correction for sparse-view tomographic reconstruction}
\label{alg:gmres_reconstruction}
\begin{algorithmic}[1]
\Require Current volume $\boldsymbol{\mu}$, measured projections
$\mathbf{Y}$, pixel-wise lengths $PL$ and voxel indices $PV$, voxel-wise lengths $VL$ and pixel indices $VP$, number of views $N$, maximum number of iterations $K$, tolerance $\varepsilon$
\Ensure Updated volume $\boldsymbol{\mu}$

\AlgSection{projectionBG}{projectionText}{1. Projection residual and initialization}

\State $\widehat{\mathbf{Y}}\gets PL_{:,p} \odot \boldsymbol{\mu}_{PV_{:,p}}$ 
\Comment{Implemented using \eqref{eq:Ax}}
\State $\mathbf{r}^{(j)}\gets\mathbf{Y}^{(j)}-\widehat{\mathbf{Y}}^{(j)}$
    \Comment{Projection residual}

\State $b_v \gets \sum_{j=1}^{N} VL_{:,v}^{(j)} \cdot r_{VP_{:,v}}^{(j)}$
    \Comment{$\mathbf{b}\gets A^{T}\mathbf{r}$, using \eqref{eq:ATx} for all views $j$ and sum}

\State $\delta\boldsymbol{\mu}_0\gets\mathbf{0}$
    \Comment{Initial correction}

\State $\mathbf{s}_0\gets\mathbf{b}$
    \Comment{Initial GMRES residual}

\State $\beta\gets\|\mathbf{b}\|_2$

\If{$\beta>0$}
    \State $\mathbf{q}_1\gets\mathbf{b}/\beta$
    \State $H \gets \mathbf{0}\in\mathbb{R}^{(K+1)\times K}$ \Comment{Initialize the Hessenberg matrix}
    \State $\mathbf{g}\gets \beta\mathbf{e}_1\in\mathbb{R}^{K+1}$

\For{$k=1$ \textbf{to} $K$}
\AlgSection{arnoldiBG}{arnoldiText}{2. Arnoldi basis construction}
        
        \State $\mathbf{w}\gets VL_{:,v} \odot (PL_{:,p} \odot (\mathbf{q}_k)_{PV_{:,p}})_{VP_{:,v}}$
            \Comment{Apply $A^{T}A$} using \eqref{Eq:ATAx} and \eqref{eq:b_gmres}. 
        
        \For{$i=1$ \textbf{to} $k$}
            \State $h_{i,k}\gets\mathbf{q}_i^{T}\mathbf{w}$
            \State $\mathbf{w}\gets \mathbf{w}-h_{i,k}\mathbf{q}_i$
        \EndFor

        \State $h_{k+1,k}\gets\|\mathbf{w}\|_2$

\AlgSection{leastSquaresBG}{leastSquaresText}{3. Least-squares solution and convergence}

        \State $H_k\gets H(1:k+1,1:k)$
        \State $\mathbf{g}_k\gets\mathbf{g}(1:k+1)$

        \State Obtain $\tilde{\mathbf{y}}$ solving $H_k\mathbf{y}=\mathbf{g}_k$  
        \Comment{ $\displaystyle \tilde{\mathbf{y}}\gets \arg\min_{\mathbf{y}\in\mathbb{R}^{k}} \left\| H_k\mathbf{y}-\mathbf{g}_k \right\|_2$}

        \State $\displaystyle \rho_k\gets \frac{\| \mathbf{g}_k- H_k\tilde{\mathbf{y}} \|_2}{\beta}$

        \If{$\rho_k\leq\varepsilon$}
            \State \textbf{break} \Comment{Check GMRES convergence}
        \EndIf

        \If{$h_{k+1,k}=0$}
            \State \textbf{break} \Comment{Check Arnoldi Convergence}
        \EndIf

        \State $\mathbf{q}_{k+1}\gets \mathbf{w}/h_{k+1,k}$

    \EndFor

\AlgSection{updateBG}{updateText}{4. Volume correction}

    \State $Q_k\gets [\,\mathbf{q}_1,\ldots,\mathbf{q}_k\,]$

    \State $\delta\boldsymbol{\mu}\gets Q_k\tilde{\mathbf{y}}$
        \Comment{Computed volume correction}

    \State $\boldsymbol{\mu}\gets \boldsymbol{\mu}+\delta\boldsymbol{\mu}$

\EndIf

\end{algorithmic}
\end{algorithm}

\subsection{LSQR}

Given the current reconstruction $\boldsymbol{\mu}$, 
LSQR computes a correction $\delta\boldsymbol{\mu}$ by approximately solving
\[
    \min_{\delta\boldsymbol{\mu}}
    \left\|\mathbf A\delta\boldsymbol{\mu}-\mathbf r\right\|_2 .
\]
To this end, the basis of Krylov-like projection space and voxel space are constructed $\{\mathbf{u_1},\ldots ,\mathbf{u_k}\}$ and $\{\mathbf{v_1},\ldots ,\mathbf{v_k}\}$ respectively. The construction is recursive with each $\mathbf{v_i}$ generated by back-projecting $\mathbf{u_{i-1}}$ and vice-versa.
Golub-Kahan bidiagonalization is initialized by
\[
    \beta_1=\|\mathbf r\|_2,\qquad
    \mathbf u_1=\frac{\mathbf r}{\beta_1},\qquad
    \alpha_1=\|\mathbf A^T\mathbf u_1\|_2,\qquad
    \mathbf v_1=\frac{\mathbf A^T\mathbf u_1}{\alpha_1}.
\]
For $i=1,\ldots,k$, the subsequent vectors are generated by
\[
\begin{aligned}
    \widetilde{\mathbf u}_{i+1}
        &=\mathbf A\mathbf v_i-\alpha_i\mathbf u_i,
    &\beta_{i+1}&=\|\widetilde{\mathbf u}_{i+1}\|_2,
    &\mathbf u_{i+1}&=
        \frac{\widetilde{\mathbf u}_{i+1}}{\beta_{i+1}},\\
    \widetilde{\mathbf v}_{i+1}
        &=\mathbf A^T\mathbf u_{i+1}-\beta_{i+1}\mathbf v_i,
    &\alpha_{i+1}&=\|\widetilde{\mathbf v}_{i+1}\|_2,
    &\mathbf v_{i+1}&=
        \frac{\widetilde{\mathbf v}_{i+1}}{\alpha_{i+1}} .
\end{aligned}
\]
After $k$ steps,
\[
    V_k=[\mathbf v_1,\ldots,\mathbf v_k],\qquad
    U_{k+1}=[\mathbf u_1,\ldots,\mathbf u_{k+1}],
\]
and
\[
    \mathbf A V_k=U_{k+1}B_k,\qquad
    B_k=
    \begin{bmatrix}
    \alpha_1      &               &        \\
    \beta_2       & \alpha_2      &        \\
                  & \ddots        & \ddots \\
                  &               & \beta_{k+1}
    \end{bmatrix}.
\]
Restricting the correction to the Krylov subspace,
$\delta\boldsymbol{\mu}_k=V_k\mathbf z_k$, and using
$\mathbf r=\beta_1\mathbf u_1
=\beta_1U_{k+1}\mathbf e_1$, gives
\[
    \|\mathbf A\delta\boldsymbol{\mu}_k-\mathbf r\|_2
    =
    \|B_k\mathbf z_k-\beta_1\mathbf e_1\|_2,
\]
since $U_{k+1}^TU_{k+1}=I$. Therefore,
\[
    \widetilde{\mathbf z}_k=
    \arg\min_{\mathbf z}
    \|B_k\mathbf z-\beta_1\mathbf e_1\|_2,
    \qquad
    \boxed{\delta\boldsymbol{\mu}_k
    =V_k\widetilde{\mathbf z}_k},
\]
after which the reconstruction is updated as
$\boldsymbol{\mu}\leftarrow
\boldsymbol{\mu}+\delta\boldsymbol{\mu}_k$.

\noindent
\textbf{Incremental LSQR reconstruction.}
LSQR can be implemented without explicitly constructing the projection-space
basis $U_{k+1}$, the voxel-space basis $V_k$, or the bidiagonal matrix $B_k$.

The iteration is initialized as
\[
    \beta_1=\|\mathbf r\|_2,\qquad
    \mathbf u_1=\frac{\mathbf r}{\beta_1},
\]
\[
    \widetilde{\mathbf v}_1=\mathbf A^T\mathbf u_1,\qquad
    \alpha_1=\|\widetilde{\mathbf v}_1\|_2,\qquad
    \mathbf v_1=\frac{\widetilde{\mathbf v}_1}{\alpha_1}.
\]
The quantities used for the incremental solution are initialized by
\[
    \mathbf w_1=\mathbf v_1,\qquad
    \bar{\rho}_1=\alpha_1,\qquad
    \bar{\phi}_1=\beta_1,\qquad
    \delta\boldsymbol{\mu}_0=\mathbf 0 .
\]

At iteration $k$, one Golub-Kahan step is performed using only forward
and backprojection:
\[
\begin{aligned}
    \widetilde{\mathbf u}_{k+1}
        &=\mathbf A\mathbf v_k-\alpha_k\mathbf u_k, &
    \beta_{k+1}
        &=\|\widetilde{\mathbf u}_{k+1}\|_2, &
    \mathbf u_{k+1}
        &=\frac{\widetilde{\mathbf u}_{k+1}}{\beta_{k+1}},\\
    \widetilde{\mathbf v}_{k+1}
        &=\mathbf A^T\mathbf u_{k+1}-\beta_{k+1}\mathbf v_k, &
    \alpha_{k+1}
        &=\|\widetilde{\mathbf v}_{k+1}\|_2, &
    \mathbf v_{k+1}
        &=\frac{\widetilde{\mathbf v}_{k+1}}{\alpha_{k+1}} .
\end{aligned}
\]
A Givens transformation eliminates the new bidiagonal coefficient
$\beta_{k+1}$ through
\[
    \rho_k=\sqrt{\bar{\rho}_k^2+\beta_{k+1}^2},
    \qquad
    c_k=\frac{\bar{\rho}_k}{\rho_k},
    \qquad
    s_k=\frac{\beta_{k+1}}{\rho_k},
\]
and generates
\[
    \theta_{k+1}=s_k\alpha_{k+1},\qquad
    \bar{\rho}_{k+1}=-c_k\alpha_{k+1},
\]
\[
    \phi_k=c_k\bar{\phi}_k,\qquad
    \bar{\phi}_{k+1}=s_k\bar{\phi}_k .
\]
The volume correction is then updated directly as
\[
    \delta\boldsymbol{\mu}_k
    =
    \delta\boldsymbol{\mu}_{k-1}
    +\frac{\phi_k}{\rho_k}\mathbf w_k
\]
with the voxel-space search direction updated according to
\[
    \mathbf w_{k+1}
    =
    \mathbf v_{k+1}
    -\frac{\theta_{k+1}}{\rho_k}\mathbf w_k .
\]
Thus, LSQR requires only the current vectors
$\mathbf u_k$, $\mathbf v_k$, and $\mathbf w_k$ and the accumulated correction
$\delta\boldsymbol{\mu}_k$, while each iteration requires one forward
projection $\mathbf A\mathbf v_k$ and one backprojection
$\mathbf A^T\mathbf u_{k+1}$.
The above implementation of LSQR with Givens rotations is shown in algorithm \ref{alg:lsqr_reconstruction}.

LSQR can be implemented with either per-view, or per epoch update, in a manner identical to that of GMRES.

\definecolor{projectionBG}{RGB}{224,239,255}
\definecolor{bidiagBG}{RGB}{231,245,233}
\definecolor{givensBG}{RGB}{255,243,224}
\definecolor{constraintBG}{RGB}{243,229,245}

\definecolor{projectionText}{RGB}{25,80,140}
\definecolor{bidiagText}{RGB}{35,105,55}
\definecolor{givensText}{RGB}{150,90,20}
\definecolor{constraintText}{RGB}{105,50,130}

\begin{algorithm}
\caption{LSQR correction for sparse-view tomographic reconstruction}
\label{alg:lsqr_reconstruction}
\begin{algorithmic}[1]

\Require Current volume $\boldsymbol{\mu}$, measured projections
$\mathbf{Y}$, computed projections $\mathbf{\widehat{Y}}$, pixel-wise lengths $\mathbf{PL}$ and voxel indices $\mathbf{PV}$, voxel-wise lengths $\mathbf{VL}$ and pixel indices $\mathbf{VP}$,
 number of views $N_{J}$, maximum number of iterations $K$, tolerance $\varepsilon$
\Ensure Updated volume $\boldsymbol{\mu}$

\AlgSection{projectionBG}{projectionText}
{1. Initialization and Initial Backprojection}

\State $\mathbf{u} = \mathbf{Y} - \mathbf{\widehat{Y}}$ \Comment{Initial residual}
\State $\beta = \| \mathbf{u} \|$
\State $\mathbf{u} = \mathbf{u} / \beta $ 

\State $\mathbf{v}_v = VL_{:,v} \odot \mathbf{u}_{VP_{:,v}}$, for all voxels $v$
\Comment{Initial backprojection, $\mathbf{v} = A^T \mathbf{u}$, as in \eqref{eq:ATx}}
\State $\alpha = \| \mathbf{v} \|$
\State $\mathbf{v} = \mathbf{v} / \alpha$

\State $\mathbf{w} = \mathbf{v}$ \Comment{Initialize search direction for Givens Rotation}
\State $\bar{\rho} = \alpha$
\State $\bar{\phi} = \beta$

\For{$k=1$ \textbf{to} $K$}
\AlgSection{bidiagBG}{bidiagText}
{2. Bidiagonalization}

        \State $\mathbf{p}_p \gets PL_{:,p} \odot \mathbf{v}_{PV_{:,p}}$, for all pixel $p$ 
        \Comment{Forward projection $\mathbf{p} = A\mathbf{v}$, as in to \eqref{eq:Ax}}
        \State $\mathbf{u} \gets \mathbf{p} - \alpha \mathbf{u}$
        \State $\beta \gets \| \mathbf{u} \|$
        \State $\mathbf{u} \gets \mathbf{u} / \beta$
        
        \State $\mathbf{q} \gets \mathbf{VL}_{:,v} \odot \mathbf{u}_{VP_{:,v}}$ \Comment{Backprojection $\mathbf{q} = A^T\mathbf{u}$}
        \State $\mathbf{v} \gets \mathbf{q} - \beta \mathbf{v}$
        \State $\alpha \gets \| \mathbf{v} \|$
        \State $\mathbf{v} \gets \mathbf{v} / \alpha$

\AlgSection{givensBG}{givensText}
{3. Orthogonal Transformation (Givens Rotation)}

        \State $\rho \gets \sqrt{\bar{\rho}^2 + \beta^2}$ \Comment{Calculate plane rotation}
        \State $c \gets \bar{\rho} / \rho$
        \State $s \gets \beta / \rho$
        \State $\theta \gets s \cdot \alpha$
        \State $\phi \gets c \cdot \bar{\phi}$
        
        \State $\bar{\rho} \gets -c \cdot \alpha$ \Comment{Update scalars for next iteration}
        \State $\bar{\phi} \gets s \cdot \bar{\phi}$

\AlgSection{updateBG}{updateText}
{4. Volume Update and Convergence}

        \State $\boldsymbol{\mu} \gets \boldsymbol{\mu} + (\phi / \rho) \mathbf{w}$ \Comment{Update volume solution}
        \State $\mathbf{w} \gets \mathbf{v} - (\theta / \rho) \mathbf{w}$ \Comment{Update search direction}
        
        \If{$|\bar{\phi}| \leq \varepsilon$}
            \State \textbf{break} \Comment{Check LSQR convergence}
        \EndIf

\EndFor

\AlgSection{constraintBG}{constraintText}
{5. Non-negativity Constraint}

    \State $\boldsymbol{\mu} \gets \max(\boldsymbol{\mu}, 0)$ \Comment{Set negative voxel values to 0}

\end{algorithmic}
\end{algorithm}

\section{Evaluation Metrics}
To evaluate the volumetric segmentation and reconstruction performance, we employ a combination of intensity-based, geometric, and projection-domain metrics.

\subsection{Volumetric Overlap and Segmentation}
We assess segmentation accuracy using the Intersection over Union (IoU). Each reconstructed volume $\mu$ is binarized automatically without ground-truth intervention. The volume is smoothed with a 3D Gaussian filter ($\sigma = 0.5$ voxels) and normalized by its $99.9$th percentile to mitigate outlier influence. We calculate a base threshold $t_1$ using Otsu's method on the entire volume to isolate the background. A second threshold $t_2$ is computed via Otsu's method applied strictly to voxels exceeding $t_1$. This isolates the vessel cores. We extract the final binary segmentation $S$ via hysteresis thresholding. Voxels exceeding $t_2$ act as seeds and grow into 26-connected regions for all adjacent voxels exceeding $0.8 t_2$. The IoU is then defined as
\begin{equation}
    \mathrm{IoU} = \frac{|S \cap G|}{|S \cup G|}
\end{equation}
where $G$ is the ground-truth vessel mask. We complement this with a soft IoU that operates directly on the normalized intensities, using the algebraic product for the intersection space:
\begin{equation}
    \mathrm{IoU}_{soft}(S, G) = \frac{\sum S \cdot G}{\sum S + \sum G - \sum S \cdot G + \epsilon}
\end{equation}

\subsection{Contrast-to-Noise Ratio (CNR)}
The CNR quantifies vessel prominence relative to background noise in the continuous reconstructed volume. It is defined as
\begin{equation}
    \mathrm{CNR} = \frac{\bar{\mu}_G - \bar{\mu}_B}{\sigma_B}
\end{equation}
where $\bar{\mu}_G$ is the mean reconstructed intensity of the ground-truth vessel voxels. The background region $B$ is the entire volume excluding a two-voxel dilation of $G$. This exclusion margin prevents partial-volume blurring at the vessel boundaries from corrupting the background mean $\bar{\mu}_B$ and standard deviation $\sigma_B$.

\subsection{Geometric and Spatial Errors}
We use the Euclidean distance transform $D(\mathbf{x})$ to compute spatial deviations. The Chamfer distance measures the mean bidirectional Euclidean distance between the active voxels of the reconstructed segmentation $S$ and the ground truth $G$:
\begin{equation}
    d_{Chamfer} = \frac{1}{2} \left( \text{mean}(D_G(S)) + \text{mean}(D_S(G)) \right)
\end{equation}

Topological accuracy is evaluated on the vessel centerlines. The reconstruction is binarized with
Otsu's threshold (after min--max normalization), giving $B_P$; the ground truth $B_G$ is binary.
The skeletons $\mathcal{S}_P = \{\mathbf{p}_i\}_{i=1}^{N_R}$ and $\mathcal{S}_G = \{\mathbf{g}_j\}_{j=1}^{N_G}$
are obtained by 3D morphological thinning, discarding branches shorter than 10 voxels.
For each $\mathbf{p}_i$ let $\mathbf{c}(i)$ be its nearest ground-truth skeleton point, and for each $\mathbf{g}_j$ let $\mathbf{c}'(j)$, obtained from the Euclidean distance transform of the other skeleton.
The skeleton positional error is the mean distance over the matches in both directions,
\begin{equation}
    E^{Pos} = \frac{\bar h}{N_R + N_G}\left(\sum_{i=1}^{N_R}\|\mathbf{p}_i-\mathbf{c}(i)\|_2
            + \sum_{j=1}^{N_G}\|\mathbf{g}_j-\mathbf{c}'(j)\|_2\right),
\end{equation}
with distances in voxels and $\bar h = (h_x h_y h_z)^{1/3}$ the mean voxel size.

The local vessel diameter at a point $\mathbf{x}$ of a skeleton is $d_P(\mathbf{x}) = 2\,D_{\bar B_P}(\mathbf{x})$,
where $D_{\bar B_P}(\mathbf{x})$ is the Euclidean distance (in mm, with the voxel spacing) from $\mathbf{x}$
to the nearest background voxel of $B_P$; $d_G$ is defined in the same way on $B_G$.
The signed diameter errors of the matched points are
$\Delta d_i = d_P(\mathbf{p}_i) - d_G(\mathbf{c}(i))$ and $\Delta d'_j = d_P(\mathbf{c}'(j)) - d_G(\mathbf{g}_j)$,
and the diameter accuracy is the mean of the two directions,
\begin{align}
    E^{d}_{\mathrm{MAE}} &= \frac12\left(\frac{1}{N_R}\sum_{i=1}^{N_R}|\Delta d_i|
                         + \frac{1}{N_G}\sum_{j=1}^{N_G}|\Delta d'_j|\right), \quad \\
    E^{d}_{\mathrm{RMSE}} &= \frac12\left(\sqrt{\frac{1}{N_R}\sum_{i=1}^{N_R}\Delta d_i^2}
                         + \sqrt{\frac{1}{N_G}\sum_{j=1}^{N_G}\Delta d'^2_j}\right).
\end{align}

\subsection{Held-Out View Similarity}
We project the reconstructed volume to a simulated C-arm view outside the acquired angular range. We compute line integrals through the reconstructed volume and compare them to the noise-free reference line integrals of the ground-truth phantom. Projection fidelity is quantified over the detector pixels using the Structural Similarity Index Measure (SSIM), the Normalized Root Mean Square Error (NRMSE), and the Pearson correlation coefficient.

\subsection*{Held-out view similarity}

One C-arm view that is not used for the reconstruction is simulated, the same for all
models: primary angle $100^\circ$, secondary angle $30^\circ$. It lies outside the acquired
primary range ($0^\circ$-$90^\circ$), and its central ray is $10^\circ$ from that of the
nearest acquired view. The reference image is the line integral of the
ground truth volume, and the prediction is the line integral of
the reconstruction along the same rays. Over the set $D$ of detector pixels whose rays cross the volume we report the mean structural similarity
\begin{equation}
  \mathrm{SSIM} = \frac{1}{|D|}\sum_{i\in D}
  \frac{(2 m^p_i m^r_i + C_1)(2 s^{pr}_i + C_2)}
       {\bigl((m^p_i)^2 + (m^r_i)^2 + C_1\bigr)\bigl((s^p_i)^2 + (s^r_i)^2 + C_2\bigr)},
\end{equation}
with local means $m$, standard deviations $s$ and covariance $s^{pr}$ in a Gaussian window
($\sigma = 1.5$ pixels), $C_1 = (0.01L)^2$, $C_2 = (0.03L)^2$ and $L$ the dynamic range of the
reference over $D$, together with
$\mathrm{NRMSE} = \lVert p_D - r_D \rVert_2 / \lVert r_D \rVert_2$ and the Pearson
correlation of $p_D$ and $r_D$.

\subsection{Statistical Rank Agreement}
To ensure method rankings remain consistent across different phantom models and are not biased by threshold parameter selection, we compute rank agreement metrics. We evaluate the method rankings using Spearman's rank correlation $\rho$ and Kendall's rank correlation coefficient $\tau_b$.

\section{Results}
\label{sec:results}~

\subsection{Patient-derived coronary artery models}
The coronary vessel geometries were derived from digitized real-patient data and represented as triangulated surface meshes in STL format. We selected five models, Heart0064, Heart0069, Heart0088, Heart0109, and Heart0148, created by the Visible Heart Laboratories in partnership with LifeSource and the University of Minnesota Anatomy Bequest Program \cite{3Dmodels}. To evaluate reconstruction quality across these models, all methods were executed using an identical set of 10 projection angles.

The following methods participated in the comparisons:
\begin{itemize}
\item Gradient Descent framework with update global step size $\alpha$, the inverse of the majorizer \ref{eq:major} to each voxel update and the Jacobi preconditioner \ref{eq:jacobi}.
Several regularizers were also tested: Total Variation (TV), L1, 0th and 1st order Tikhonov (L2) were tested for different $\lambda$ values across the five models. 
\item Algebraic Reconstruction Techniques (ART, SART, and SIRT) for different values of the relaxation parameter, $\lambda$.

\item Conjugate Gradient, LSQR per-view and per-epoch (all views) update strategies, and GMRES laso per-view and per-epoch (all views) update strategies. These algorithms require many internal iterations and carry high computational overhead. 

\end{itemize}

\subsubsection{Method parameterization}
The parameterization of the methods includes the value of learning rate $\alpha$, the relaxation factor of the regularizers $\lambda$ and parameters that control vessel segmentation which affect IoU. To maintain a fair comparison, we used model 0148 to obtain parameter values with reference to ground truth that these parameters were applied to the remaining 4 3D models.

\begin{figure}[h]
    \centering
    \includegraphics[width=1\linewidth]{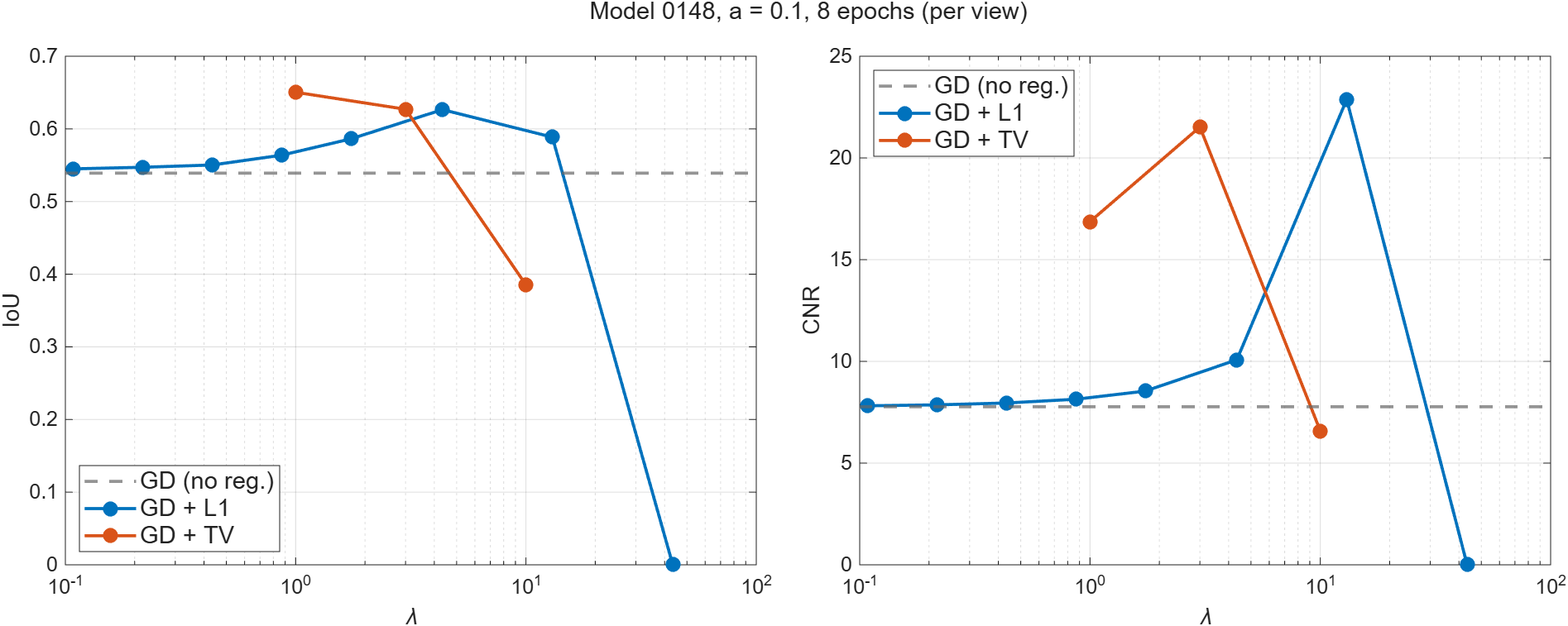}
    \caption{Evaluation of the regularization parameter $\lambda$ on model 0148 for gradient descent. The plots display the Intersection over Union (IoU) and Contrast-to-Noise Ratio (CNR) across various $\lambda$ values for L1 and Total Variation (TV) regularization to determine the optimal parameters applied to the remaining models.}
    \label{fig:best_lambda}
\end{figure}

\paragraph{\textbf{Gradient Descend variants}}
Figure \ref{fig:best_lambda} illustrates the evaluation of the ideal $\lambda$ parameter for gradient descent. 
The parameter tuning was initially conducted on model 0148, and the resulting optimal values were subsequently applied to the remaining models. The plots demonstrate the impact of varying $\lambda$ on the Intersection over Union (IoU) and Contrast-to-Noise Ratio (CNR). 
For L1 regularization, $\lambda$ was tested across the values [0.081, 0.16, 0.32, 0.65, 1.3, 3.2, 9.7, 32], yielding the best IoU at $\lambda = 3.2$ and the highest CNR at $\lambda = 9.7$. Similarly, for Total Variation (TV) regularization, tested $\lambda$ values included [1, 3, 10], which resulted in the best IoU being achieved with $\lambda = 0.1$ and the best CNR with $\lambda = 0.3$.

\begin{figure}[H]
    \centering
    \includegraphics[width=0.85\linewidth]{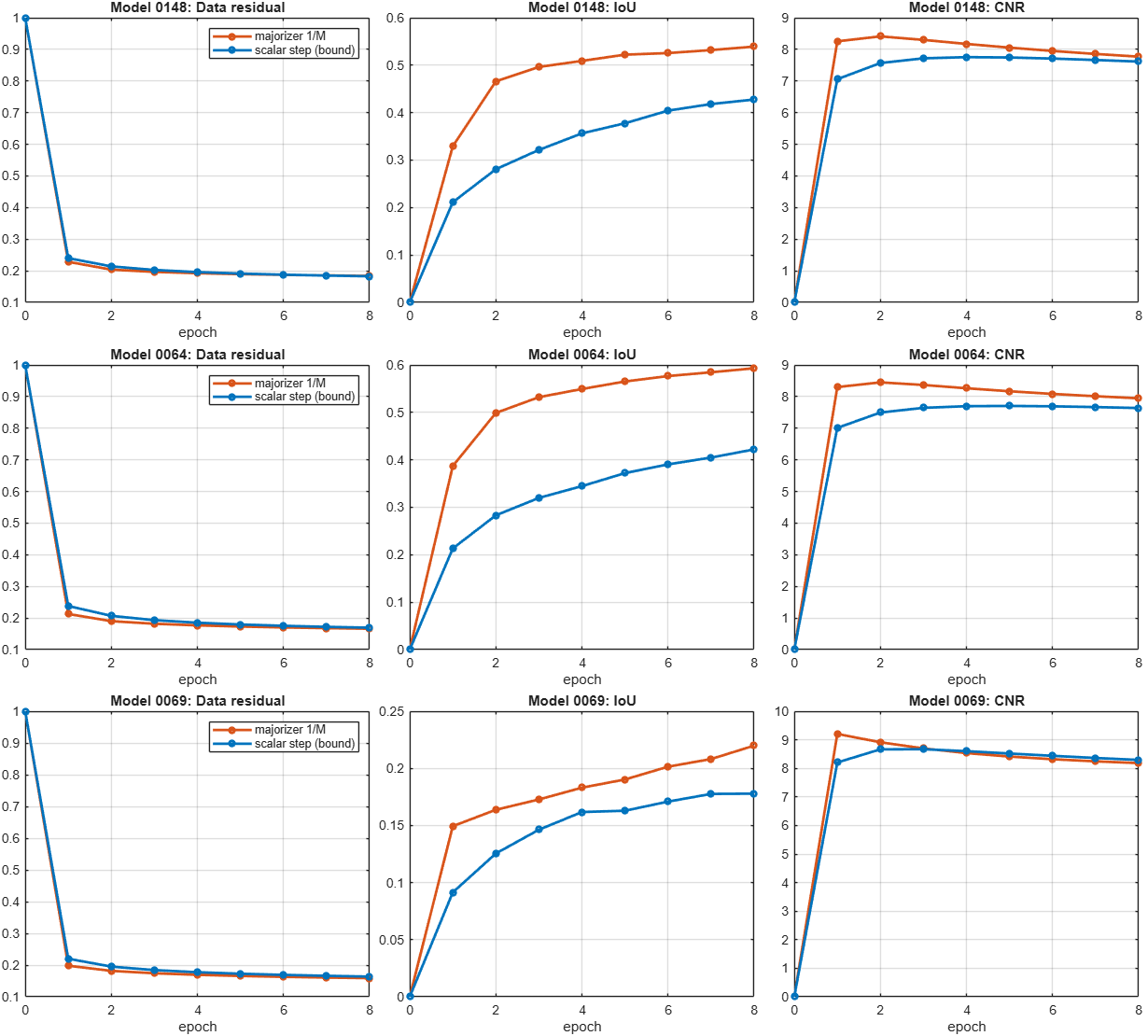}
    \caption{Evaluation of step size selection in unregularized gradient descent reconstructions for models 0064, 0069, and 0148 (rows). The columns track the residual, IoU, and CNR across epochs. The majorizer inverse step size, depicted in the orange curve, outperforms our Hessian approximation, blue curve, in both IoU and CNR, despite converging to a comparable residual.}
    \label{fig:bound_vs_maj}
\end{figure}

For the Gradient Descent framework, we also investigated the Lipschitz gradient step size \ref{eq:lr_maj} versus the inverse of the majorizer \ref{eq:major} to each voxel update. 
Figure \ref{fig:bound_vs_maj} illustrates the volume reconstructions for models 0064, 0069, and 0148 using unregularized gradient descent. The standard learning rate \ref{eq:standard_lr} provided stable updates. 
The Jacobi preconditioner \ref{eq:jacobi} did not reach a solution within the $8$-epoch limit due to slow convergence rates. While both step sizes minimize the residual similarly, the majorizer approach yields better Intersection over Union (IoU) and Contrast-to-Noise Ratio (CNR), as well as faster intermediate convergence. 
Consequently, the remainder of the experiments were performed using the majorizer step size.

\paragraph{\textbf{Algebraic Reconstruction Techniques}}
We evaluated the Algebraic Reconstruction Techniques (ART, SART, and SIRT) for $8$ epochs.
For standard ART, the ray-by-ray update strategy is very sensitive to the relaxation parameter, $\lambda$, with marked error metric deterioration as  $\lambda $ increases from $1.0$, to $\lambda = 2.0$. For instance, at $\lambda = 2.0$, IoU dropped from $0.573$ to $0.295$, CNR and skeleton exhibited similar behavior.
SART improved error metrics with higher relaxation parameters. It achieved its best reconstruction at $\lambda = 2.0$ with an IoU of $0.543$ and a Chamfer distance of $0.5756$. It completed the process in $63$ seconds.
SIRT also requires a higher relaxation parameter, $\lambda = 2.0$. However SIRT showed a performance inconsistency between IoU (middle performance), CNR (2nd best) and SSIM (worst performance) across all five models. 

Overall, the per-view update strategy of SART offered a better balance of geometric accuracy and computational speed.

\paragraph{\textbf{Krylov subspace solvers}}
The tested Krylov subspace solvers (Conjugate Gradient-GD-, LSQR, and GMRES) require many internal iterations and carry high computational overhead. Therefore we evaluated them over $4$ global epochs.

The CG method was applied using an all-views (per-epoch) update scheme. Despite executing 44 internal passes per epoch, it failed to recover the underlying vessel geometry effectively, yielding IoU of only $0.218$. Furthermore, it suffered from geometric distortion, evidenced by a high Chamfer distance of $3.8283$ and the highest Skeleton Center Error recorded in our experiments ($5.4314$). 

LSQR was tested using both per-view and per-epoch update , allocating 24 passes per epoch. While LSQR performed marginally better than CG, the computational cost was high, requiring roughly $716$ seconds to complete. The per-view formulation achieved a slightly better volumetric overlap (IoU of $0.326$) compared to the all-views approach (IoU of $0.305$). However, the all-views approach proved much more effective at reducing outlier noise, dropping the Chamfer distance to $1.8234$ (down from $3.0268$ in the per-view setup) and improving the CNR to $6.4$. Even with these slight improvements, LSQR failed to approach the topological fidelity achieved by the simpler, regularized gradient methods.

Finally, GMRES was evaluated for 44 internal passes, proving to be by far the most computationally demanding scheme tested. Running GMRES took between 1370 and $1397$ seconds, however it does not produce meaningful results. The per-view update yielded a poor IoU of $0.271$ and a degraded Chamfer distance of $4.0051$. Similarly to LSQR, shifting GMRES to a per-epoch update smoothed out the worst of the spatial artifacts—improving the Chamfer distance to $1.4441$ and the CNR to $7.1$, but the overall results remained completely poor with an IoU of $0.290$.

Overall we obtained results from 11 method parameterizations. Table \ref{tab:ranks-5-iou} shows the achieved IOU from each method parameterization and the ranking of the method (1 being the best) for each one of the 3D models. Tables \ref{tab:ranks-5-cnr} and \ref{tab:ranks-5-ssim} depict the same results for CNR and SSIM on the held-out view.
Table~\ref{tab:reconstruction_schemes} summarizes the rankings from the previous three tables by reporting mean values across the five 3D models for IoU, SSIM, Chamfer distance, skeleton centerline and diameter error. The volume noise parameter was held constant at $\nu=0.1$ during all these experiments. All GD variants, as well as ART, SART and SIRT were evaluated over 8 epochs. In to keep the comparisons fair, for the remaining metheds that have internal iterations we reported the achieved metrics over the number of passes (forward and back-projections). Specifically, GD and GMRES were allowed to evolve for 4 epochs with 11 internal passes per epoch. Also LSQR converged for 4 epochs with 6 passes per epoch.
Time measurements were measured on a Lenovo laptop with an AMD Ryzen AI 7 450 CPU (8 cores, 16 threads), 32 GB RAM and Windows 11, in MATLAB R2026a running single-process on the CPU with implicit multithreading (8 threads).


As summarized in Table \ref{tab:reconstruction_schemes}, Gradient Descent regularized with Total Variation (GD + TV) and GD + L1 proved the most effective approach, in terms of both CNR and SSIM. Other techniques exhibited distinct performance trade-offs. For example, SIRT ($\lambda=2$) ranked second for noise suppression (CNR) but eleventh for held-out projection structural preservation (SSIM). Finally, the Krylov subspace solvers occupied the lowest ranks across all image quality and segmentation metrics despite their computational overhead.
Both regularization schemes suppressed artifacts and maintained the topology of the volumes.

\begin{table}[H]
  \centering
  \footnotesize
  \caption{IoU and rank (1 = best) of each method on the five models ($a = 0.1$), sorted by the mean rank. Parameters selected on model 0148 and used unchanged on the other models. On models 0069 and 0088 the IoU of the unregularized methods is limited by the segmentation threshold, which includes background noise.}
  \label{tab:ranks-5-iou}
  \setlength{\tabcolsep}{3pt}
  \begin{tabular}{lccccccccccc}
    \toprule
    & \multicolumn{2}{c}{Model 0064} & \multicolumn{2}{c}{Model 0069} & \multicolumn{2}{c}{Model 0088} & \multicolumn{2}{c}{Model 0109} & \multicolumn{2}{c}{Model 0148} & \\
    \cmidrule(lr){2-3}\cmidrule(lr){4-5}\cmidrule(lr){6-7}\cmidrule(lr){8-9}\cmidrule(lr){10-11}
    Method & IoU & Rank & IoU & Rank & IoU & Rank & IoU & Rank & IoU & Rank & Mean rank \\
    \midrule
    GD + TV ($\lambda = 0.1$) & 0.673 & 1 & 0.666 & 1 & 0.565 & 2 & 0.766 & 1 & 0.650 & 1 & 1.20 \\
    GD + $\text{L}_1$ ($\lambda = 3.2$) & 0.661 & 2 & 0.396 & 4 & 0.424 & 3 & 0.729 & 2 & 0.627 & 2 & 2.60 \\
    ART ($\lambda=1$) & 0.650 & 3 & 0.505 & 2 & 0.571 & 1 & 0.675 & 4 & 0.563 & 4 & 2.80 \\
    SART ($\lambda=2$) & 0.647 & 4 & 0.381 & 5 & 0.389 & 4 & 0.723 & 3 & 0.601 & 3 & 3.80 \\
    Nes. GD + TV ($\beta=0.9$) & 0.552 & 6 & 0.473 & 3 & 0.330 & 6 & 0.577 & 5 & 0.545 & 5 & 5.00 \\
    SIRT ($\lambda=2$) & 0.514 & 7 & 0.309 & 6 & 0.365 & 5 & 0.525 & 6 & 0.513 & 7 & 6.20 \\
    GD no reg. & 0.593 & 5 & 0.220 & 7 & 0.238 & 7 & 0.449 & 7 & 0.539 & 6 & 6.40 \\
    LSQR (per view) & 0.458 & 8 & 0.156 & 8 & 0.165 & 9 & 0.421 & 8 & 0.364 & 10 & 8.60 \\
    LSQR (per epoch) & 0.395 & 9 & 0.138 & 10 & 0.161 & 10 & 0.394 & 9 & 0.383 & 8 & 9.20 \\
    GMRES (per epoch) & 0.349 & 11 & 0.148 & 9 & 0.169 & 8 & 0.325 & 10 & 0.373 & 9 & 9.40 \\
    GMRES (per view) & 0.390 & 10 & 0.132 & 11 & 0.147 & 11 & 0.321 & 11 & 0.291 & 11 & 10.80 \\
    CG (per epoch) & 0.324 & 12 & 0.063 & 12 & 0.089 & 12 & 0.273 & 12 & 0.267 & 12 & 12.00 \\
    \bottomrule
  \end{tabular}
\end{table}

\begin{table}[H]
  \centering
  \footnotesize
  \caption{CNR and rank (1 = best) of each method on the five models ($a = 0.1$), sorted by the mean rank. Parameters selected on model 0148 and used unchanged on the other models.}
  \label{tab:ranks-5-cnr}
  \setlength{\tabcolsep}{3pt}
  \begin{tabular}{lccccccccccc}
    \toprule
    & \multicolumn{2}{c}{Model 0064} & \multicolumn{2}{c}{Model 0069} & \multicolumn{2}{c}{Model 0088} & \multicolumn{2}{c}{Model 0109} & \multicolumn{2}{c}{Model 0148} & \\
    \cmidrule(lr){2-3}\cmidrule(lr){4-5}\cmidrule(lr){6-7}\cmidrule(lr){8-9}\cmidrule(lr){10-11}
    Method & CNR & Rank & CNR & Rank & CNR & Rank & CNR & Rank & CNR & Rank & Mean rank \\
    \midrule
    GD + TV ($\lambda = 0.1$) & 19.4 & 1 & 30.3 & 1 & 11.9 & 1 & 26.7 & 1 & 16.9 & 1 & 1.00 \\
    Nes. GD + TV ($\beta=0.9$) & 12.3 & 2 & 18.2 & 2 & 9.4 & 3 & 15.7 & 2 & 11.9 & 2 & 2.20 \\
    SIRT ($\lambda=2$) & 10.8 & 3 & 12.2 & 3 & 10.5 & 2 & 12.7 & 3 & 11.2 & 3 & 2.80 \\
    GD + $\text{L}_1$ ($\lambda = 3.2$) & 10.4 & 4 & 10.9 & 4 & 9.0 & 4 & 12.1 & 4 & 10.1 & 4 & 4.00 \\
    SART ($\lambda=2$) & 8.8 & 5 & 9.4 & 5 & 7.5 & 5 & 10.1 & 5 & 8.3 & 5 & 5.00 \\
    GD no reg. & 7.9 & 6 & 8.2 & 7 & 6.8 & 6 & 9.0 & 6 & 7.8 & 6 & 6.20 \\
    GMRES (per epoch) & 6.8 & 8 & 7.5 & 8 & 6.2 & 7 & 8.0 & 7 & 7.1 & 7 & 7.40 \\
    ART ($\lambda=1$) & 7.5 & 7 & 8.7 & 6 & 5.8 & 8 & 3.3 & 12 & 6.7 & 8 & 8.20 \\
    LSQR (per epoch) & 6.2 & 9 & 6.8 & 9 & 5.4 & 9 & 7.2 & 8 & 6.2 & 9 & 8.80 \\
    LSQR (per view) & 5.4 & 10 & 6.4 & 10 & 4.7 & 10 & 6.1 & 9 & 4.8 & 10 & 9.80 \\
    GMRES (per view) & 4.9 & 11 & 5.9 & 11 & 4.2 & 11 & 5.5 & 10 & 4.3 & 12 & 11.00 \\
    CG (per epoch) & 4.7 & 12 & 5.4 & 12 & 4.1 & 12 & 5.3 & 11 & 4.4 & 11 & 11.60 \\
    \bottomrule
  \end{tabular}
\end{table}

\begin{table}[H]
  \centering
  \footnotesize
  \caption{Held-out view SSIM and rank (1 = best) of each method on the five models ($a = 0.1$), sorted by the mean rank. Parameters selected on model 0148 and used unchanged on the other models.}
  \label{tab:ranks-5-ssim}
  \setlength{\tabcolsep}{3pt}
  \begin{tabular}{lccccccccccc}
    \toprule
    & \multicolumn{2}{c}{Model 0064} & \multicolumn{2}{c}{Model 0069} & \multicolumn{2}{c}{Model 0088} & \multicolumn{2}{c}{Model 0109} & \multicolumn{2}{c}{Model 0148} & \\
    \cmidrule(lr){2-3}\cmidrule(lr){4-5}\cmidrule(lr){6-7}\cmidrule(lr){8-9}\cmidrule(lr){10-11}
    Method & SSIM & Rank & SSIM & Rank & SSIM & Rank & SSIM & Rank & SSIM & Rank & Mean rank \\
    \midrule
    GD + TV ($\lambda = 0.1$) & 0.758 & 1 & 0.682 & 1 & 0.800 & 1 & 0.857 & 1 & 0.897 & 1 & 1.00 \\
    Nes. GD + TV ($\beta=0.9$) & 0.724 & 2 & 0.644 & 3 & 0.781 & 2 & 0.831 & 2 & 0.883 & 2 & 2.20 \\
    GD no reg. & 0.713 & 3 & 0.649 & 2 & 0.780 & 3 & 0.817 & 3 & 0.878 & 3 & 2.80 \\
    SART ($\lambda=2$) & 0.702 & 4 & 0.638 & 4 & 0.769 & 4 & 0.807 & 4 & 0.864 & 4 & 4.00 \\
    ART ($\lambda=1$) & 0.686 & 5 & 0.612 & 5 & 0.744 & 5 & 0.787 & 5 & 0.860 & 5 & 5.00 \\
    LSQR (per epoch) & 0.674 & 6 & 0.611 & 6 & 0.744 & 5 & 0.785 & 6 & 0.859 & 6 & 5.80 \\
    LSQR (per view) & 0.640 & 7 & 0.557 & 7 & 0.671 & 7 & 0.742 & 7 & 0.825 & 7 & 7.00 \\
    GMRES (per view) & 0.615 & 8 & 0.531 & 8 & 0.640 & 9 & 0.719 & 8 & 0.806 & 8 & 8.20 \\
    CG (per epoch) & 0.590 & 10 & 0.524 & 9 & 0.639 & 10 & 0.696 & 9 & 0.777 & 9 & 9.40 \\
    GD + $\text{L}_1$ ($\lambda = 3.2$) & 0.611 & 9 & 0.508 & 10 & 0.661 & 8 & 0.692 & 10 & 0.750 & 10 & 9.40 \\
    GMRES (per epoch) & 0.472 & 11 & 0.394 & 11 & 0.505 & 11 & 0.545 & 11 & 0.595 & 11 & 11.00 \\
    SIRT ($\lambda=2$) & 0.428 & 12 & 0.325 & 12 & 0.401 & 12 & 0.488 & 12 & 0.534 & 12 & 12.00 \\
    \bottomrule
  \end{tabular}
\end{table}

\begin{sidewaystable}[!htbp]
\centering
\caption{\textbf{All reconstruction schemes, $a = 0.1$}}
\label{tab:reconstruction_schemes}
\vspace{0.5em}
\small
\setlength{\tabcolsep}{4.5pt}
\begin{tabular}{l l c c c c c c c c c}
\toprule
\textbf{Method} & \textbf{Update} & \textbf{Epochs} & \textbf{Passes} & \textbf{IoU} & \textbf{SSIM} & \textbf{Chamfer Dist.} & \textbf{Skel Center Err.} & \textbf{Skel Diam Err.} & \textbf{CNR} & \textbf{Time (s)} \\
\midrule
GD no reg. & perView & 8 & 8 & 0.451 & 0.767 & 0.8242 & 1.9222 & 0.6279 & 8.0 & 48 \\
GD + $\text{L}_1$ ($\lambda = 3.2$) & perView & 8 & 8 & 0.649 & 0.644 & 0.3113 & 1.6109 & 0.6451 & 17.0 & 119 \\
GD + TV ($\lambda = 0.1$) & perView & 8 & 8 & 0.648 & 0.799 & 0.3732 & 1.3105 & 0.3413 & 21.0 & 48 \\
Nes. GD + TV ($\beta =0.9$) & perEpoch & 8 & 8 & 0.495 & 0.773 & 1.3750 & 1.0460 & 0.4840 & 13.5 & 240 \\
\addlinespace
ART ($\lambda = 0.5$) & perView & 8 & 8 & 0.500 & 0.750 & 0.8581 & 1.7388 & 0.7034 & 7.3 & 97 \\
ART ($\lambda = 1.0$) & perView & 8 & 8 & 0.573 & 0.738 & 0.8047 & 1.2583 & 0.7764 & 7.6 & 97 \\
ART ($\lambda = 1.5$) & perView & 8 & 8 & 0.560 & 0.724 & 1.0503 & 1.1488 & 0.8424 & 7.6 & 91 \\
ART ($\lambda = 2.0$) & perView & 8 & 8 & 0.295 & 0.688 & 3.8202 & 3.7578 & 0.8410 & 5.6 & 88 \\
\addlinespace
SART ($\lambda = 0.5$) & perView & 8 & 8 & 0.390 & 0.770 & 1.0021 & 2.2234 & 0.5906 & 8.1 & 63 \\
SART ($\lambda = 1.0$) & perView & 8 & 8 & 0.451 & 0.767 & 0.8258 & 1.9185 & 0.6291 & 8.0 & 63 \\
SART ($\lambda = 1.5$) & perView & 8 & 8 & 0.493 & 0.763 & 0.7238 & 1.6022 & 0.6438 & 8.1 & 63 \\
SART ($\lambda = 2.0$) & perView & 8 & 8 & 0.543 & 0.756 & 0.5756 & 1.3051 & 0.6443 & 8.8 & 63 \\
\addlinespace
SIRT ($\lambda = 0.5$) & perEpoch & 8 & 8 & 0.186 & 0.745 & 1.9452 & 3.9981 & 1.0172 & 7.9 & 97 \\
SIRT ($\lambda = 1.0$) & perEpoch & 8 & 8 & 0.253 & 0.757 & 1.4319 & 2.8539 & 0.8018 & 8.2 & 94 \\
SIRT ($\lambda = 1.5$) & perEpoch & 8 & 8 & 0.286 & 0.762 & 1.2852 & 2.7373 & 0.7376 & 8.3 & 93 \\
SIRT ($\lambda = 2.0$) & perEpoch & 8 & 8 & 0.445 & 0.435 & 0.5782 & 1.3161 & 0.5327 & 11.4 & 92 \\
\addlinespace
CG (per epoch) & perEpoch & 4 & 44 & 0.218 & 0.645 & 3.8283 & 5.4314 & 0.7058 & 4.8 & 151 \\
LSQR (per view) & perView & 4 & 24 & 0.326 & 0.687 & 3.0268 & 2.9987 & 0.7696 & 5.6 & 716 \\
LSQR (per epoch) & perEpoch & 4 & 24 & 0.305 & 0.735 & 1.8234 & 3.3830 & 0.6465 & 6.4 & 716 \\
GMRES (per view) & perView & 4 & 44 & 0.271 & 0.662 & 4.0051 & 3.3872 & 0.8161 & 5.0 & 1397 \\
GMRES (per epoch) & perEpoch & 4 & 44 & 0.290 & 0.502 & 1.4441 & 3.2574 & 0.6632 & 7.1 & 1370 \\
\bottomrule
\end{tabular}
\end{sidewaystable}

The quantitative results across the five test models are detailed in 
Table \ref{tab:rank-agreement-5} shows the correlation between method ranking for the 11 parameterized methods for IoU, CNR, SSIM and NRMSE error metrics.  The method rankings demonstrated high stability regardless of the specific vessel geometry, yielding Kendall's $W$ scores above $0.919$ for all metrics. 
\begin{table}[H]
  \centering
  \caption{Agreement of the method rankings between the five models ($n = 13$ methods, $a = 0.1$): mean and range of the pairwise Spearman $\rho$ and Kendall $\tau_b$ over the 10 model pairs, and Kendall's $W$ for the five rankings jointly. IoU depends on the segmentation threshold; the other metrics do not. SSIM and NRMSE are computed on the held-out view.}
  \label{tab:rank-agreement-5}
  \begin{tabular}{lccc}
    \toprule
    Metric & Spearman $\rho$, mean (range) & Kendall $\tau_b$, mean (range) & Kendall's $W$ \\
    \midrule
    IoU (threshold) & 0.953 (0.918-0.978) & 0.856 (0.795-0.923) & 0.962 \\
    CNR & 0.953 (0.885-0.995) & 0.892 (0.821-0.974) & 0.963 \\
    SSIM & 0.977 (0.940-0.995) & 0.933 (0.872-0.974) & 0.982 \\
    NRMSE & 0.898 (0.808-0.995) & 0.795 (0.641-0.974) & 0.919 \\
    \bottomrule
  \end{tabular}
\end{table}

The convergence behavior of the evaluated methods is detailed in Figure \ref{fig:all_conv}. Each sub-figure contains three plots tracking the  residual, IoU, and CNR against the number of projector passes. GD and algebraic techniques ran for $8$ epochs, while both LSQR variants ran for $24$ passes, and CG and GMRES solvers ran for $44$ total passes.

For the residuals, the regularized GD algorithms, ART, and SART show a rapid drop during early iterations, reducing normalized residual error below $0.25$ within the first few epochs. SIRT exhibits a slower residual drop due to its per-epoch update scheme.  Both LSQR variants decrease similarly fast within $5$ passes, but slower than GD, ART and SART. Among the other Krylov solvers, GMRES per-view requires over $10$ passes to drop close to $0.20$ and then remains constant, while the per-epoch variant decreases very slowly throughout $44$ passes. CG initially reduces the residual to $\approx 0.20-0.25$ at pass $11$ but subsequently diverges, rising above $0.40$.   

As for IoU, GD + TV, GD + L1, ART ($\lambda=1$), and SART ($\lambda=2$) display steep upward IoU trajectories in the initial passes. GD regularized with TV achieved the highest final IoU. SIRT displays non-monotonic progression across epochs. The Krylov subspace methods converge much slower across their passes, plateauing under $0.40$.   

CNR curves highlight clear separation among the regularization strategies. GD with TV shows a continuous upward curve, reaching the highest CNR scores. SIRT and GD with L1 reach moderate CNR values and unregularized GD and standard algebraic techniques achieve lower CNR. The Krylov subspace solvers exhibit weak noise suppression overall. LSQR and GMRES variants reach a plateau below CNR of $8$, while CG peaks at pass $11$ before steadily dropping to a CNR of $\approx 5$. 

\begin{figure}[H]
    \centering
    \begin{subfigure}[b]{1\textwidth}
        \centering
        \includegraphics[width=1\linewidth]{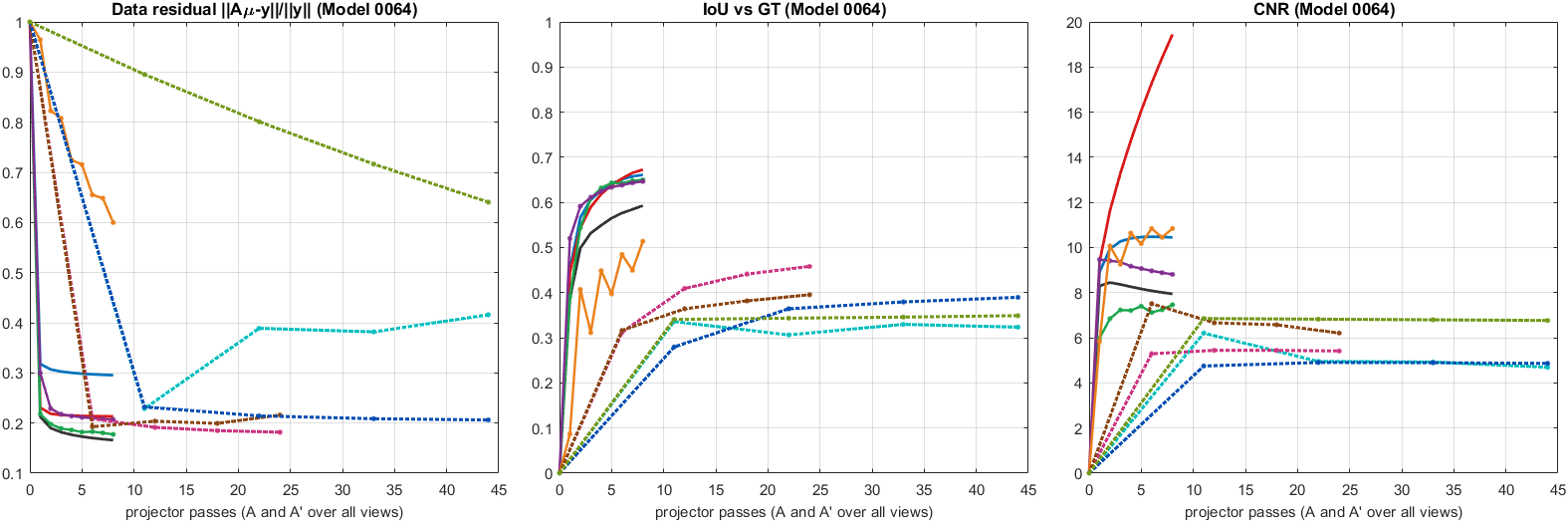}
        \caption{Model $0064$}
        \label{fig:0064_conv}
    \end{subfigure}
    
    \vspace{1em} 
    
    \begin{subfigure}[b]{1\textwidth}
        \centering
        \includegraphics[width=1\linewidth]{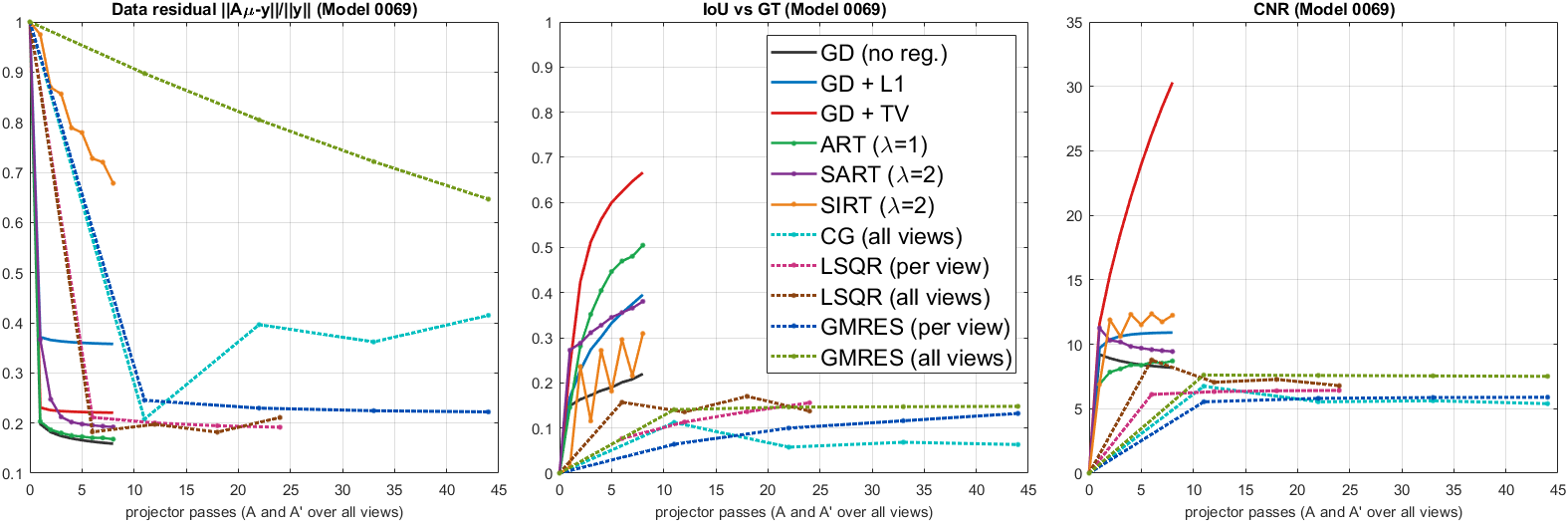}
        \caption{Model $0069$}
        \label{fig:0069_conv}
    \end{subfigure}
    
    \vspace{1em}
    
    \caption{Convergence behavior of the different methods for the five models. Displayed left to right: (i) Data Residual, (ii) IoU, (iii) CNR. (Continued on next page)}
    \label{fig:all_conv_part1}
\end{figure}

\begin{figure}[H]
    \ContinuedFloat 
    \centering
        \begin{subfigure}[b]{1\textwidth}
        \centering
        \includegraphics[width=1\linewidth]{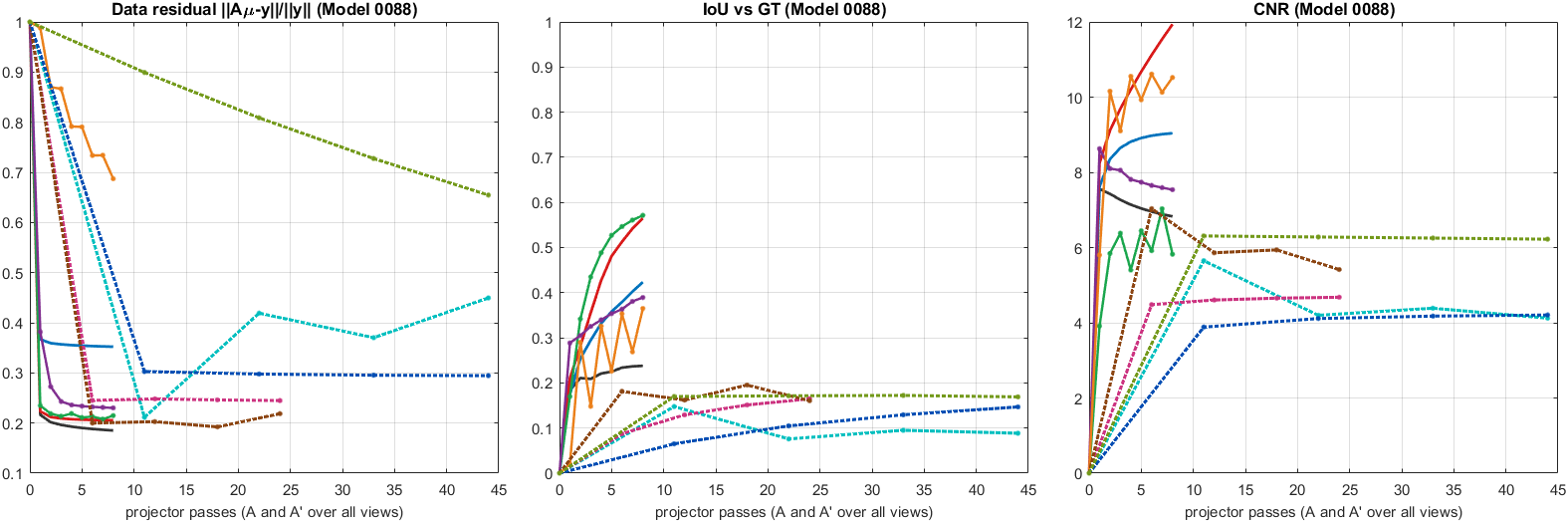}
        \caption{Model $0088$}
        \label{fig:0088_conv}
    \end{subfigure}
    
    \begin{subfigure}[b]{1\textwidth}
        \centering
        \includegraphics[width=1\linewidth]{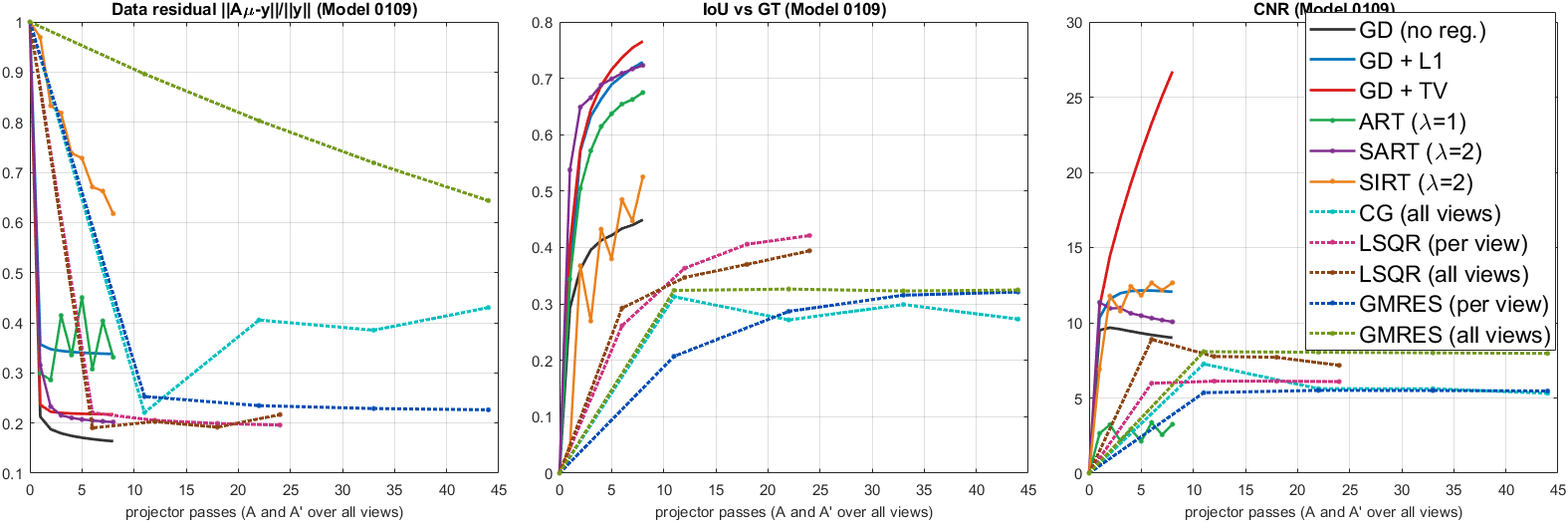}
        \caption{Model $0109$}
        \label{fig:0109_conv}
    \end{subfigure}
    
    \vspace{1em}
    
    \begin{subfigure}[b]{1\textwidth}
        \centering
        \includegraphics[width=1\linewidth]{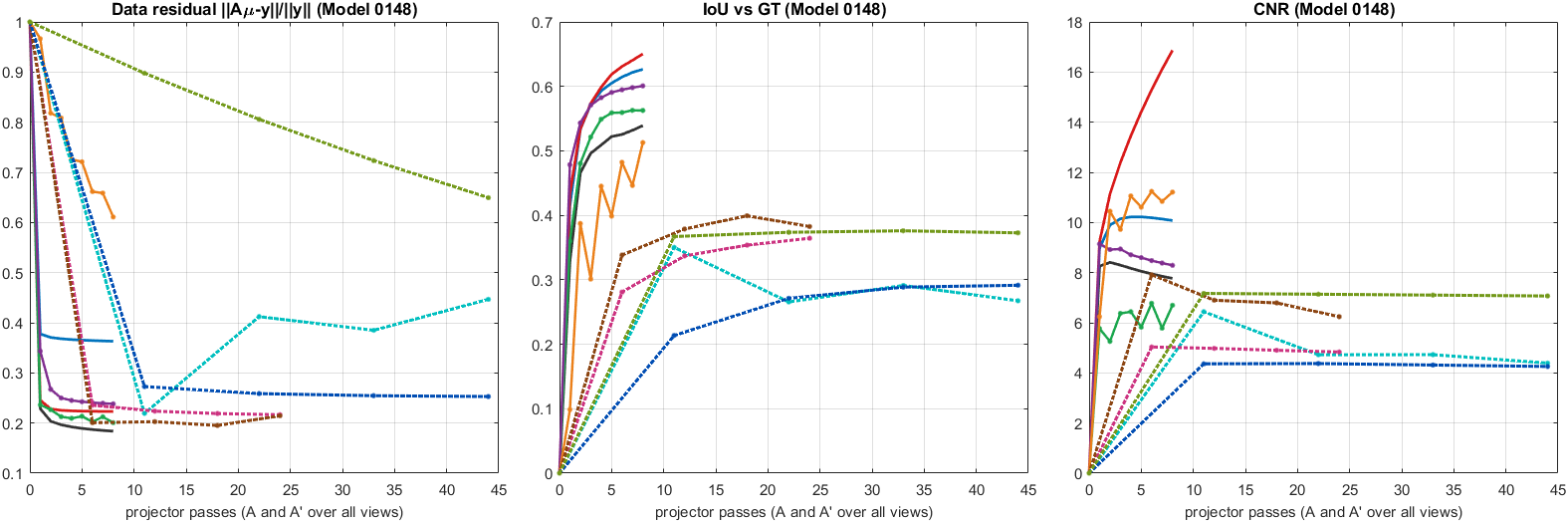}
        \caption{Model $0148$}
        \label{fig:0148_conv}
    \end{subfigure}
    
    \caption{Convergence behavior of the different methods for the five models (continued).}
    \label{fig:all_conv}
\end{figure}

\begin{figure}
    \centering
    \includegraphics[width=1\linewidth]{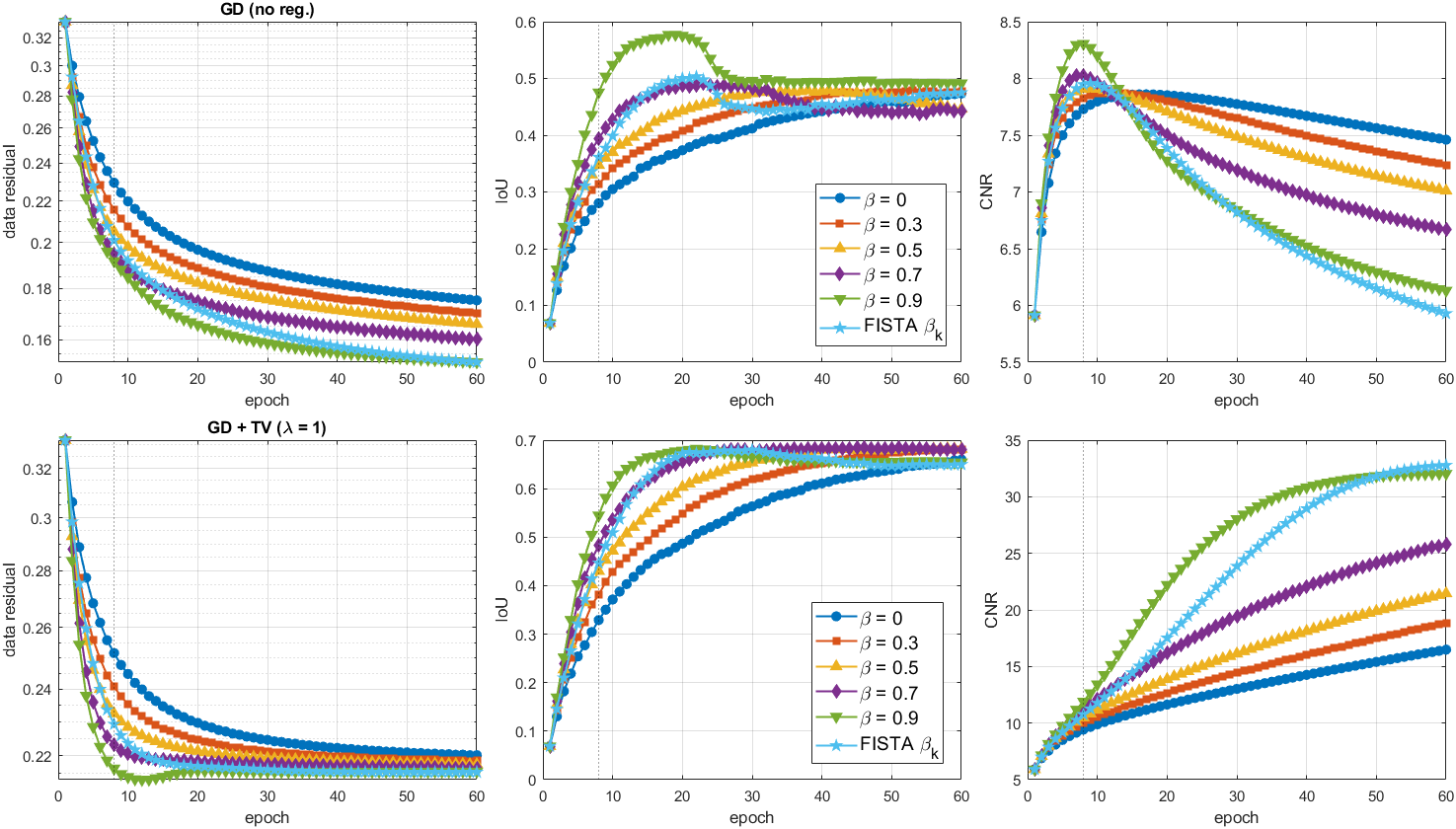}
    \caption{Data residual (left), IoU (middle), and CNR (right) for unregularized (top) versus TV-regularized (bottom) gradient descent with Nesterov momentum parameter $\beta \in \{0, 0.3, 0.5, 0.7, 0.9\}$.}
    \label{fig:nesterov}
\end{figure}

Figure \ref{fig:nesterov} shows the evolution of data residual, IoU and CNR for model 0148 for GD (with majorizer) with TV ($\lambda=1$ for several values of the momentum coefficient $\beta$. Since Nesterov's accelerated GD updates per epoch, $\beta=0$ corresponds to GD per epoch. The second row shows the same graphs for the unregularized GD (with majorizer). The optimal value of $\beta$ is within [0.5 - 0.8] range. 

\begin{figure}[H]
    \centering
    \begin{subfigure}[b]{0.85\textwidth}
        \centering
        \includegraphics[width=\linewidth]{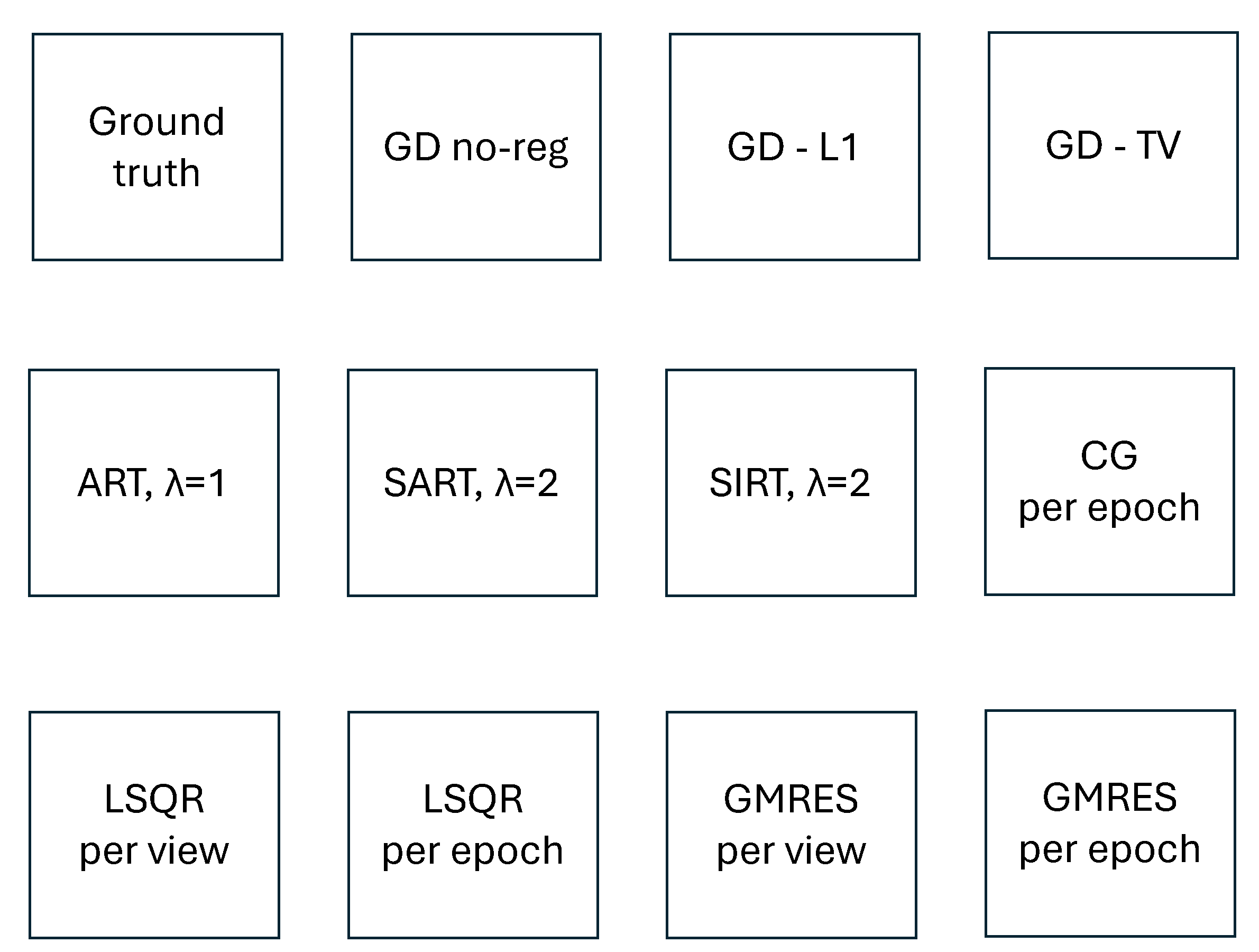}
        \caption{Legend mapping all the figures to the corresponding reconstruction technique}
        \label{fig:leg}
    \end{subfigure}
    \hfill
    \begin{subfigure}[b]{0.85\textwidth}
        \centering
        \includegraphics[width=\linewidth]{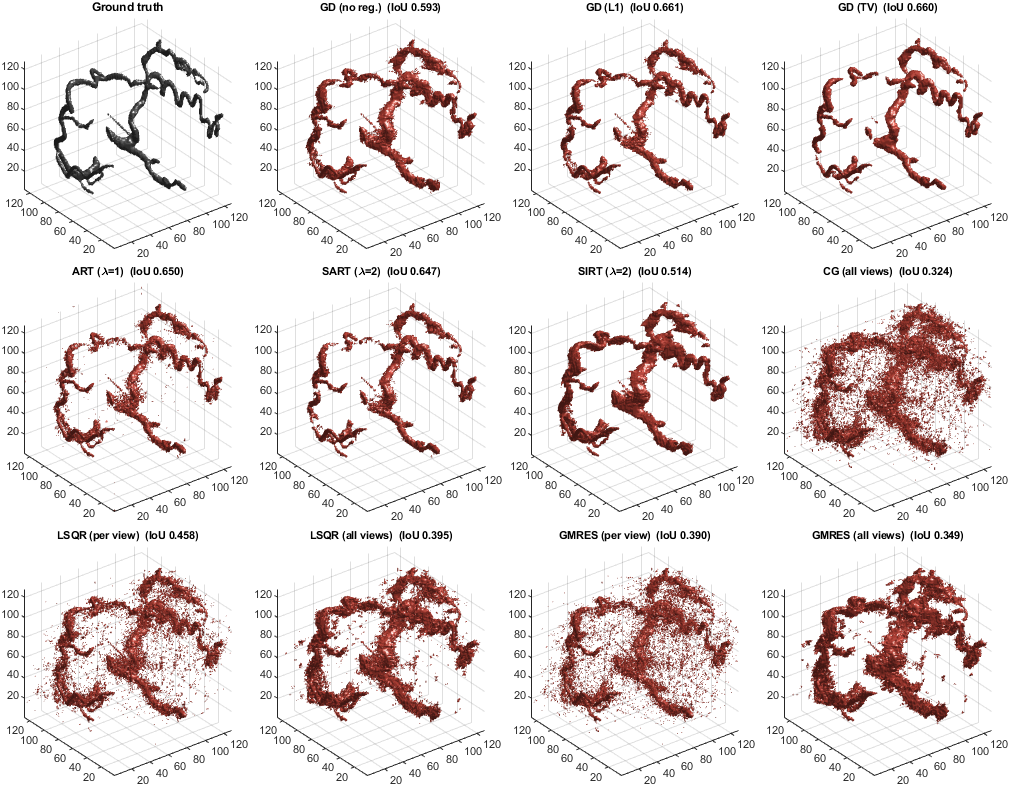}
        \caption{Model $0064$}
        \label{fig:0064_recon}
    \end{subfigure}
    \caption{Reconstruction analysis: Overview and Model $0064$.}
    \label{fig:group1_recon}
\end{figure}

\begin{figure}[H]
    \centering
    \begin{subfigure}[b]{0.85\textwidth}
        \centering
        \includegraphics[width=\linewidth]{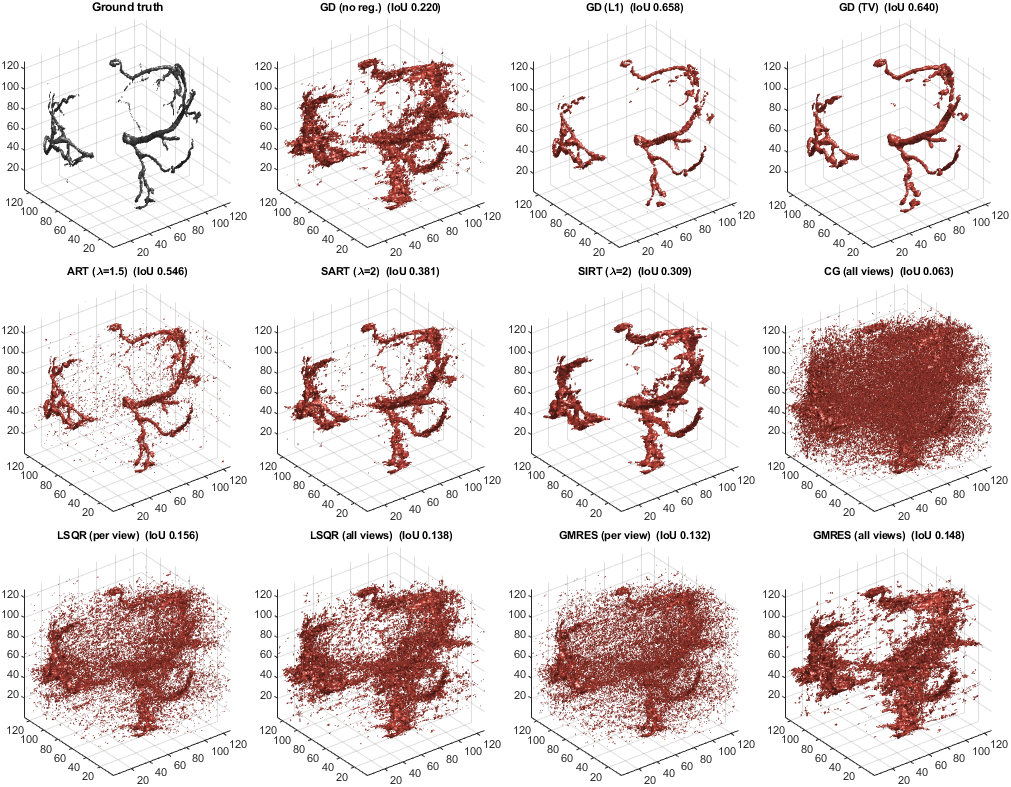}
        \caption{Model $0069$}
        \label{fig:0069_recon}
    \end{subfigure}
    \hfill
    \begin{subfigure}[b]{0.85\textwidth}
        \centering
        \includegraphics[width=\linewidth]{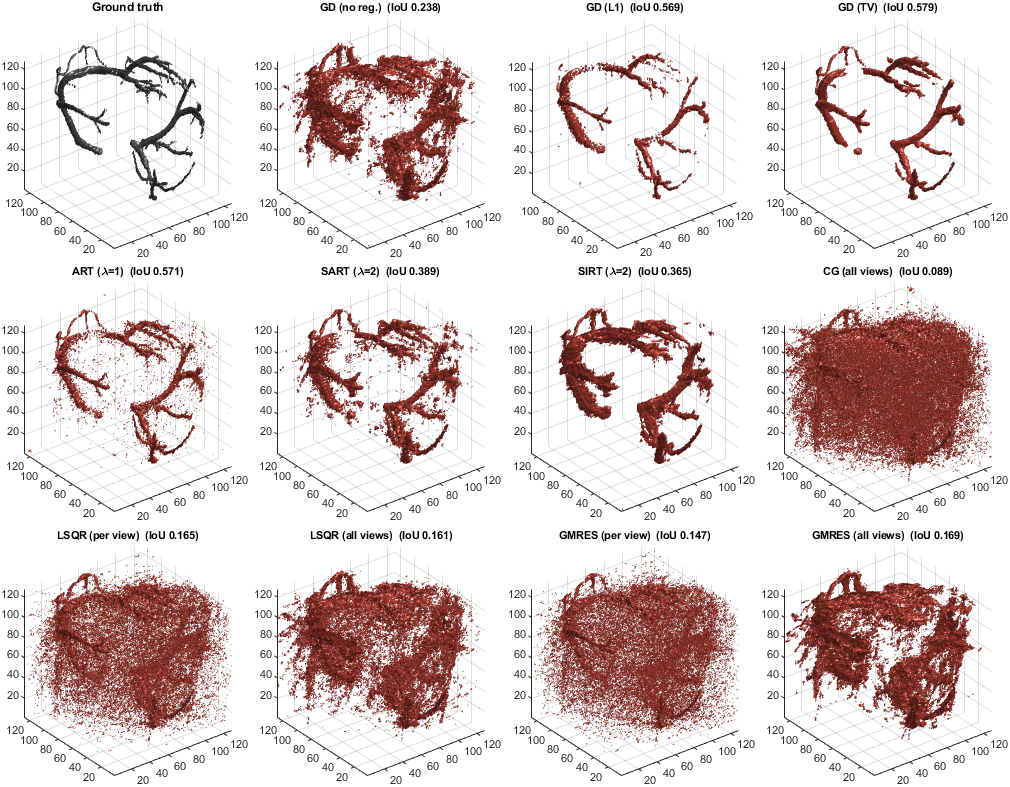}
        \caption{Model $0088$}
        \label{fig:0088_recon}
    \end{subfigure}
    \caption{Reconstruction results for models $0069$ and $0088$.}
    \label{fig:group2_recon}
\end{figure}

\begin{figure}[H]
    \centering
    \begin{subfigure}[b]{0.85\textwidth}
        \centering
        \includegraphics[width=\linewidth]{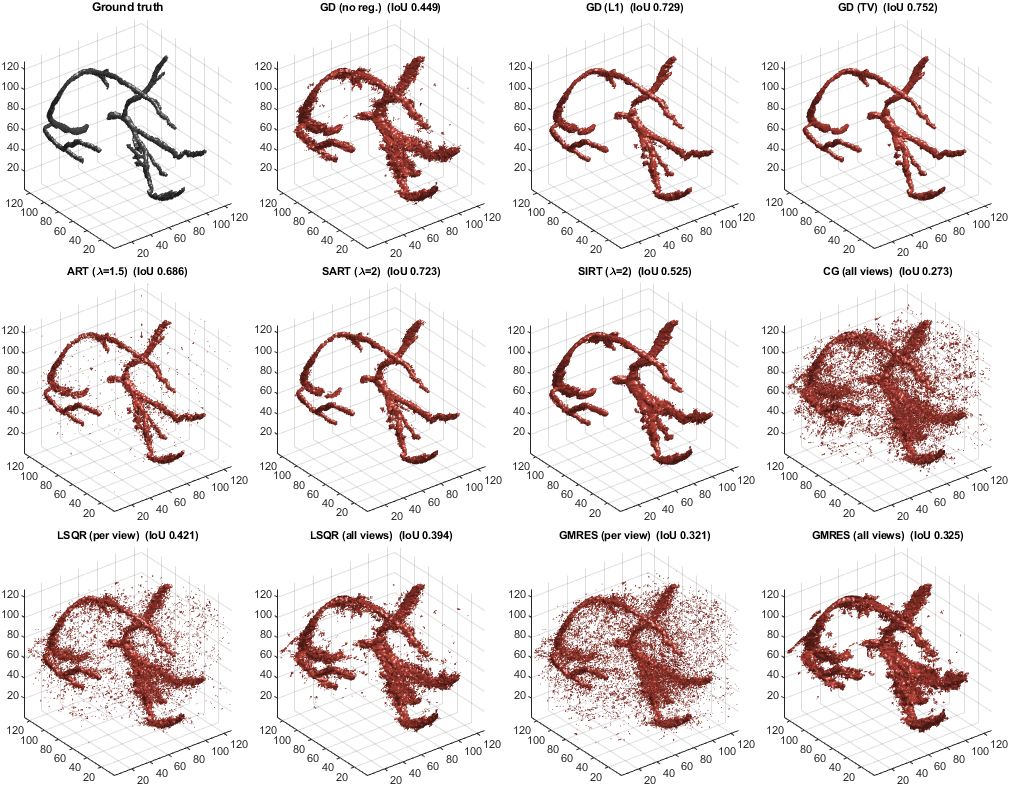}
        \caption{Model $0109$}
        \label{fig:0109_recon}
    \end{subfigure}
    \hfill
    \begin{subfigure}[b]{0.85\textwidth}
        \centering
        \includegraphics[width=\linewidth]{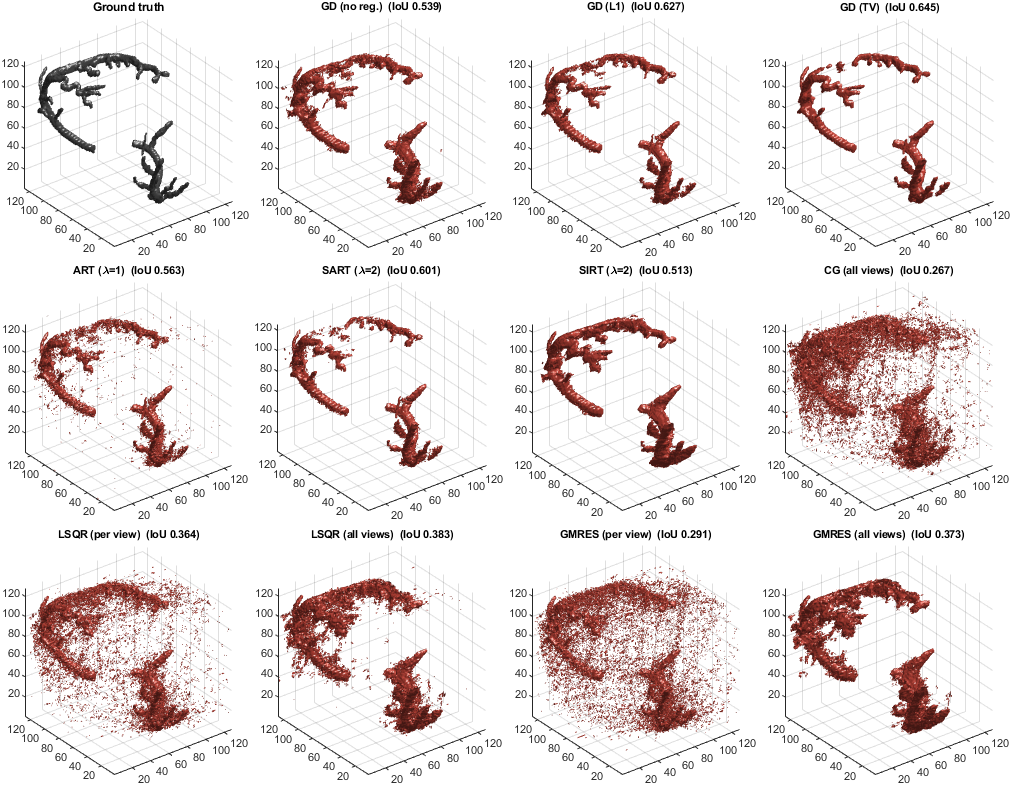}
        \caption{Model $0148$}
        \label{fig:0148_recon}
    \end{subfigure}
    \caption{Reconstruction analysis: Models $0109$ and $0148$.}
    \label{fig:group3_recon}
\end{figure}

\begin{figure}[H]
    \centering
    \begin{subfigure}[b]{0.85\textwidth}
        \centering
        \includegraphics[width=\linewidth]{Images/Legend.png}
        \caption{Legend mapping all the figures to the corresponding reconstruction technique}
        \label{fig:leg_1}
    \end{subfigure}
    \hfill
    \begin{subfigure}[b]{0.85\textwidth}
        \centering
        \includegraphics[width=\linewidth]{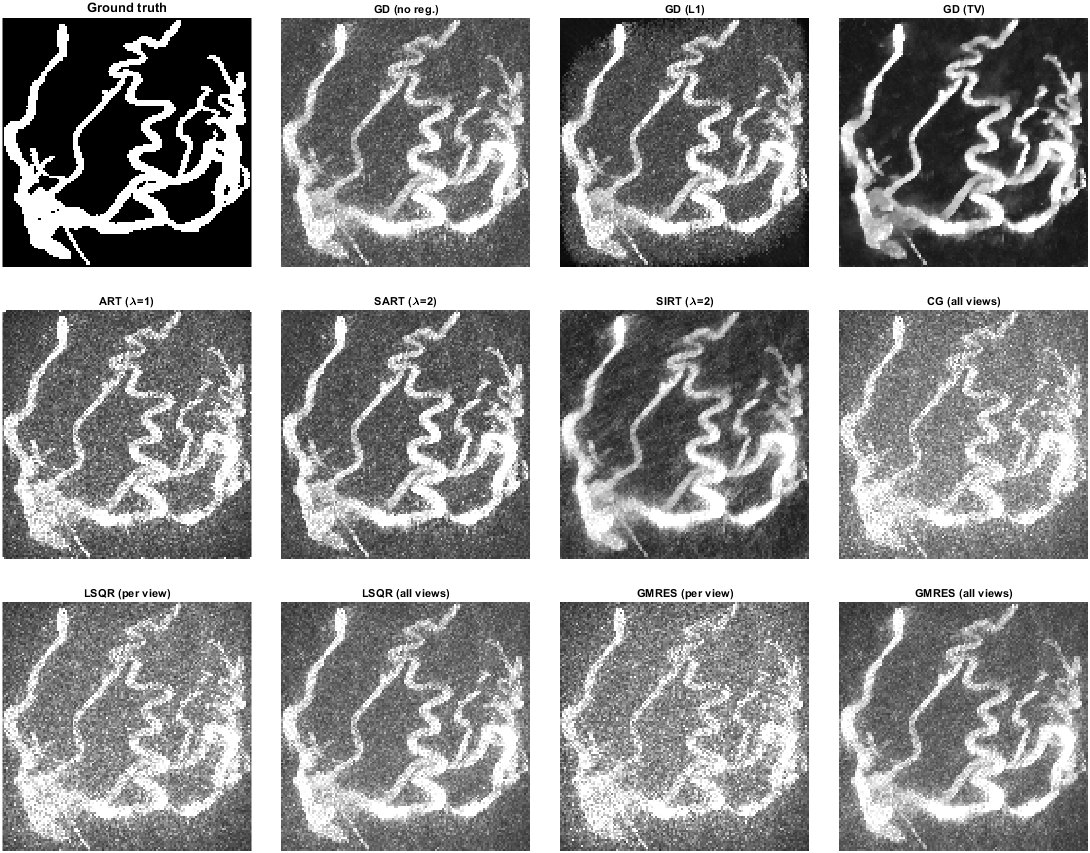}
        \caption{Model $0064$}
        \label{fig:0064_proj}
    \end{subfigure}
    \caption{Held-out Projections: Legend and Model $0064$.}
    \label{fig:group1_proj}
\end{figure}

\begin{figure}[H]
    \centering
    \begin{subfigure}[b]{0.85\textwidth}
        \centering
        \includegraphics[width=\linewidth]{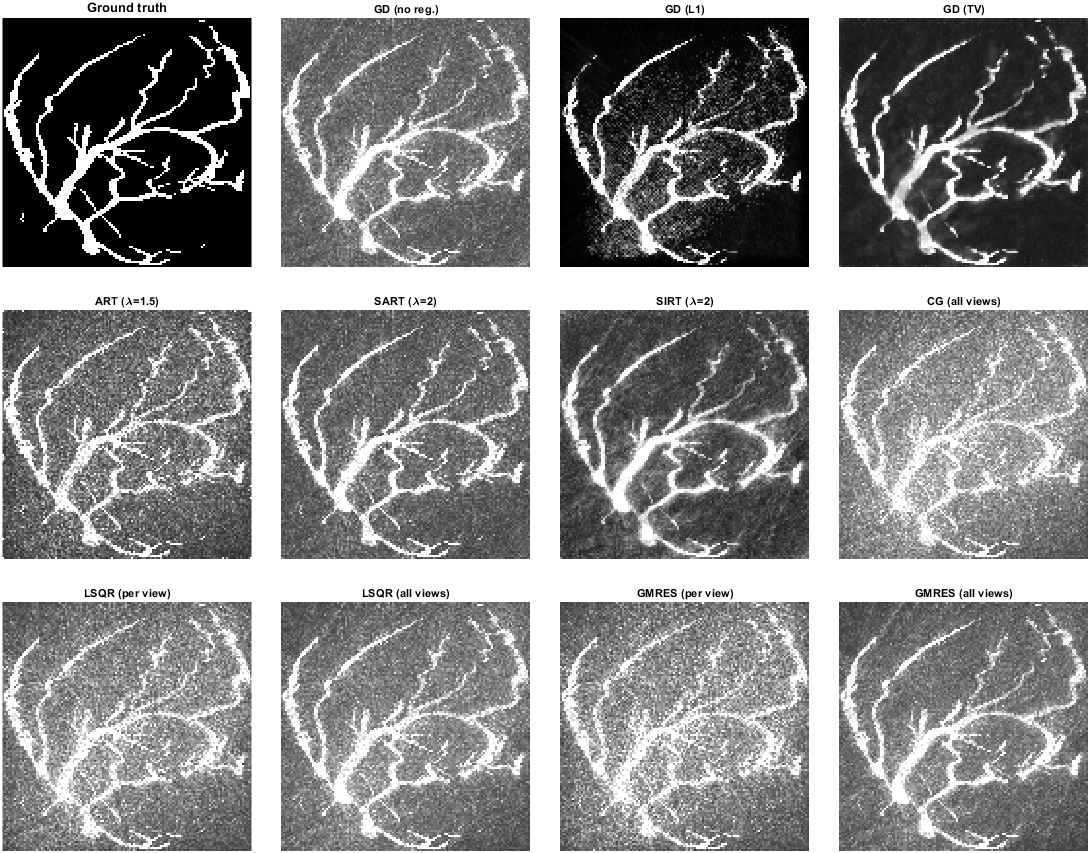}
        \caption{Model $0069$}
        \label{fig:0069_proj}
    \end{subfigure}
    \hfill
    \begin{subfigure}[b]{0.85\textwidth}
        \centering
        \includegraphics[width=\linewidth]{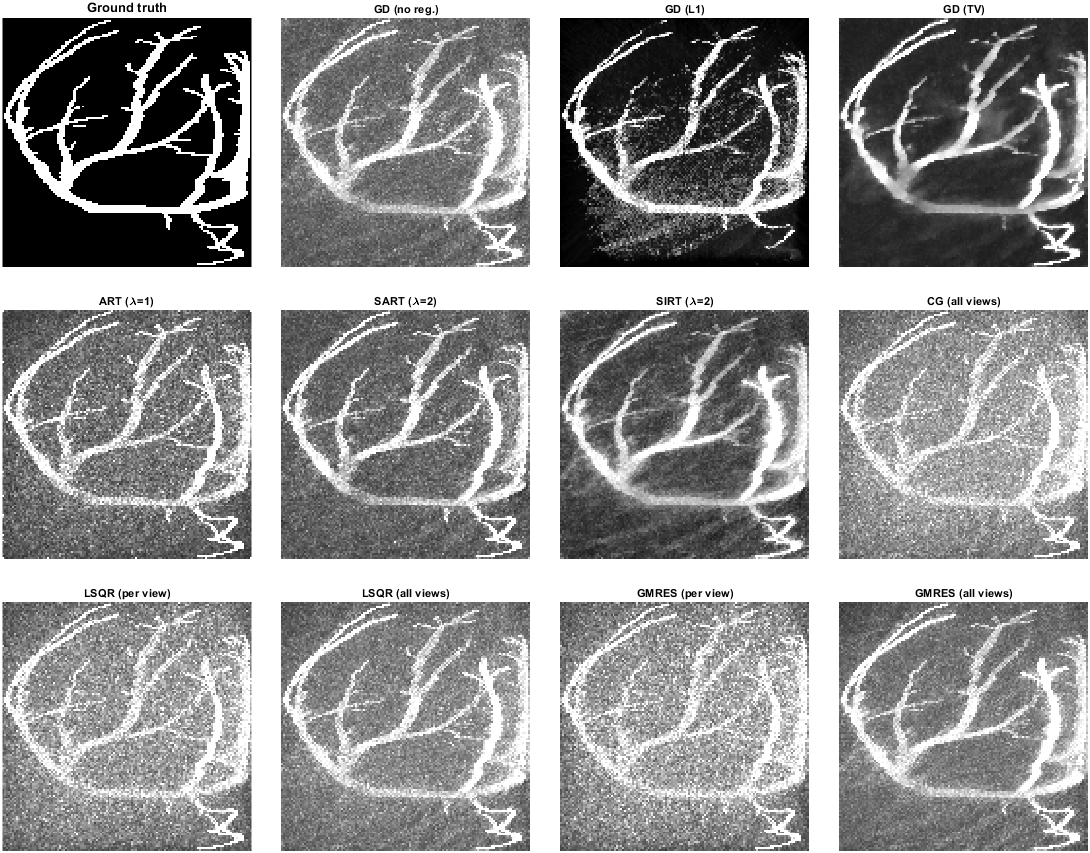}
        \caption{Model $0088$}
        \label{fig:0088_proj}
    \end{subfigure}
    \caption{Held-out Projections: Models $0069$ and $0088$.}
    \label{fig:group2_proj}
\end{figure}

\begin{figure}[H]
    \centering
    \begin{subfigure}[b]{0.85\textwidth}
        \centering
        \includegraphics[width=\linewidth]{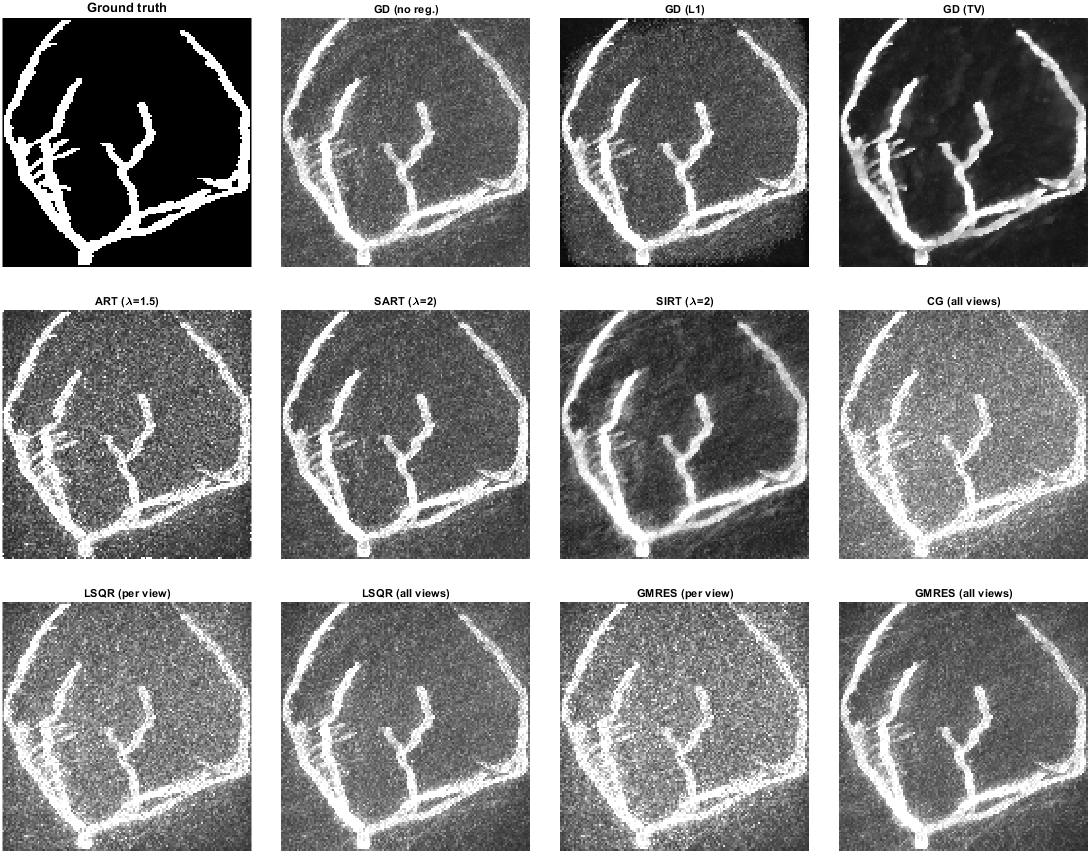}
        \caption{Model $0109$}
        \label{fig:0109_proj}
    \end{subfigure}
    \hfill
    \begin{subfigure}[b]{0.85\textwidth}
        \centering
        \includegraphics[width=\linewidth]{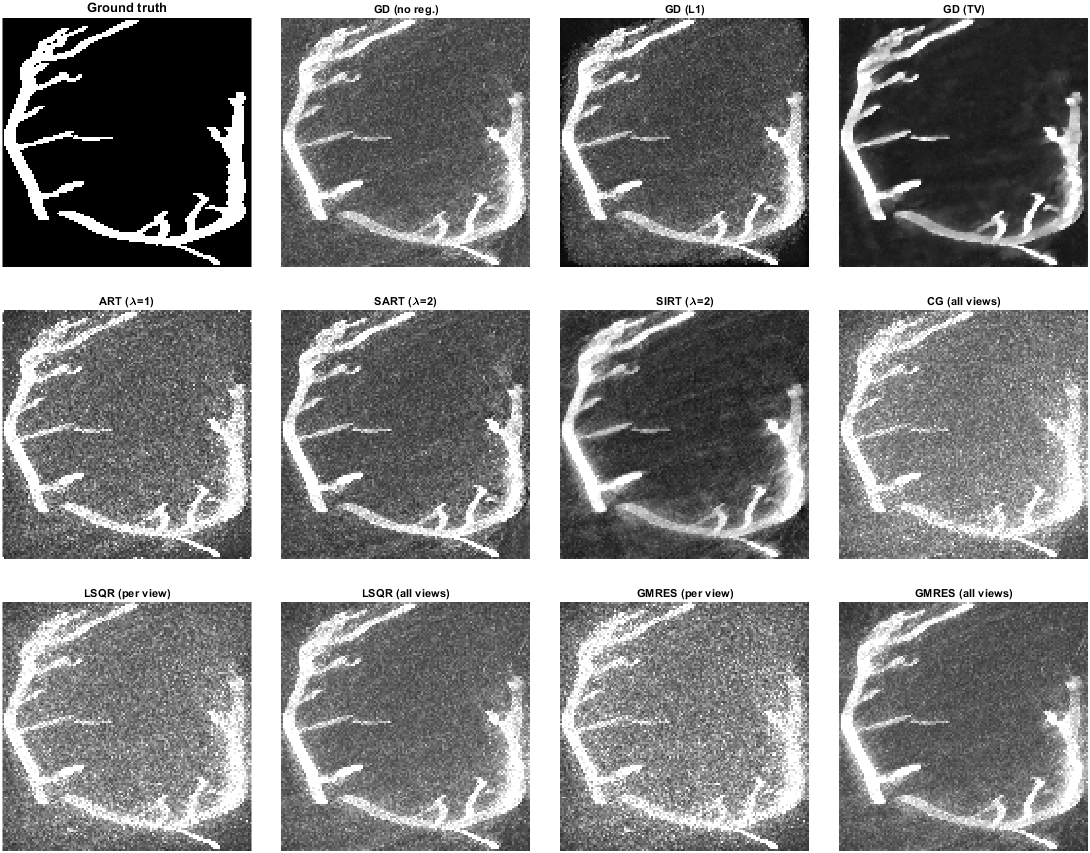}
        \caption{Model $0148$}
        \label{fig:0148_proj}
    \end{subfigure}
    \caption{Held-out Projections: Models $0109$ and $0148$.}
    \label{fig:group3_proj}
\end{figure}

A visual comparison of the methods is provided across ten figures. For each of the five test models, we present the 3D reconstructed vessel volumes next to the ground truth for 11  evaluated parameterized reconstruction methods, displayed in Figures \ref{fig:group1_recon} to \ref{fig:group3_recon}. The visual appearance of the 3D structures correlates with the achieved IoU, displayed on top of each sub-figure.
For the same parameterized reconstruction methods, one
parallel projection of the corresponding volumes is displayed at an orientation not used during reconstruction, on Figures \ref{fig:group1_proj} to \ref{fig:group3_proj}. These visual results align with the quantitative data, such as the SSIM evaluations.
They highlight the structural accuracy of the regularized gradient methods and the noisy reconstruction of the Krylov solvers.

\subsubsection{Volume Noise}
To further investigate the sensitivity and robustness of the reconstruction algorithms against volumetric perturbations, we evaluated the effect of varying the noise level $\nu \in [0, 0.5]$ in Eq.\eqref{eq_vol_noise} at increments of $0.05$. Reconstructions were compared across three gradient descent optimization schemes: unregularized GD, $L_1$-regularized GD, and TV-regularized GD.

\begin{figure}[H]
    \centering
    \includegraphics[width=1\linewidth]{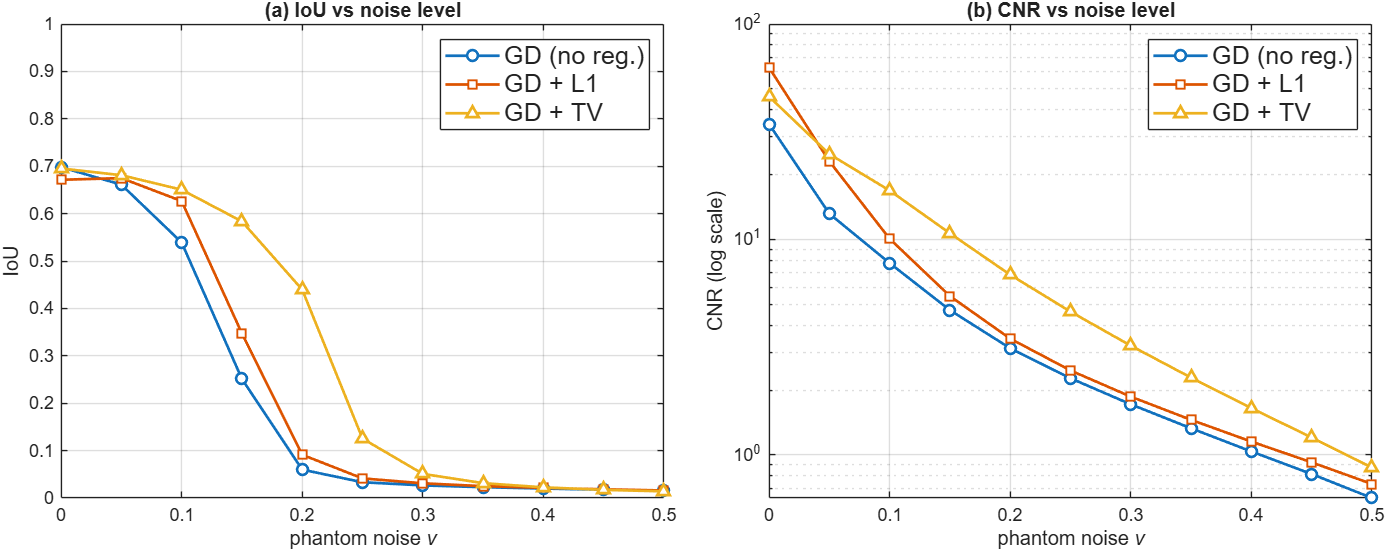}
    \caption{The IoU achieved by three versions of GD for model 0148 under different values of volume noise parameters $\nu.$}
    \label{fig:noise_iou_cnr}
\end{figure}

In the noise-free case ($\nu = 0$), all methods achieved peak reconstruction quality with an Intersection over Union ($\text{IoU}$) of $\approx 0.7$. As expected, increasing $\nu$ led to a monotonic reduction in $\text{IoU}$. At $\nu = 0.1$, which has been the baseline noise level selected for the preceding experiments, the reconstruction $\text{IoU}$ dropped to $\text{IoU} \approx 0.6$. When the perturbation level reached $\nu = 0.2$, $\text{GD} + \text{TV}$ demonstrated superior robustness by maintaining an $\text{IoU} \approx 0.45$, whereas the IoU of unregularized GD and $\text{GD} + L_1$ dropped to below $0.1$. Beyond $\nu \ge 0.3$, the $\text{IoU}$ of all three schemes approaches zero.
Similar trends are observed for CNR, as illustrated in Figure~\ref{fig:noise_iou_cnr}.

\subsection{Real phantom data}
In order to extend the evaluation of our methods we utilized data from a static phantom containing vessel-like structures \cite{Vukicevic2018} (dataset I).
The input consists of nine X-ray angiography cine runs of a static coronary phantom, acquired on a GE DL system (512×512 pixels, 0.39 mm detector pixels, 8-bit).
Each run provides primary and secondary C-arm angles, the source-to-detector distance (850–976 mm) and the source-to-patient distance (720 mm) in its DICOM header.
Frames were averaged within each run, excluding frames affected by the recursive-filter transient or, in IM0, by table and magnification changes.
Line integrals were obtained as $p = \operatorname{log}(I_bg/I)$, where the background $I_bg$ was estimated by morphological closing within a collimator mask.

The acquisition geometry was modelled as a C-arm with rotation with the primary and secondary angles $\alpha,\beta$ following the DICOM definition.
The unknown angle signs, image orientation and per-view detector offsets were estimated from the projections alone.
Detector offsets were then refined on the full-resolution grid until the updates fell below 0.25 px.
Volumes of $144^3$ voxels (0.9 mm) were reconstructed from eight views, with a ninth view (IM0) held out for validation.
Exact ray–voxel system matrices were used for each view.
The reconstruction applies sequential per-view gradient descent with a diagonal majorizer step and a non-negativity constraint, optionally regularized by L1 or isotropic total variation.
Reconstructions were segmented by a fixed hysteresis-Otsu rule that uses no ground truth.
They were compared with the voxelized MSCT surface, rigidly aligned by ICP, in terms of IoU, CNR, centerline distance and diameter error.
Vessel diameters were also measured independently of any threshold, as the full width at half maximum of intensity profiles perpendicular to the centerline.
Fidelity was further assessed by SSIM, NRMSE and correlation between the held-out projection and the reprojection of each reconstruction.

Limitations of the data
The angiograms are 8-bit, display-processed (DISP) or dynamically range-compressed (DRM) images, and the vessel cores are clipped to zero.
The unknown intensity transfer function therefore prevents recovery of true line integrals.
The DICOM header lacks the angle convention, the image orientation, the detector offsets and the object pose. The runs also span two sessions with an apparent object displacement, and the first cine dicom file (IM0) includes table and magnification changes not recorded in its single-valued header.
Only eight views with a limited angular span are available, and there are no flat-field, dark or exposure-normalization data.
The MSCT reference surface lies at about 18\% of the air-to-peak contrast and overestimates the vessel diameter by about 1.4 mm (27\%), while the CT itself saturates at 3071 Housfield Units -HU.

\begin{figure}{H}
    \centering
    \includegraphics[width=1\linewidth]{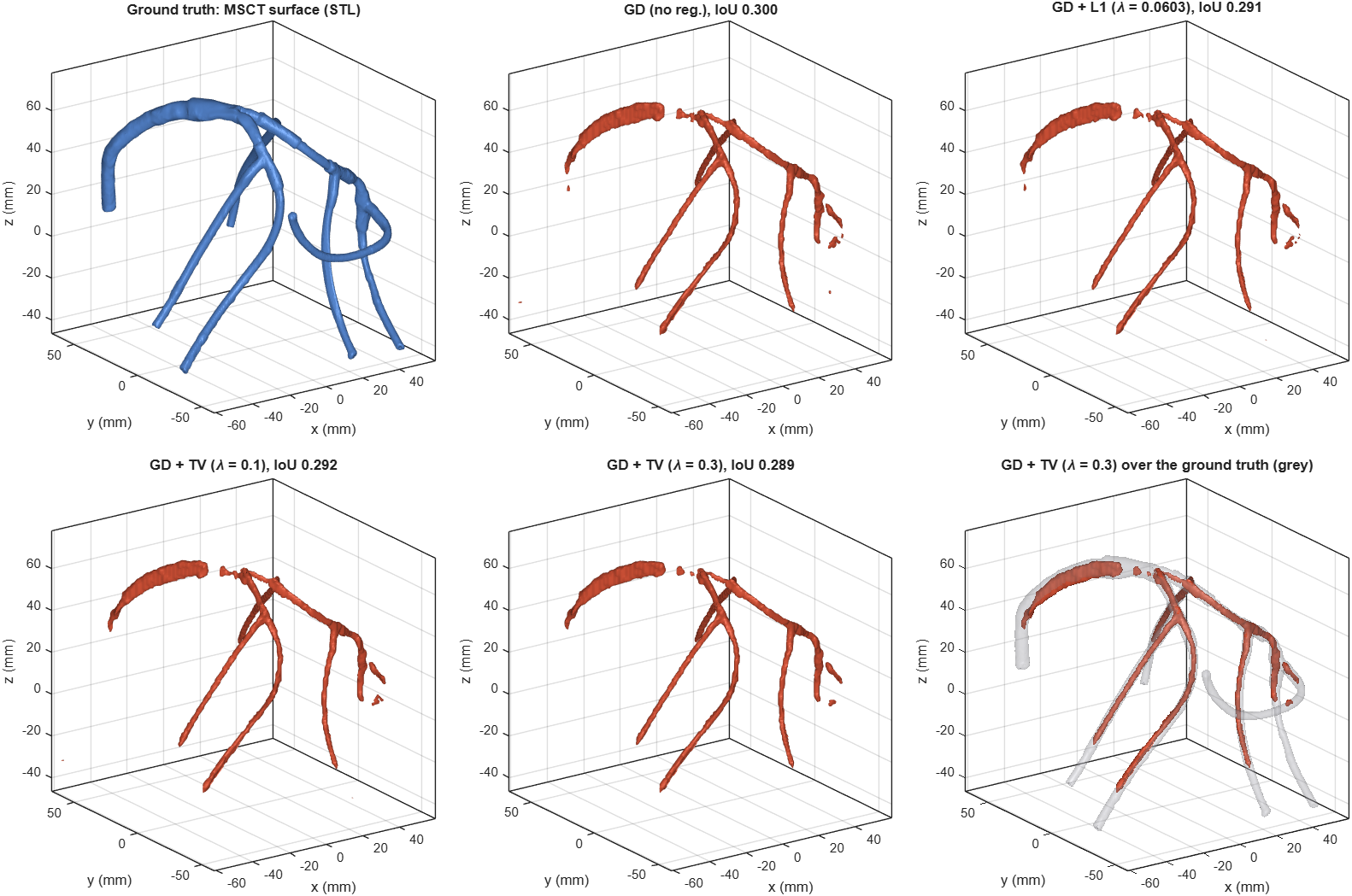}
    \caption{3D volume representations comparing the ground-truth STL (top-left, blue) with four reconstructions: no regularizer (top-middle), $L_1$ regularization (top-right), and TV regularization with two different $\lambda$ values (bottom-left and bottom-middle). The bottom-right panel displays the final reconstructed vessel overlaid on the ground truth (gray).}
    \label{fig:ph_iso}
\end{figure}

\begin{figure}{H}
    \centering
    \includegraphics[width=1\linewidth]{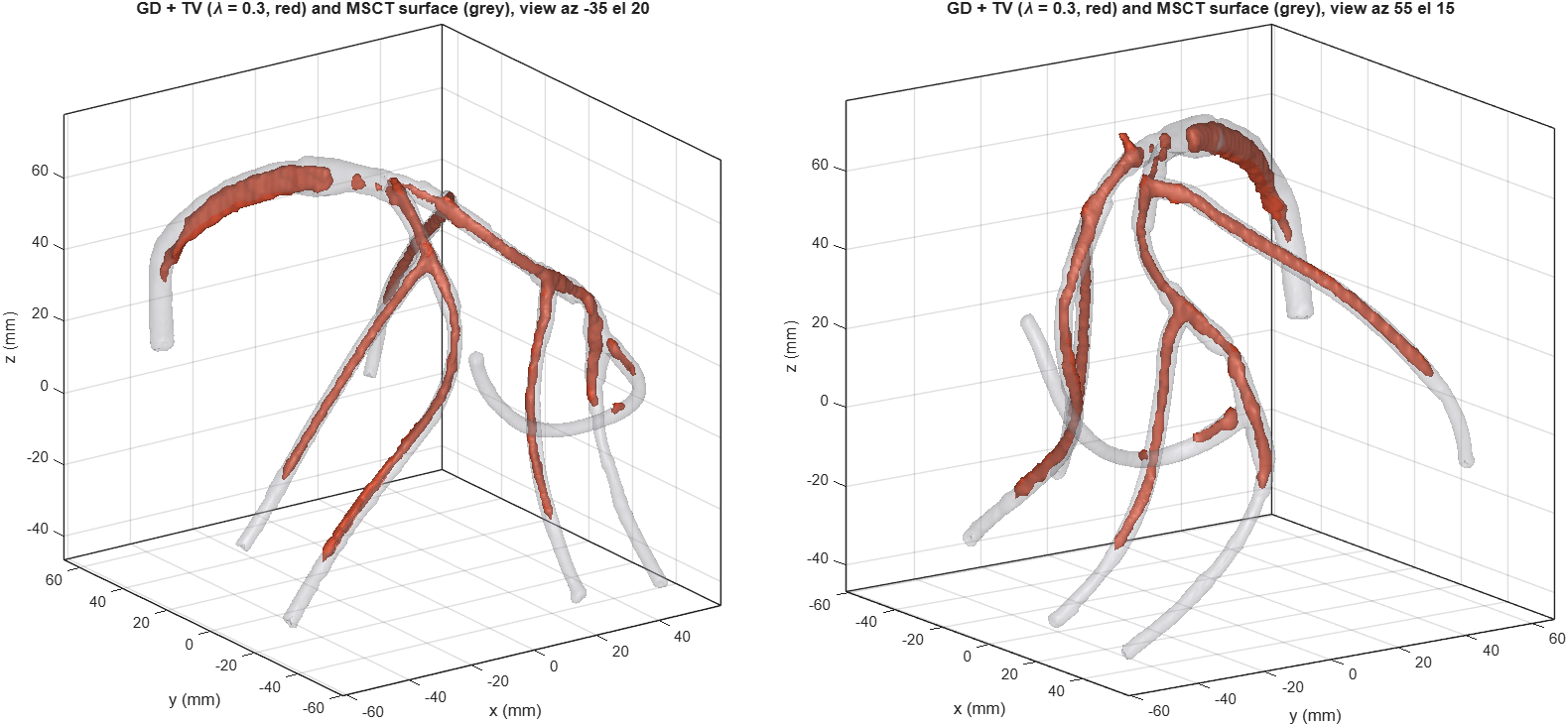}
    \caption{TV-regularized ($\lambda = 0.3$) gradient descent reconstruction (colored) overlaid on the ground truth (gray) at two viewing angles: azimuth -53°, elevation 20° (left) and azimuth 55°, elevation 15° (right).}
    \label{fig:ph_iso_overlay}
\end{figure}

We compared the ground-truth against our results of the gradient descent reconstruction of the phantom. For the reconstruction we investigated the effects of L1 and TV regularization. Figure~\ref{fig:ph_iso} displays the 3D volume of the reconstructed data alongside the ground-truth STL geometry. To display spatial fidelity in detail, Figure~\ref{fig:ph_iso_overlay} overlays the TV-regularized result ($\lambda = 0.3$) onto the ground truth across two viewing angles.

\begin{figure}{H}
    \centering
    \includegraphics[width=1\linewidth]{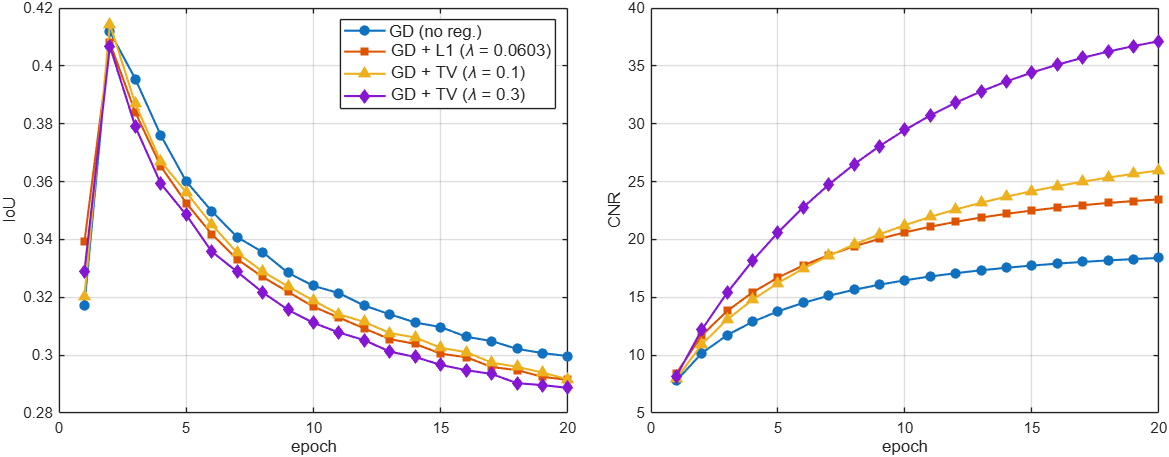}
    \caption{Evaluation of IoU (left) and CNR (right) across epochs for four reconstruction approaches: unregularized GD, GD + $L_1$, and GD + TV regularization ($\lambda = 0.1$ and $\lambda = 0.3$).}
    \label{fig:ph_iou_cnr}
\end{figure}

\begin{figure}[H]
    \centering
    \includegraphics[width=1\linewidth]{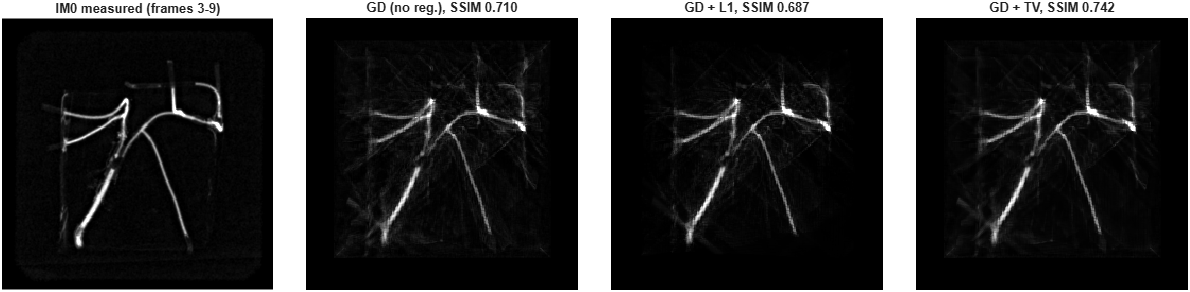}
    \caption{Comparison of the ground-truth held-out projection against projections of the reconstructions using unregularized gradient descent (GD), GD with $L_1$ regularization, and GD with Total Variation (TV) regularization ($\lambda = 0.3$).}
    \label{fig:placeholder}
\end{figure}

Unregularized gradient descent (GD) yielded the highest accuracy overall. The quantitative metrics in Figure~\ref{fig:ph_iou_cnr} confirm this result, as unregularized GD achieves superior Intersection over Union (IoU) and Contrast-to-Noise Ratio (CNR) values throughout the epochs. Finally, Figure~\ref{fig:placeholder} evaluates model performance on forward projections at a held-out angle, where the unregularized GD projection aligns most closely with the ground truth.

\section{Conclusions and Further Work}
This work presented a unified matrix-free framework for three-dimensional reconstruction from sparse X-ray projections and a comparative numerical evaluation ofalgebraic, gradient-based, and Krylov subspace methods.
The framework employs efficient ray--voxel intersection tables for forward and backprojection, avoiding explicit construction of the large system matrix.

Experiments on five patient-derived coronary artery models demonstrated that gradient descent (GD) with diagonal Hessian majorization and total variation (TV) regularization provided the most favorable overall reconstruction quality.
It consistently achieved high volumetric overlap, contrast-to-noise ratio (CNR), and structural similarity (SSIM) at an unseen projection angle, while preserving vessel geometry. The diagonal majorizer also improved convergence compared with the conventional scalar learning
rate and Jacobi preconditioning. The noise sensitivity experiments further demonstrated the superior robustness of TV regularization under increasing volumetric perturbations.

Among the algebraic reconstruction methods, ART and SART provided competitive geometric accuracy at relatively low computational cost, although their performance depended on the relaxation parameter. SIRT exhibited effective noise suppression but slower convergence and inferior held-out projection similarity.
The investigated Krylov methods, including CG, LSQR, and GMRES, generally produced lower geometric accuracy at substantially greater computational cost under the selected iteration budgets. These findings concern the evaluated implementations and do not preclude improvements through Krylov regularization or alternative stopping
strategies.

Nesterov momentum accelerated GD for selected parameter values, but did not consistently improve upon the per-view TV-regularized GD results, highlighting the sensitivity of accelerated schemes to their momentum parameters and update schedules. The update per-epoch, which is mandadtory for Nesterov's GD includes a substantial drwaback compared to pero-view updates.

Evaluation on experimentally acquired coronary phantom projections provided an additional assessment under realistic acquisition conditions. In contrast to the
simulated coronary models, unregularized GD achieved the best IoU, whereas the GD-TV ($\lambda=0.3$) achieved the best CNR. the limitations of the phantom data, as previously indicated may be the cause for this observation. 

Overall, the results highlight the importance of regularization, gradient scaling, and update scheduling in severely underdetermined tomographic reconstruction.
Future work will investigate adaptive regularization and stopping criteria, accelerated and regularized
Krylov methods, and the influence of projection-angle selection on the recovery of vascular geometry. We will also investigate the hybridization of these methods with deep learning approaches.

\clearpage
\bibliographystyle{etna}       
\bibliography{references}
\end{document}

%% file: table_bibliography.tex

\begingroup
\footnotesize
\setlength{\tabcolsep}{3pt}
\renewcommand{\arraystretch}{1.18}
\setlength{\LTleft}{0pt}
\setlength{\LTright}{0pt}

\begin{longtable}{@{}
  p{0.13\linewidth}
  p{0.23\linewidth}
  p{0.14\linewidth}
  p{0.23\linewidth}
  p{0.19\linewidth}
  @{}}

\caption{Representative studies on iterative tomographic
reconstruction, sparse-view CT, and numerical optimization.}
\label{tab:literature_review}\\

\toprule
\textbf{Reference} &
\textbf{Reconstruction method} &
\textbf{Sparse-view setting} &
\textbf{Experimental data} &
\textbf{Evaluation criteria} \\
\midrule
\endfirsthead

\multicolumn{5}{c}{
\textit{Table \thetable{} (continued)}}\\
\toprule
\textbf{Reference} &
\textbf{Reconstruction method} &
\textbf{Sparse-view setting} &
\textbf{Experimental data} &
\textbf{Evaluation criteria} \\
\midrule
\endhead

\midrule
\multicolumn{5}{r}{\textit{Continued on next page}}\\
\endfoot

\bottomrule
\endlastfoot

\multicolumn{5}{@{}l}{
\textit{Sparse-view CT and imaging-oriented methods}}\\
\addlinespace

Sidky and Pan
\cite{Sidky2008}
&
ASD-POCS; constrained TV minimization,
adaptive gradient descent, positivity.
&
Yes; sparse and limited-angle 3D CBCT.
&
Simulated 3D Defrise disk and FORBILD
jaw phantoms.
&
Data discrepancy, TV, convergence,
visual reconstruction quality.
\\ \addlinespace

Choi et al.
\cite{choi2010compressed}
&
Nesterov-based first-order optimization;
TV and weighted data fidelity.
&
Yes; undersampled 3D CBCT.
&
Digital 3D phantom and experimentally
acquired anthropomorphic phantom data.
&
CNR, intensity profiles,
objective convergence, image quality.
\\ \addlinespace

Bian et al.
\cite{Bian2010}
&
TV-based ASD-POCS; comparison
with FDK and iterative methods.
&
Yes; 15--96 views.
&
Measured CBCT projections from
physical cylindrical and head phantoms.
&
UQI, mutual information, CNR,
noise spectrum, lesion detectability.
\\ \addlinespace

Sidky et al.
\cite{Sidky2012}
&
Chambolle--Pock primal--dual method;
convex data-fidelity and TV formulations.
&
Yes; undersampled 2D fan-beam CT.
&
Simulated low-dose breast CT
phantom with Poisson noise.
&
Convergence, noise characteristics,
ROI and visual image assessment.
\\ \addlinespace

Sabat{\'e} Landman et al.
\cite{SabateLandman2023}
&
Krylov methods: CGLS, LSQR, LSMR
and GMRES variants; regularization.
&
Yes; includes sparse-view 3D CBCT.
&
Synthetic CT, measured head-phantom
CBCT and industrial micro-CT.
&
Residual and reconstruction errors,
convergence, computational time,
visual quality.
\\ \addlinespace

Chillar{\'o}n P{\'e}rez et al.
\cite{ChillaronPerez2024}
&
LSMB, LSQR, LSMR; block variants,
regularization and acceleration.
&
Yes; few-view 2D CT.
&
Clinical CT images with simulated
undersampled projection data.
&
PSNR, SSIM, convergence,
iterations, execution time.
\\ \addlinespace

Erdo\u{g}an and Fessler
\cite{Erdogan1999}
&
Separable paraboloidal surrogates
(SPS); ordered-subsets transmission
reconstruction.
&
Not specifically sparse-view CBCT.
&
PET transmission phantom and
patient measurements.
&
Objective convergence,
reconstruction and segmentation
accuracy.
\\ \addlinespace

Biguri et al.
\cite{Biguri2016}
&
TIGRE: GPU-based FDK, algebraic,
Krylov and TV reconstruction methods.
&
Supports sparse-view CBCT.
&
Simulated 3D phantoms and measured
CBCT demonstration data.
&
Computation time, numerical
performance, reconstruction quality.
\\ \midrule

\multicolumn{5}{@{}l}{
\textit{Numerical analysis and inverse-problem methods}}\\
\addlinespace

Andersen and Hansen
\cite{AndersenHansen2014}
&
Generalized row-action methods;
relaxed incremental proximal-gradient
algorithms and regularization.
&
Not specifically sparse-view.
&
Simulated 2D tomographic
reconstruction problems.
&
Convergence, reconstruction error,
effects of relaxation and regularization.
\\ \addlinespace

Hansen and J{\o}rgensen
\cite{HansenJorgensen2018}
&
AIR Tools II: ART, SIRT and related
algebraic methods; relaxation and
stopping rules.
&
Supports sparse-view geometries.
&
Synthetic CT and other
tomographic test problems.
&
Convergence, stopping criteria,
accuracy, computational performance.
\\ \addlinespace

Chung and Gazzola
\cite{ChungGazzola2019}
&
Flexible Krylov methods;
Golub--Kahan/Arnoldi,
hybrid $\ell_p$ regularization.
&
Not specifically sparse-view.
&
Simulated tomographic reconstruction
and image-deblurring problems.
&
Reconstruction error, convergence,
regularization performance,
computational efficiency.
\\ \addlinespace

J{\o}rgensen and Sidky
\cite{JorgensenSidky2015}
&
Sparsity-regularized optimization;
$\ell_1$ and TV minimization;
phase-diagram analysis.
&
Yes; 1--26 views in small-scale tests.
&
Synthetic sparse images and
2D walnut-derived phantoms.
&
Relative reconstruction error,
recovery probability,
phase-transition diagrams.
\\
\midrule
\multicolumn{5}{@{}l}{
\textit{Coronary artery reconstruction methods}}\\
\addlinespace

Liu et al.
\cite{Liu:16}
&
Regularized iterative reconstruction;
sparsity constraints and 3D/2D
motion compensation.
&
Yes; 8 gated projections
in principal experiments.
&
Simulated rotational angiograms
from a dynamic 4D XCAT
coronary phantom.
&
Morphological overlap (MMO),
relative radius error (RRE),
reprojection correlation,
computation time.
\\ \addlinespace

Neubauer et al.
\cite{neubauer2010clinical}
&
Fully automated ECG-gated
iterative 3D reconstruction
from rotational angiography.
&
Yes; angularly undersampled
cardiac-gated acquisition.
&
Clinical rotational coronary
angiograms from 23 patients.
&
Image quality, lesion detection,
QCA vessel diameters and
segment lengths.
\\

Hansis et al.
\cite{Hansis2008}
&
Iterative sparse-object reconstruction with
$\ell_1$ regularization and morphological
background suppression.
&
Yes; 5--10 projections per cardiac phase.
&
Simulated coronary artery phantoms;
preliminary clinical ECG-gated
rotational angiography.
&
Reconstructed attenuation coefficients,
vessel radius error, image contrast,
visual reconstruction quality.
\\ \addlinespace

\end{longtable}
\endgroup

%% file: gd_clean.tex
\subsection{Gradient Descent}

\subsubsection{Gradient Formulation}
Let us define the residual vector $\mathbf{r} = \mathbf{Y} - \mathbf{\hat{Y}}$. We formulate the least square error as the objective function:
\begin{equation}
E(\boldsymbol{\mu}) = \sum_{p} \left(Y_p - \hat{Y}_p(\boldsymbol{\mu}) \right )^2 = \frac{1}{2}\lVert\mathbf{r}\rVert^2
    \label{eq:of}
\end{equation}

Taking the partial derivative with respect to any $\mu_v$, we obtain:
\begin{align}
\frac{\partial E}{\partial \mu_v} = - \sum_{p} r_p \cdot A_{pv}
\end{align}

Using the equation \eqref{eq:of}  and standard matrix differentiation The following Lemma is obtained.


\begin{lemma}
The gradient and the Hessian of $E(\boldsymbol{\mu})$ with respect to $\boldsymbol{\mu}$ are given by:
\begin{equation}
\nabla_{\boldsymbol{\mu}} E(\boldsymbol{\mu}) = -\mathbf{A}^T( \mathbf{Y} - \mathbf{A}\boldsymbol{\mu}), \quad \mathbf{H}_{E(\boldsymbol{\mu})} = \mathbf{A}^T \mathbf{A}
\end{equation}

\end{lemma}

The voxel intensities get updated iteratively. To ensure physical non-negativity, we apply constraint at each step:

\begin{align}
    \mu_v^{(t)} &= \max \left(\mu_v^{(t-1)} - \alpha_v \cdot \frac {\partial E}{\partial_{\mu_v}}, 0 \right)
    = \max \left( \mu_v^{(t-1)} + \alpha_v \sum_{p} r_p \cdot A_{pv}, 0 \right)
    \label{eq:mu_update_no_reg}
\end{align}
where $\boldsymbol{\alpha} \in \mathbb{R}^{N_v}$ is the vector defining the voxel-wise step of the update. The above update can expressed in matrix form as:
\begin{align}
\boldsymbol{\mu}^{(t)} &= \max \left( \boldsymbol{\mu}^{(t-1)} - \boldsymbol{\alpha} \odot \nabla_{\boldsymbol{\mu}} E\left(\boldsymbol{\mu}^{(t-1)}\right), 0 \right) = \max \left( \boldsymbol{\mu}^{(t-1)} + \boldsymbol{\alpha} \odot \mathbf{A}^T(\mathbf{Y} - \mathbf{A}\boldsymbol{\mu}^{(t-1)}), 0 \right).
\label{eq:upRule}
\end{align}



\paragraph{Gradient Formulation without constructing $\mathbf{A}$}

The gradient and the update rule can be formulated with our auxiliary lookup tables \ref{sec:aux_tab}.
For every voxel $v$ of volume $\boldsymbol{\mu}$
\begin{equation}
    {\mu}_v^{(t)} = \max \left( {\mu}_v^{(t-1)} + {\alpha}_v \cdot (\mathbf{VL}_{:,v})^T \cdot \mathbf{r}_{\mathbf{VP}_{:,v}}, 0 \right)
\end{equation}

Alternatively, the correction of the volume can be expressed in relation to the pixel matrices. For projection pixel (ray) $p$, let
\begin{equation}
    \mathcal{V}_p = \left\{ v \;:\; A_{pv} \neq 0 \right\}
\end{equation}
denote the set of voxels intersected by the corresponding ray. The volume correction for all voxels intersected by the same ray can be updated simultaneously as:
\begin{equation}
    \Delta{\mu}_{\mathcal{V}_p} = \Delta{\mu}_{\mathcal{V}_p} + (\mathbf{PL}_{:,p})^T \cdot r_p
\end{equation}

After all rays have been processed:
\begin{equation}
    \boldsymbol{\mu}^{(t+1)} = \max \left ( \boldsymbol{\mu}^{(t)} + \boldsymbol{\alpha} \odot \Delta\boldsymbol{\mu}, 0 \right)
\end{equation}

\subsubsection{Gradient Scaling}
We examined different approaches to adjust the step size of the updates to improve optimization, specifically we implemented a standard learning rate, a Jacobi preconditioner, and a diagonal majorizer.

\paragraph{\textbf{Learning Rate ${\alpha}$}}
The learning rate $\alpha$ in \eqref{eq:mu_update_no_reg} is chosen as the inverse of the Lipschitz constant of the gradient, which is given by the largest eigenvalue $L = \lambda_{\max}$ of the Hessian $\mathbf{H}=\mathbf{A}^T \mathbf{A}$:
\begin{equation}
    \alpha \approx \frac{1}{\lambda_{\max}},
    \quad
    \lambda_{\max} = \lambda_{\max}(\mathbf{A}^T \mathbf{A}).
\end{equation}

We approximate $\alpha$ without the Hessian and the system matrix $A$ as such:

\begin{itemize}
    \item \textbf{Trace based bound for $\lambda_{max}$}: 
    A straightforward way to bound the spectral properties of $A^T A$ is to examine its diagonal entries. As a first attempt, we assumed that $\lambda_{max}<tr(H).$ In order to calculate the trace of the Hessian (without constructing it) we utilize our implementation as follows:

    \begin{align}
        diag(A A^T)_p = \sum_v (A_{pv})^2 = \sum_{k} ({PL}_{kp})^2, \\
        diag(H)=diag(A^T A)_v = \sum_p (A_{pv})^2 = \sum_{k} ({VL}_{kv})^2
    \end{align}
    
    The trace, $ tr(\mathbf{H}) = \sum_p diag(H)_p$ provides a very loose upper bound on $\lambda_{\max}$ and is overly conservative in practice.

    \item \textbf{Norm-based bound for $\lambda_{max}$:}
    Instead, we use matrix norms to approximate a bound for the largest eigenvalue. Specifically, the spectral norm satisfies the following inequality:

    \begin{equation}
        \sigma_{\max}(A) \le \sqrt{\lVert \mathbf{A}\rVert_1 \lVert \mathbf{A}\rVert_\infty}
    \end{equation}
    
    Considering that the largest eigenvalue of the Hessian $A^T A$ is equal to the square of the largest singular value of $A$, we obtain the bound:
    \begin{equation}
        \lambda_{\max}(\mathbf{A}^T \mathbf{A}) = \sigma_{\max}^2(\mathbf{A}) \le \|\mathbf{A}\|_1 \lVert\mathbf{A}\rVert_\infty
    \end{equation}

    where the matrix norms are: 
    
    \begin{align}
        \lVert \mathbf{A} \rVert_1 &= \max_v \sum_p |A_{pv}|=\max_v \sum_{k} |{VL}_{kv}|, \\
        \lVert \mathbf{A} \rVert_\infty &= \max_p \sum_v |A_{pv}|=\max_p \sum_{k} |{PL}_{kp}|.
    \end{align}

    So finally the max eigenvalue $\lambda_{\max}$ is approximated as:
    \begin{equation}
        \lambda_{\max} \le \|\mathbf{A}|\|_1 \cdot \|\mathbf{A}\|_\infty
        = \max_v \sum_{k} |{VL}_{kv}|
        \cdot \max_p \sum_{k} |{PL}_{kp}|
        \label{eq:standard_lr}
    \end{equation}
\end{itemize}    
    
We evaluated the performance of this estimation on a system matrix of dimensions 
$\mathbf{A} \in \mathbb{R}^{64^2 \times 32^3}$. The exact maximum eigenvalue computed directly from the system matrix is $\lambda_{\max} = 1.7292 \times 10^4$, whereas our proposed method yields an estimated value of $\hat{\lambda}_{\max} = 1.9512 \times 10^4$. This corresponds to a relative estimation error of $12.84\%$, confirming that the proposed approach provides a reliable approximation while significantly reducing computational overhead.

After calculating $a$ we construct the vector $\boldsymbol{\alpha}=[\alpha,\ldots,\alpha]^T \in \mathbb{R}^{N_v}$

\paragraph{\textbf{Jacobi preconditioning}}
The Jacobi preconditioner uses only the diagonal elements of the
Hessian,
\begin{equation}
\mathbf{D} = \operatorname{diag}(\mathbf{H}),
\qquad D_{v} = H_{vv} = \sum_i A_{iv}^{2},
\end{equation}
and it can be expressed as:
\begin{equation}
    D_{v} = \sum_{k} VL_{kv} ^2
    \label{eq:jacobi}
\end{equation}
so $\boldsymbol{\alpha}$ is formulated as ${\alpha}_v = D_v^{-1}.$ If $D_v=0$ then $\alpha_v=0.$ 



\paragraph{\textbf{Diagonal majorization}}
Let 
\begin{equation} \mathbf{1} = \begin{bmatrix} 1 & 1 & \cdots & 1 \end{bmatrix}^{T} \in\mathbb{R}^{N_v}, \end{equation} 
as the all-ones column vector. 
The \(v\)-th element of $\mathbf{H} \mathbf{1}$ vector is equal to the sum of row \(v\) of the Hessian: \begin{equation} \left[\mathbf{H}\mathbf{1}\right]_v = \sum_{u=1}^{N_v}H_{vu}. \end{equation}
A more conservative diagonal matrix is obtained from
\begin{equation}
\mathbf{m}
=
\mathbf{H}\mathbf{1}
=
\mathbf{A}^{T}\mathbf{A}\mathbf{1},
\qquad
\mathbf{M}
=
\operatorname{diag}(\mathbf{m}).
\end{equation}
Thus, applying 
\begin{equation}
M_{vv}=\sum_u H_{vu}=\sum_i A_{iv}\left(\sum_u A_{iu}\right).
\label{eq:major}
\end{equation}

Matrix $\mathbf{M}$ can be generated without constructing $\mathbf{H}$ and $\mathbf{A}.$ More specifically, let $S_p$ denote the total length of ray $p$. Equivalently, $S_p = \mathbf{A}\textbf{1}$, which using \eqref{eq:Ax} becomes $\sum_{k=1}^{K_1} \text{PL}_{k,p}$. Furthermore, let $l_{p,v}$ be the segment length of ray $p$ inside voxel $v$, where $l_{p,v} = \text{PL}_{k,p}$ if $\text{PV}_{k,p} = v$, and $0$ otherwise. The diagonal element $M_{vv}$ is given by:

\begin{equation}
    M_{vv} = \sum_{p=1}^{N_p} l_{p,v} S_p
\end{equation}

For a nonnegative projection matrix \(\mathbf{A}\), $\mathbf{M}-\mathbf{H}$ is PSD, or equivalently matrix \(\mathbf{M}\) satisfies
\begin{equation}
\mathbf{M}-\mathbf{H}\succeq 0,
\end{equation}
and therefore provides a diagonal majorizer of the Hessian. The
resulting step size becomes 
\begin{equation}
    \boldsymbol{\alpha} = \operatorname{diag} ( \mathbf{M}^{-1}). 
    \label{eq:lr_maj}
\end{equation}

For any vector $\mathbf{x}\in\mathbb{R}^{N_v}$,
\begin{align}
\mathbf{x}^T(M-H)\mathbf{x}
&= \sum_v\left(\sum_u H_{vu}\right)x_v^2
   -\sum_{v,u}H_{vu}x_vx_u \nonumber\\
&= \frac{1}{2}\sum_{v,u}H_{vu}(x_v-x_u)^2.
\end{align}
Because $H=A^TA$ and $A\geq 0$ elementwise, $H_{vu}\geq 0$.
Consequently,
\begin{equation}
\mathbf{x}^T(M-H)\mathbf{x}\geq 0
\end{equation}
for all $\mathbf{x}\in\mathbb{R}^{N_v}$, and hence
\(
M-H\succeq 0.
\)
Thus, $M$ is a diagonal majorizer of $H$.

In our reconstruction experiments, Jacobi preconditioning did not
provide satisfactory convergence, whereas the diagonal majorizer
\(\mathbf{M}=\operatorname{diag}(\mathbf{H}\mathbf{1})\) resulted in
substantially more stable and effective updates.

\input{nesterov}

\subsubsection{Optimization}

To prevent overfitting to noise and outliers, we extend the update rule \eqref{eq:mu_update_no_reg} by introducing a regularization penalty. This penalty suppresses excessively large or oscillatory parameter values, favoring simpler solutions that more accurately reflect the underlying vascular anatomy. 
Table \ref{tab:reg} summarizes each regularization term and Lipschitz gradients, and algorithm \ref{alg:gd_l1} demonstrates the explicit workflow logic for gradient descent regularized with $L_1$.

\paragraph{\textbf{L1 Regularizer}}
$L_1$ regularization enforces sparsity within the parameter space by adding a penalty proportional to the absolute magnitude of the weights,  yielding reconstructions resilient to noise \cite{choi2010compressed}.
The objective function $J(\boldsymbol{\mu})$ combining the unregularized loss $E(\boldsymbol{\mu})$ \eqref{eq:of} with an L1 penalty, where $R(\boldsymbol{\mu})=\lVert \boldsymbol{\mu} \rVert_1 = \sum_v \lvert \mu_v \rvert$ is:
\begin{equation}
    J(\boldsymbol{\mu}) = E(\boldsymbol{\mu}) + R(\boldsymbol{\mu}) = E(\boldsymbol{\mu})+ \lambda \lVert \boldsymbol{\mu} \rVert_1,
\end{equation}

The gradient of the objective with respect to $\boldsymbol{\mu}$ is
\begin{equation}
    \nabla_{\boldsymbol{\mu}} J = \frac{\partial E}{\partial \boldsymbol{\mu}} + \frac{\partial}{\partial \boldsymbol{\mu}} \left( \lambda \lVert \boldsymbol{\mu} \rVert_1 \right) = - \sum_{p} r_p \cdot A_{pv} + \lambda \cdot \text{sgn}(\mu_v).
\end{equation}

The Lipschitz gradient for the step size $\alpha$ is $L =\lVert \mathbf{A}\rVert_2^2$ and stems exclusively from the smooth data term. Thus $\alpha$ is bound: $0<\alpha < \frac{2}{\lVert \mathbf{A}\rVert_2^2}.$

The resulting voxel-based update rule is with a non-negativity constraint becomes:
\begin{align}
    \mu_v^{t} &= \max \left( \mu_v^{t-1} + \alpha_v \left( \sum_{p} r_p \cdot A_{pv} - \lambda \cdot \text{sgn}(\mu_v) \right), 0 \right) \\
    &=  \max \left( \mu_v^{t-1} + \alpha_v \left( {\mathbf{VL}_{:,v} \cdot \mathbf{r}_{{\mathbf{VP}}_{:,v}}} - \lambda \cdot \text{sgn}(\mu_v) \right), 0 \right)
    \label{eq:l1}
\end{align}

\paragraph{\textbf{Zero-order Tikhonov Regularization (L2)}}

L2 regularization, or 0-order Tikhonov regularization \cite{tikhonov1963solution}, is employed to
stabilize ill-posed 3D tomographic reconstruction large coefficient magnitudes to the objective function. This suppresses artifact amplification and forces a physically realistic solution \cite{Lou2010}.
The L2 regularization term is the functional  
\begin{equation}
    R(\mu) = \lVert \boldsymbol{\mu} \rVert_2^2 = \int_{\Omega} \lvert \mu(x) \rvert^2 \, dx.
    \label{eq:L2}
\end{equation}
The partial derivative of the complete objective function $J(\boldsymbol{\mu})$ with respect to voxel $\mu_v$ is:
\begin{equation}
    \frac{\partial J}{\partial \mu_v} = -\sum_{p} \left( Y_p - \hat{Y}_p \right) A_{pv} + \lambda \mu_v
    \label{eq:l2_full_grad}
\end{equation}
The Lipschitz gradient is derived as $\lVert \mathbf{A}\rVert_2^2 + \lambda$ and the resulting update rule is:
\begin{align}
    \mu_v^{t} &= \max \left( \mu_v^{t-1} + \alpha_v \left( \sum_{p} r_p \cdot A_{pv} - \lambda \mu_v^{t-1} \right), 0 \right) \\
    &= \max \left( \mu_v^{t-1} + \alpha_v \left( \mathbf{VL}_{:,v} \cdot \mathbf{r}_{\mathbf{\mathbf{VP}_{:,v}}} - \lambda \mu_v^{t-1} \right), 0 \right).
    \label{eq_l2}
\end{align}

\paragraph{\textbf{First-order Tikhonov regularizer}}

First-order Tikhonov regularization penalizes the spatial gradients of the image, promoting piecewise smoothness. For a discrete volume $\boldsymbol{\mu}$, where voxel $v$ corresponds to 3D grid coordinates $(i,j,k)$, the penalty is defined as:
\begin{equation}
    R(\boldsymbol{\mu}) = \sum_{v} \lVert \nabla \mu_v \rVert_2^2 = \sum_{i,j,k} \left[ (\nabla_x \mu_{i,j,k})^2 + (\nabla_y \mu_{i,j,k})^2 + (\nabla_z \mu_{i,j,k})^2 \right].
\end{equation}

The partial derivative of this regularization term with respect to voxel $\mu_{i,j,k}$ can be expressed using a discrete $3 \times 3 \times 3$ Laplacian kernel $K$, which approximates the continuous Laplacian operator $\Delta$:
\begin{equation}
    \frac{\partial R}{\partial \mu_{i,j,k}} = -2 (\boldsymbol{\mu} \ast K)_{i,j,k} = -2 \Delta \mu_v.
\end{equation}
The discrete kernel $K \in \mathbb{R}^{3 \times 3 \times 3}$ is defined across its three depth slices as:
\begin{equation*}
    K_{\cdot,\cdot,1} = \begin{bmatrix} 0 & 0 & 0 \\ 0 & 1 & 0 \\ 0 & 0 & 0 \end{bmatrix}, \quad
    K_{\cdot,\cdot,2} = \begin{bmatrix} 0 & 1 & 0 \\ 1 & -6 & 1 \\ 0 & 1 & 0 \end{bmatrix}, \quad
    K_{\cdot,\cdot,3} = \begin{bmatrix} 0 & 0 & 0 \\ 0 & 1 & 0 \\ 0 & 0 & 0 \end{bmatrix}.
\end{equation*}

By expressing the convolution with the Laplacian as Toeplitz matrix-vector multiplication to derive the Lipschitz gradient for the calculation of the stepping size $\alpha$ as $L=\lVert \mathbf{A}\rVert_2^2 + 12\lambda$, absorbing the constant factor of $2$ into the regularization parameter $\lambda$, and noting that gradient descent subtracts the negative gradient, the resulting voxel-wise update rule is:
\begin{align}
    \mu_v^{(t)} &= \max \left( \mu_v^{(t-1)} + \alpha_v \left( \sum_{p} r_p \cdot A_{pv} + \lambda \left( \boldsymbol{\mu}^{(t-1)} \ast K \right)_v \right), 0 \right) \\
    &= \max \left( \mu_v^{(t-1)} + \alpha_v \left( \mathbf{VL}_{:,v} \cdot \mathbf{r}_{\mathbf{\mathbf{VP}_{:,v}}} + \lambda \left( \boldsymbol{\mu}^{(t-1)} \ast K \right)_v \right), 0 \right).
    \label{eq_1t}
\end{align}

\paragraph{\textbf{Total Variation (TV) Regularizer}}

Total variation (TV) regularization is designed to preserve edges while suppressing noise, making it particularly appropriate for vascular reconstruction \cite{Lou2010}. 
\begin{equation}
    TV(\boldsymbol{\mu}) = \sum_{v} \lVert \nabla \mu_v \rVert_2 = 
    \sum_{v} \sqrt{(\nabla_x \mu_v)^2 + (\nabla_y \mu_v)^2 + (\nabla_z \mu_v)^2  +\epsilon^2},
\end{equation}
where $\epsilon$ is a small positive constant.
The gradient of the Total Variation functional with respect to a single voxel $\mu_v$ is given by the negative divergence of the normalized gradient:
\begin{equation}
    \frac{\partial \text{TV}}{\partial \mu_v} = - \nabla \cdot \left( \frac{\nabla \mu_v}{\sqrt{\lVert \nabla \mu_v \rVert^2 + \epsilon^2}} \right).
\end{equation}

The Lipschitz gradient for this regularizer depends on $\mu$ and cannot be calculated explicitly. It can be bounded by a more conservative $\lVert \mathbf{A} \rVert_2^2 +\frac{12\lambda}{\epsilon}.$ However for small values of $\epsilon$, Lipschitz gradient becomes large, thus experimentation with values of $\alpha$ exceeding the bound is advisable.

Subtracting this derivative in the gradient descent optimization yields the voxel-wise update rule:
\begin{align}
    \mu_v^{(t)} &= \max \left( \mu_v^{(t-1)} + \alpha \left( \sum_{p} r_p \cdot A_{pv} + \lambda \nabla \cdot \left( \frac{\nabla \mu_v^{(t-1)}}{\sqrt{\lVert \nabla \mu_v^{(t-1)} \rVert^2 + \varepsilon^2}} \right) \right), 0 \right) \\
    &= \max \left( \mu_v^{(t-1)} + \alpha \left( \mathbf{VL}_{:,v} \cdot \mathbf{r_{VP_{:,v}}} + \lambda \nabla \cdot \left( \frac{\nabla \mu_v^{(t-1)}}{\sqrt{\lVert \nabla \mu_v^{(t-1)} \rVert^2 + \varepsilon^2}} \right) \right), 0 \right) 
    \label{eq_tv}
\end{align}

\begin{table}[htbp]
\centering
\caption{Regularizitaion term and corresponding Lipschitz gradient bound for step size $\alpha$}
\begin{tabular}{lll}
    \hline
    \textbf{Regularizer} & \textbf{Regularizer} & \textbf{Lipschitz gradient} \\
    \hline
    L1 & $-\lambda$ & $\lambda \cdot \operatorname{sgn}(\mu_v)$ \\
    $0^{th}-$order Tikhonov & $-\lambda{\mu}_v$ & $\lVert \mathbf{A}\rVert_2^2 + \lambda$ \\
    $1^{st}-$order Tikhonov & $\lambda(\mathbf{\mu} \ast K_{\Lambda}^{3D})_v$ & $\lVert \mathbf{A} \rVert_2^2 +12\lambda$ \\
    Total Variation & $\left[\text{div}\left(\frac{\nabla \mathbf{\mu}}{\sqrt{|\nabla \mathbf{\mu}|^2 + \epsilon^2}}\right)\right]_v$ & $\lVert \mathbf{A} \rVert_2^2 +\frac{12\lambda}{\epsilon}$ \\
    \label{tab:reg}
\end{tabular}
\end{table}

\begin{algorithm}
\caption{Gradient Descent Update with L1 Regularization}
\label{alg:gd_l1}
\begin{algorithmic}[1]
\Require Projection data $\mathbf{Y}$, regularization parameter $\lambda$, current volume $\boldsymbol{\mu}$, matrices $\mathbf{PV}^{(j)}, \mathbf{PL}^{(j)}, \mathbf{VL}^{(j)}$ for all views $j$
\For{$e = 1$ \textbf{to} $N_{epoch}$}
    \For{$j = 1$ \textbf{to} $N_{J}$}
        \State $\Delta\mathbf{\mu} \gets \mathbf{0} \in \mathbb{R}^{N_v \times 1}$ \Comment{Initialize gradient accumulator}
        \State $\mathbf{M} \gets \text{mprecon}(PV^{(j)}, PL^{(j)}, N_v)$ \Comment{Compute diagonal preconditioner}
        \State $N_{p} \gets$ number of pixels in the detector
        
        \For{$p = 1$ \textbf{to} $N_{p}$}
            \State $\mathbf{idx} \gets PV_{:,p}^{(j)} > 0$ \Comment{Boolean mask for valid voxels}
            \If{$\sum \mathbf{idx} > 0$}
                \State $\mathbf{v} \gets \mathbf{PV}_{\mathbf{idx}, p}^{(j)}$ \Comment{Indices of voxels traversed by ray $p$}
                \State $\mathbf{Q} \gets \mathbf{PL}_{\mathbf{idx}, p}^{(j)}$ \Comment{Segment lengths inside voxels}
                
                \State $\hat{Y}_p \gets \sum_k \mathbf{\mu}_{\mathbf{v}(k)} Q_k =\mu_{\mathbf{v}} \cdot\mathbf{Q}$ \Comment{Forward projection for ray $p$}
                \State $r_p \gets  Y_p^{(j)} - \hat{Y}_p$ \Comment{Calculate residual (y\_hat - y)}
                \State $d\mathbf{v} \gets r_p \cdot \mathbf{Q}$ \Comment{Gradient step for this ray}
                \State $\Delta\mathbf{\mu}_{\mathbf{v}} \gets \Delta\mathbf{\mu}_{\mathbf{v}} + \left( d\mathbf{v} \oslash \mathbf{M}_{\mathbf{v}} \right)$ \Comment{Accumulate preconditioned gradient}
            \EndIf
        \EndFor
        
        \State $\Delta\boldsymbol{\mu} \gets \Delta\boldsymbol{\mu} - \lambda \cdot \mathbf{M}$ \Comment{Add global L1 penalty}
        \State $\boldsymbol{\mu} \gets \boldsymbol{\mu} + \Delta\boldsymbol{\mu}$ \Comment{Apply gradient descent step}
        \State $\boldsymbol{\mu} \gets \max(\boldsymbol{\mu}, 0)$ \Comment{Enforce non-negativity constraint}
    \EndFor
\EndFor
\State \Return $\mathbf{\mu}$
\end{algorithmic}
\end{algorithm}

%% file: nesterov.tex
\paragraph{\textbf{Nesterov acceleration with diagonal majorization}}

To accelerate the convergence of gradient descent (GD), we investigate Nesterov's momentum scheme, which uses information from two consecutive volume estimates to construct an extrapolated solution \cite{BeckTeboulle2009}. Let
$\boldsymbol{\mu}^{(t)}$ and $\boldsymbol{\mu}^{(t-1)}$
denote the two most recent reconstructions. The extrapolated volume is defined as
\begin{equation}
\mathbf{z}^{(t)}
=
\boldsymbol{\mu}^{(t)}
+
\beta_t
\left(
\boldsymbol{\mu}^{(t)}
-
\boldsymbol{\mu}^{(t-1)}
\right),
\label{eq:nesterov_extrapolation}
\end{equation}
where $\beta_t$ is the momentum coefficient.

Unlike conventional GD, the projection residual and the corresponding correction are evaluated at the extrapolated volume $\mathbf{z}^{(t)}$, rather than at $\boldsymbol{\mu}^{(t)}$:
\begin{equation}
\mathbf{r}_z^{(t)}
=
\mathbf{Y}-\mathbf{A}\mathbf{z}^{(t)}.
\end{equation}

Using the diagonal Hessian majorizer $\mathbf{M}=\operatorname{diag}(\mathbf{m})$, defined previously, the accelerated GD update  is expressed as
\begin{equation}
\boldsymbol{\mu}^{(t+1)}
=
\max\left\{
\mathbf{0},
\mathbf{z}^{(t)}
+
\left[
\mathbf{A}^T\mathbf{r}_z^{(t)}
\right]
\oslash\mathbf{m}
\right\},
\label{eq:nesterov_tv}
\end{equation}
where $\oslash$ denotes element-wise division.
The non-negativity constraint is subsequently enforced.

For a constant momentum coefficient $\beta_t=\beta$,
setting $\beta=0$ recovers conventional GD, whereas
$\beta=1$ yields the extrapolation
$\mathbf{z}^{(t)}=
2\boldsymbol{\mu}^{(t)}-\boldsymbol{\mu}^{(t-1)}$.

